\documentclass[reqno]{amsart}

\usepackage[utf8]{inputenc} 
\usepackage[T1]{fontenc}    
\usepackage{url}            
\usepackage{booktabs}       
\usepackage{caption}
\usepackage{amsfonts}       
\usepackage{nicefrac}       
\usepackage{microtype}      
\usepackage[]{xcolor}         

\usepackage{amsmath,amssymb,amsfonts,mathrsfs}
\usepackage{nicefrac}
\usepackage{algorithmic}
\usepackage{algorithm}
\usepackage{graphicx}

\usepackage{multirow}
\usepackage{placeins}

\usepackage{wrapfig}

\usepackage{booktabs}

\usepackage{tikz}
\usetikzlibrary{positioning, arrows.meta, calc, matrix}

\usepackage{xcolor}

\definecolor{strblue}{HTML}{0F1ED2}
\definecolor{strred}{HTML}{E61E8C}

\usepackage{natbib}

\usepackage{comment}

\usepackage{aliascnt}
\usepackage{hyperref}
\usepackage{cleveref}

\theoremstyle{definition}

\newtheorem{theorem}{Theorem}[section]
\crefname{theorem}{Theorem}{Theorems}
\Crefname{theorem}{Theorem}{Theorems}

\newaliascnt{definition}{theorem}
\newtheorem{definition}[definition]{Definition}
\aliascntresetthe{definition}
\crefname{definition}{Definition}{Definitions}
\Crefname{definition}{Definition}{Definitions}

\newaliascnt{proposition}{theorem}
\newtheorem{proposition}[proposition]{Proposition}
\aliascntresetthe{proposition}
\crefname{proposition}{Proposition}{Propositions}
\Crefname{proposition}{Proposition}{Propositions}

\newaliascnt{lemma}{theorem}
\newtheorem{lemma}[lemma]{Lemma}
\aliascntresetthe{lemma}
\crefname{lemma}{Lemma}{Lemmas}
\Crefname{lemma}{Lemma}{Lemmas}

\newaliascnt{corollary}{theorem}

\aliascntresetthe{corollary}
\crefname{corollary}{Corollary}{Corollaries}
\Crefname{corollary}{Corollary}{Corollaries}

\newaliascnt{remark}{theorem}
\newtheorem{remark}[remark]{Remark}
\aliascntresetthe{remark}
\crefname{remark}{Remark}{Remarks}
\Crefname{remark}{Remark}{Remarks}

\newaliascnt{example}{theorem}
\newtheorem{example}[example]{Example}
\aliascntresetthe{example}
\crefname{example}{Example}{Examples}
\Crefname{example}{Example}{Examples}

\newaliascnt{assumption}{theorem}
\newtheorem{assumption}[assumption]{Assumption}
\aliascntresetthe{assumption}
\crefname{assumption}{Assumption}{Assumptions}
\Crefname{assumption}{Assumption}{Assumptions}

\numberwithin{equation}{section}

\crefname{equation}{Equation}{Equations}
\crefname{figure}{Figure}{Figures}
\crefname{table}{Table}{Tables}
\crefname{algorithm}{Algorithm}{Algorithms}

\crefname{section}{}{}
\creflabelformat{section}{\S #2#1#3}
\crefname{subsection}{}{}
\creflabelformat{subsection}{\S #2#1#3}
\crefname{appendix}{Appendix}{Appendixes}

\DeclareMathOperator*{\argmin}{arg\,min} 
\DeclareMathOperator*{\argmax}{arg\,max} 

\DeclareMathOperator{\Proj}{Proj}

\usepackage{pifont}

\title[Inverse Mixed-Integer Programming]{Inverse Mixed-Integer Programming:\\ Learning Constraints then Objective Functions}

\author{%
  Akira Kitaoka	
}
\address{NEC Corporation, 1753 Shimonumabe, Nakahara-ku, Kawasaki, Kanagawa, Japan }
\email{akira-kitaoka@nec.com}


\begin{document}

\maketitle

\begin{abstract}%
    Data-driven inverse optimization for mixed-integer linear programs (MILPs), which seeks to learn an objective function and constraints consistent with observed decisions, is important for building accurate mathematical models in a variety of domains, including power systems and scheduling. However, to the best of our knowledge, existing data-driven inverse optimization methods primarily focus on learning objective functions under known constraints, and learning both objective functions and constraints from data for MILPs remains largely unexplored.
    In this paper, we propose a two-stage approach for a class of inverse optimization problems in which the objective is a linear combination of given feature functions and the constraints are parameterized by unknown functions and thresholds. Our method first learns the constraints and then, conditioned on the learned constraints, estimates the objective-function weights.
    On the theoretical side, we provide finite-sample guarantees for solving the proposed inverse optimization problem. To this end, we develop statistical learning tools for pseudo-metric spaces under sub-Gaussian assumptions and use them to derive a learning-theoretic framework for inverse optimization with both unknown objectives and constraints. On the experimental side, we demonstrate that our method successfully solves inverse optimization problems on scheduling instances formulated as ILPs with up to 100 decision variables.
\end{abstract}

\section{Introduction}

A mixed-integer linear program (MILP) is a standard modeling framework for decision-making problems arising in scheduling \citep{suzuki2019TV,kolb2017learning,kumar2019automating,suenaga2024attempt}, power systems \citep{birge2017inverse}, and transportation planning \citep{bertsimas2015data}, among others.
In practice, however, the quantities required to build such mathematical models---including objective-function weights and constraint thresholds (e.g., processing times, due dates, and resource limits)---are not always known a priori.
This motivates the study of \emph{data-driven inverse optimization problems} (DDIOPs), which aim to identify objective and/or constraint parameters from observed data \citep[cf.][]{ahuja2001inverse,heuberger2004inverse,chan2019inverse,chan2023inverse}.
In this paper, we consider a setting in which a state $s$ (environment/input) and an expert's optimal solution $\hat{x}^*(s)$ under that state are observed, and we seek to learn the parameters of the objective function and constraints that explain the observations.
A schematic illustration of the input/output of data-driven inverse optimization is provided in \Cref{fig:IOP_image}.

\begin{figure}[ht]
\centering

\scalebox{0.8}{%
\newlength{\InputBoxWidth}
\setlength{\InputBoxWidth}{4.8cm}
\addtolength{\InputBoxWidth}{2em}
\newlength{\OutputBoxWidth}
\setlength{\OutputBoxWidth}{\dimexpr3\InputBoxWidth/4\relax}

\begin{tikzpicture}[
  font=\small,
  box/.style={
    draw, rounded corners=3pt, align=left,
    inner sep=6pt, line width=0.8pt
  },
  iobox/.style={
    draw, rounded corners=3pt, align=center,
    inner xsep=4pt, inner ysep=2pt,
    line width=0.8pt,
    text width=2.4cm,
    minimum height=1.4cm
  },
  blueboxin/.style={box, draw=strblue, text width=\InputBoxWidth},
  blueboxout/.style={box, draw=strblue, text width=\OutputBoxWidth},
  arrow/.style={-{Stealth[length=2.6mm]}, line width=1.1pt, draw=strblue},
]

\node[blueboxin, anchor=north west] (in1) at (0,0) {%
\textbf{Input 1: Model information}\\
$\bullet$ Objective-function coefficients $\theta$: \\
\hspace{1.4em}\textcolor{strred}{unknown}\\
$\bullet$ Constraint parameters $\phi$: \\
\hspace{1.4em}known: release times $r_j$\\
\hspace{1.4em}\textcolor{strred}{unknown}: processing times $p_j$,\\
\hspace{5.6em} due dates $d_j$
\hfill (etc.)
};

\node[blueboxin, anchor=north west] (in2) at ($(in1.south west)+(0,-4mm)$) {%
\textbf{Input 2: Observed schedule}\\
$\bullet$ Schedule $\hat{x}$
};

\coordinate (inmid-old) at ($0.5*(in1.east)+0.5*(in2.east)$);
\coordinate (inmid-new) at ($(in1.north east)!0.5!(in2.south east)$);
\coordinate (inmid) at ($(inmid-old)!0.5!(inmid-new)$);

\node[iobox, anchor=west] (io) at ($(inmid)+(7mm,-2.75mm)$) {%
\textbf{Inverse Optimization}
};

\node[blueboxout, anchor=west] (out) at ($(io.east)+(7mm,0)$) {%
\textbf{Output: Estimates}\\
$\bullet$ Objective-function \\ coefficients $\theta^{\sup}$\\
$\bullet$ Constraint thresholds \\(parameters) $\phi^{\sup}$
};

\draw[arrow] (in1.east) -| (io.north);
\draw[arrow] (in2.east) -| (io.south);
\draw[arrow] (io.east) -- (out.west);

\end{tikzpicture}%
} 

\caption{Inputs and outputs of the DDIOP (schematic). Unknown parameters (objective-function coefficients and constraint parameters) are estimated so that the observed schedule $\hat{x}$ becomes an optimal solution.}
\label{fig:IOP_image}
\vspace{-\intextsep}
\end{figure}

Known methods on inverse optimization provide only limited approaches for learning \emph{both} objective and constraint parameters from observed data.
In particular, for convex programs (including linear programs), there exist works that learn both the objective and the constraints \citep{aswani2018inverse}, and for linear programs (LPs), there exist works that learn the constraints (the feasible region) \citep[e.g.][]{chan2020inverse,ghobadi2020inferring,ren2025inverse}.
However, for DDIOPs for MILPs, learning both the objective and the constraints is computationally challenging, and many existing methods are limited to learning either only the objective or only the constraints.
Moreover, no generalization error analysis (learning theory) has been given for inverse optimization that learns both the constraints and the objective in MILPs, and it is still an important challenge to develop methods that solve DDIOPs for MILPs in a practical and theoretically justified manner.

In this paper, we propose a two-stage approach for an important problem class that includes DDIOPs for MILPs, where we \textbf{learn the constraints first and then learn the objective}.
Specifically, we consider settings in which the objective function is given by a linear combination of features, $\theta^\top f(x,s)$, and the constraints are expressed using componentwise threshold parameters $\phi$ as $h(x,s) + \phi \le 0$ (a representative instance of the more general constraint parameterization introduced in this paper).
Under this structure, we can efficiently construct a parameter $\phi^{\sup}$ (the componentwise largest parameter among those keeping the observed solutions feasible) that makes the constraints as tight as possible while keeping the observed solution $\hat{x}^*(s)$ feasible.
We then fix the resulting $\phi^{\sup}$ and learn the objective weights $\theta$ by minimizing the suboptimality loss \citep{ren2025inverse}.
For learning the objective function, we can incorporate existing algorithms for suboptimality-loss minimization \citep[e.g.][]{barmann2017emulating,barmann2018online,Kitaoka-2023-convergence-IRL,kitaoka2024fast}.

On the theory side, we first show that under conditions such as a finite state set, the proposed two-stage method can solve the DDIOP exactly, i.e., it can reproduce the observed optimal solutions.
Furthermore, we observe that the natural ``distance'' arising in generalization analysis for inverse optimization is in general a pseudo-metric rather than a metric.
We therefore develop sub-Gaussian learning theory in the setting of pseudo-metric spaces.
Using this result, we bound the generalization error of the solution obtained by minimizing the suboptimality loss, providing a learning-theoretic foundation for inverse optimization.

In our experiments, we apply the proposed method to two scheduling problems formulated as integer linear programs (ILP), and demonstrate that the DDIOP is solved in the sense that the observed solutions are reproduced as optimal solutions of the forward problem after learning.
To the best of our knowledge, there is no existing method that, under the same assumptions and within our setting, learns both constraint and objective parameters for MILPs.
Accordingly, instead of baseline comparisons, we evaluate performance using the mean squared error (MSE) of features and the suboptimality loss.

\paragraph{Contributions}
(i) We formulate a class of DDIOPs for MILPs in which the objective is given as a linear combination of features and the constraints admit an order-consistent (i.e., monotone) parameterization.
(ii) We propose a two-stage algorithm consisting of constraint learning and objective learning via suboptimality-loss minimization. Each stage employs existing components (constraint templates and suboptimality-loss minimization); the novelty lies in the coupling of the two stages via the tightest feasibility-preserving parameter $\phi^{\sup}$ and in its theoretical guarantee, stated in (iii), that this two-stage separation still exactly reproduces the observed optimal solutions.
(iii) Under finite distributions, we guarantee that the proposed method can exactly reproduce the observed optimal solutions (imitation).
(iv) We develop sub-Gaussian generalization analysis in pseudo-metric spaces with explicit constants (\Cref{sec:apx_learning_theory}), and leveraging it, derive generalization error bounds for inverse optimization (suboptimality-loss minimization) (\Cref{theo:subpoptimality_learning_E_pl_short,theo:subpoptimality_learning_prob_pl_short}).
(v) We empirically demonstrate exact reproduction of observed solutions as optima in scheduling ILPs.

We summarize the comparison between existing approaches and our method in \Cref{tab:compare_inverse_optimization}.

\begin{table*}[th]
    \vspace{-\intextsep}
  \caption{Summary comparison of inverse optimization methods with a guarantee that the parameters can be learned (C: capable of learning constraints; O: capable of learning objectives).}
  \label{tab:compare_inverse_optimization}
  \centering
  \setlength{\tabcolsep}{4pt}
  \begin{tabular}{p{2cm}ccp{0.70\textwidth}}
    \toprule
    Forward problem & C & O & References \\
    \midrule
    MILP & \checkmark & \checkmark &
    \textcolor{strblue}{\textbf{Ours}} \\
    MILP & \checkmark & $\times$ &
    \citet{kolb2017learning}, \citet{kumar2019automating}, \citet{suenaga2024attempt} \\
    MILP & $\times$ & \checkmark &
    \citet{barmann2017emulating,barmann2018online},
    \citet{gollapudi2021contextual},
    \citet{besbes2021online,besbes2025contextual},
    \citet{Kitaoka-2023-convergence-IRL,kitaoka2024fast},
    \citet{zattoni2024learning},
    \citet{sakaue2025online}, \citet{elhedhli2025inverse}, \citet{bulut2021complexity}, \citet{Kitaoka-2023-imitation-WIRL} \\
    ILP & $\times$ & \checkmark &
    \citet{suzuki2019TV} \\
    LP & \checkmark & \checkmark &
    \citet{aswani2018inverse} \\
    LP & \checkmark & $\times$ &
    \citet{chan2020inverse},
    \citet{ghobadi2020inferring},
    \citet{ren2025inverse} \\
    \bottomrule
  \end{tabular}
  \vspace{-\intextsep}
\end{table*}

\section{Related Work}

\Cref{sec:related_works_appendix} contains a companion related-work section, which surveys the statistical learning theory underlying our generalization error analysis and the literature on the characterization of lattice homomorphisms.

\paragraph{Inverse optimization algorithms}

In addition to the overview in the introduction, the method of \citet{aswani2018inverse} is computationally intractable, \citet{chan2020inverse,ghobadi2020inferring} rely on LP duality and hence do not directly apply to MILPs, and \citet{bulut2021complexity,elhedhli2025inverse} learn only the objective of MILPs. Neural ILP layers \citep{paulus2021comboptnet,nandwani2022solver} learn both without learning guarantees and are not used as baselines (\Cref{sec:related_works_appendix}).%

\paragraph{Loss functions and constraint learning}
The suboptimality loss was introduced by \citet{mohajerin2018data} and extended by \citet{ren2025inverse} to a form that also enables learning constraints by adding a penalty for constraint violation; our loss $\ell^{\mathrm{sub},\lambda}$ belongs to this line of work.
\citet{ren2025inverse} learn constraints under a known objective, but do not learn the objective and provide no generalization analysis. Our constraint learning uses constraint templates \citep{kolb2017learning,kumar2019automating,suenaga2024attempt}, which can handle integer and mixed-integer constraints.

\paragraph{Objective learning and generalization analysis}
Inverse optimization methods for learning objectives of MILPs include suboptimality-loss-based methods in the offline setting \citep{suzuki2019TV,Kitaoka-2023-convergence-IRL,Kitaoka-2023-imitation-WIRL,kitaoka2024fast,zattoni2024learning} and the online setting \citep{barmann2017emulating,barmann2018online,besbes2021online,besbes2025contextual,gollapudi2021contextual,sakaue2025online}; \citet{kitaoka2024fast,kitaoka2026explicit} gave explicit iteration-count bounds, but both assume known constraints.
Existing generalization analyses for inverse optimization and predict-then-optimize \citep{elbalghiti2019generalization,li2025fenchel,sun2023maximum,mohajerin2018data,aswani2018inverse} either target learning only the objective or remain at asymptotic consistency; no finite-sample generalization analysis has been given for inverse optimization that learns \emph{both} the constraints and the objective in MILPs, which this work provides.

\section{Background}
\label{sec:background}

In this section, we provide the necessary background to introduce our proposed method. The probability simplex is defined as
$
    \Delta^{D-1} := \left\{ \theta \in \mathbb{R}_{\geq 0}^{D} \middle| \sum_{i=1}^D \theta _i =1 \right\}
$,
where $\mathbb{R}_{\geq 0} := [0 , \infty)$.

\subsection{Inverse Optimization}
\label{sec:inverse_optimization}

Let $\mathcal{X} \subset \mathbb{R}^{D_{\mathcal{X}}}$ be a nonempty subset, and let $\mathcal{S}$ be a nonempty set. Let $\widetilde{\Theta}$ denote a parameter space. Let
$\tilde{f} \colon \mathcal{X} \times \widetilde{\Theta} \times \mathcal{S} \to \mathbb{R}$
and, for $j=1,\ldots,J$, let
$g_j \colon \mathcal{X} \times \widetilde{\Theta} \times \mathcal{S} \to \mathbb{R}$.
Define $g := (g_1,\ldots,g_J)$.
For a given $s \in \mathcal{S}$ and parameter $\widetilde{\theta} \in \widetilde{\Theta}$, we define an optimal solution of the forward problem as follows:
\begin{equation}
    \displaystyle
    x^{*} (\widetilde{\theta} , s)\in 
    \mathbf{FOP} \left( \widetilde{\theta}, s \right)
    =
    \argmax_{x \in \mathcal{X}} 
    \left\{ 
        \tilde{f} (x, \widetilde{\theta}, s)
    \, \middle| \,
        g (x, \widetilde{\theta}, s) \leq 0
    \right\}
    .
    \label{eq:FOP}
\end{equation}

Assume that there exists an (unknown) parameter $\widetilde{\theta}^{\mathrm{true}}$ and that we are given observed $\mathcal{X}$-valued data over $\mathcal{S}$ of the form $\hat{x}^{*} (s) = x^{*}(\widetilde{\theta}^{\mathrm{true}}, s)$. The \textit{data-driven inverse optimization problem} (DDIOP) associated with the forward problem \Cref{eq:FOP} is the problem of finding a parameter $\widetilde{\theta} \in \widetilde{\Theta}$ such that, for any $s \in \mathcal{S}$,
\begin{equation}
    \hat{x}^{*} (s) \in \mathbf{FOP} (\widetilde{\theta} , s)
    .
    \label{eq:IOP}
\end{equation}

For each $\widetilde{\theta} \in \widetilde{\Theta}$ and $s \in \mathcal{S}$, define the feasible region by
\[
    \mathcal{X} (\widetilde{\theta} , s)
    := \left\{ 
        x \in \mathcal{X} 
    \, \middle| \, 
        g (x , \widetilde{\theta} , s ) \leq 0
    \right\}.
\]
Hereafter, we omit arguments that are not used, whenever this causes no confusion.

For $u \in \mathbb{R}$, define the ReLU function by $\mathrm{ReLU}(u) := \max(u,0)$. Let $\lambda \in \mathbb{R}_{\geq 0}$ be a constant. As a criterion for DDIOP, we consider the suboptimality loss $\ell^{\mathrm{sub}, \lambda} \colon \mathcal{X} \times \widetilde{\Theta} \times \mathcal{S} \to \mathbb{R}_{\geq 0}$ \citep[cf.][]{ren2025inverse}, defined by
\begin{align*}
        & \ell^{\mathrm{sub}, \lambda} \left( x, \widetilde{\theta} ,s \right) \\
        & := \mathrm{ReLU} \left( \max_{x^{\star} \in \mathcal{X} ( \widetilde{\theta} , s )} \tilde{f} ( x^{\star} , \widetilde{\theta} , s) - \tilde{f} (x, \widetilde{\theta},  s ) \right) + \lambda \sum_{j=1}^J \mathrm{ReLU} \left( g_j ( x , \widetilde{\theta} , s) \right).
    \end{align*}
The suboptimality loss satisfies the following property.

\begin{proposition}[Cf. {\citet[Proposition~2.1]{ren2025inverse}}]
    \label{prop:suboptimality_imitation}
    Let $\lambda > 0$ be a constant, and let $x \in \mathcal{X}$, $\widetilde{\theta} \in \widetilde{\Theta}$, and $s \in \mathcal{S}$. Then,
    $x \in \mathbf{FOP}(\widetilde{\theta}, s)$ holds if and only if $\ell^{\mathrm{sub}, \lambda}(x,\widetilde{\theta},s)=0$.
\end{proposition}

The above proposition reduces to \citet[Proposition~2.1]{ren2025inverse} when $J=1$, and it can be proved in the same manner as in \citet[Proposition~2.1]{ren2025inverse}.

Let a probability distribution $\mathbb{P}_{\mathcal{S}}$ on $\mathcal{S}$ be given, and let $S$ be an $\mathcal{S}$-valued random variable distributed according to $\mathbb{P}_{\mathcal{S}}$. As a formulation of DDIOP, we consider the following optimization problem, which seeks $\widetilde{\theta} \in \widetilde{\Theta}$ satisfying
\begin{equation}
    \min_{\widetilde{\theta} \in \widetilde{\Theta} }
    \mathbb{E} \ell^{\mathrm{sub}, \lambda} ( \hat{x}^{*} (S), \widetilde{\theta} ,S ) 
    \label{eq:IOP_E_suboptimality_loss_lambda}
\end{equation}
By \Cref{prop:suboptimality_imitation}, when $\lambda > 0$, any parameter $\widetilde{\theta}^* \in \widetilde{\Theta}$ for which \Cref{eq:IOP_E_suboptimality_loss_lambda} attains the value $0$ satisfies \Cref{eq:IOP} for almost every $s \in \mathcal{S}$.

\subsection{Problem Setting}

In this paper, we study a fundamental and important class of forward problems---including mixed-integer linear programs (MILPs)---whose objective functions can be expressed as linear combinations of piecewise-linear functions. As a setting for DDIOP for MILP, we suppose the following assumption.

\begin{assumption}
    \label{assu:WIRL}
    Let the weight set $\Theta \subset \mathbb{R}^{D}$ be a bounded, closed, and convex set.
    Let $\Phi$ be a nonempty set that serves as the parameter space for the constraints.
    For $i=1,\ldots,D$, let $f_i \colon \mathcal{X} \times \mathcal{S} \to \mathbb{R}$, and define $f := (f_1,\ldots,f_D)$.
    For $j=1,\ldots,J$, let $g_j \colon \mathcal{X} \times \Phi \times \mathcal{S} \to \mathbb{R}$.
    For any $s\in\mathcal{S}$ and $\phi\in\Phi$, each component of the maps $f(\bullet,s)\colon \mathbb{R}^{D_{\mathcal{X}}} \to \mathbb{R}^D$ and $g(\bullet,\phi,s)\colon \mathbb{R}^{D_{\mathcal{X}}} \to \mathbb{R}^J$ is piecewise linear.
    For any $s\in\mathcal{S}$, the map $f(\bullet,s)\colon \mathbb{R}^{D_{\mathcal{X}}} \to \mathbb{R}^D$ is $L_f$-Lipschitz.
    For any $\phi \in \Phi$ and $s \in \mathcal{S}$, the feasible region $\mathcal{X} (\phi , s)$ is compact. Moreover, for any $\phi \in \Phi$ and $s \in \mathcal{S}$, $\mathcal{X} (\phi , s)$ is a finite union of polytopes, i.e., the union of finitely many convex polytopes (convex hulls of finite sets) \citep{kitaoka2024fast}. For example, the set $P \cap (\mathbb{Z}^{k_1} \times \mathbb{R}^{k_2})$ of mixed-integer points of a convex polytope $P$ is a finite union of polytopes.
    The observed data $\hat{x}^{*} \colon \mathcal{S} \to \mathcal{X}$ are generated by unknown parameters $\theta^{\mathrm{true}} \in \Theta$ and $\phi^{\mathrm{true}} \in \Phi$ in the sense that, for any $s \in \mathcal{S}$, $\hat{x}^{*}(s) = x^*(\theta^{\mathrm{true}}, \phi^{\mathrm{true}}, s)$.
    Moreover, $\sup_{s \in \mathcal{S}} ( \max_{x \in \mathcal{X} (\phi^{\mathrm{true}} , s)} f_i (x , s) - \min_{x \in \mathcal{X} (\phi^{\mathrm{true}} , s)} f_i (x , s) ) < \infty$ for each $i = 1 , \ldots , D$.
\end{assumption}

Under \Cref{assu:WIRL}, for a given $s \in \mathcal{S}$ and parameters $\theta \in \Theta$ and $\phi \in \Phi$, the forward problem \Cref{eq:FOP} is given by
\begin{equation}
    \displaystyle
    x^{*} (\theta ,\phi, s)
    \in
    \mathbf{FOP} (\theta, \phi, s)
    =
    \argmax_{x \in \mathcal{X}} 
    \left\{ 
        \theta^{\top} f (x,s)
    \, \middle| \,
        g (x, \phi, s) \leq 0
    \right\}
    .
    \label{eq:FOP_linear}
\end{equation}
Moreover, given the observed data $\hat{x}^{*} \colon \mathcal{S} \to \mathcal{X}$, the DDIOP \Cref{eq:IOP} amounts to identifying the objective weight $\theta \in \Theta$ and the constraint parameter $\phi \in \Phi$ such that, for any $s \in \mathcal{S}$,
\begin{equation}
    \hat{x}^{*} (s) \in \mathbf{FOP} (\theta, \phi, s)
    .
    \label{eq:IOP_linear}
\end{equation}
In what follows, we write $\ell^{\mathrm{sub},\lambda} (x , \theta , \phi , s)$ for the suboptimality loss with $\widetilde{\theta} = (\theta , \phi)$ substituted, and abbreviate $\ell^{\mathrm{sub},\lambda} (\theta , \phi , s) := \ell^{\mathrm{sub},\lambda} (\hat{x}^{*} (s) , \theta , \phi , s)$.

As a technical condition, we further suppose the following assumption.

\begin{assumption}
    \label{assu:WIRL-uniqueness}
    For any state $s\in\mathcal{S}$, the feature vector $f\bigl(x^*(\theta^{\mathrm{true}},\phi^{\mathrm{true}},s), s\bigr)$ is uniquely determined.
    That is, $f (\mathbf{FOP} (\theta^{\mathrm{true}} , \phi^{\mathrm{true}} , s ) , s)$ is a singleton.
\end{assumption}

\Cref{assu:WIRL-uniqueness} is natural in view of \Cref{lem:Psi_set_is_almost_Phi} (\Cref{sec:apx_imitation_theory}): when the state set $\mathcal{S}$ is finite, for $\Theta = \Delta^{D-1}$ and fixed $\phi$, it holds for almost every $\theta \in \Delta^{D-1}$.

In addition, regarding the constraint structure, we assume the following.

\begin{assumption}
    \label{assu:IOP_std}
    $\Phi$ is a subset of a Euclidean space, and the triple $(\Phi,\land,\lor)$ is a lattice with respect to the componentwise $\min$ and $\max$ (\Cref{sec:lattice}).
    For any $x \in \mathcal{X}$ and $s \in \mathcal{S}$, the map
    $g(x,\bullet,s)\colon \Phi \to \mathbb{R}^J$
    is a lattice homomorphism.
    For any $s \in \mathcal{S}$, the set $\{ \phi \in \Phi \mid g ( \hat{x}^{*} (s) , \phi , s ) \leq 0 \}$ has a maximal element.
\end{assumption}
For the definition of a lattice, see \Cref{sec:lattice}.
We give concrete examples satisfying \Cref{assu:IOP_std}. The characterization of lattice homomorphisms used in the examples (\Cref{theo:characterization_lattice_homomorphism}) and its proof are given in \Cref{sec:characterization_lattice_homomorphism}.
\begin{example}
    \label{exa:lattice_homomorphism_affine}
    Let $\Phi = \{ 0 \}^{J^0} \times \Phi^+ \times \Phi^-$.
    Let $h^0 \colon \mathcal{X} \times \mathcal{S} \to \mathbb{R}^{J^0}$,
    $h^+ \colon \mathcal{X} \times \mathcal{S} \to \mathbb{R}^{J^+}$, and
    $h^- \colon \mathcal{X} \times \mathcal{S} \to \mathbb{R}^{J^-}$, and set $J = J^0 + J^+ + J^-$.
    Define $g \colon \mathcal{X} \times \Phi \times \mathcal{S} \to \mathbb{R}^J$ by, for $\phi=(0,\phi^+,\phi^-)\in \Phi$,
    \begin{align*}
        g (x , \phi , s )
        & = (g^0 (x , \phi , s ) , g^+ (x , \phi , s ), g^- (x , \phi , s ))
        \\
        & := (h^0 (x , s ) , h^+ (x , s ) + \phi^+ , h^- (x , s ) + \phi^- )
    \end{align*}
    Since each component of $g$ is monotone nondecreasing in a certain single coordinate of $\Phi$, it follows from \Cref{theo:characterization_lattice_homomorphism} that $g$ is a lattice homomorphism.
\end{example}

\begin{example}
    \label{exa:lattice_homomorphism_affine_check}
    Using $\Phi$ and $g$ from \Cref{exa:lattice_homomorphism_affine}, define
    $\check{\Phi}^- := (-\Phi^-)$,
    $\check{\Phi} := \{ 0 \}^{J^0} \times \Phi^+ \times \check{\Phi}^-$, and, for $x \in \mathcal{X}$, $\check{\phi}^- \in \check{\Phi}^-$, and $s \in \mathcal{S}$, define
    $\check{g}^- (x , \check{\phi}^- , s ) := - g^{-} (x , - \check{\phi}^- , s) = \check{\phi}^- - h^- (x , s)$,
    and set $\check{g} = (g^0 , g^+ , \check{g}^- )$.
    Since each component of $\check{g}$ is monotone nondecreasing in a certain single coordinate of $\check{\Phi}$, it follows from \Cref{theo:characterization_lattice_homomorphism} that $\check{g}$ is also a lattice homomorphism.
    Then, the constraint set described by the constraint map $\check{g}$ can be written as follows: for
    $\check{\phi} = (0 , \phi^+ , \check{\phi}^-) \in \check{\Phi}$ and $s \in \mathcal{S}$,
    \[
        \mathcal{X} (\check{\phi} , s) = \left\{ x \in \mathcal{X} \middle|
            h^0 (x , s ) \leq 0 ,
            h^+ (x , s ) \leq - \phi^+ ,
            h^- (x , s ) \geq \check{\phi}^-
        \right\}.
    \]
\end{example}

\subsection{Learning Constraints}

For a subset $\mathcal{S}^\prime \subset \mathcal{S}$, let the constraint-parameter estimate $\phi^{\sup}(\mathcal{S}^\prime)\in\Phi$ satisfy
\begin{equation}
    \phi^{\sup} (\mathcal{S}^\prime)
    \in 
    \argmax_{\phi \in \Phi} \left\{ \phi  \, \middle| \,  g \bigl( \hat{x}^{*} (s) ,\phi , s \bigr) \leq 0 \text{ for } s \in \mathcal{S}^\prime \right\}
    \notag
\end{equation}
where  $\max$ denotes maximization in the sense of $\leq$-maximality.
Under \Cref{assu:IOP_std}, since $g ( \hat{x}^{*} (s) , \bullet , s)$ is a lattice homomorphism, the set on the right-hand side is closed under $\lor$. Hence a maximal element, if it exists, is the maximum and is uniquely determined.
If it is clear from the context, we write $\phi^{\sup} = \phi^{\sup}(\mathcal{S})$.

\begin{proposition}
    \label{prop:constraint_learning_lattice_map}
    Suppose \Cref{assu:IOP_std}.
    Let $\mathcal{S}^\prime \subset \mathcal{S}$ be a nonempty finite set. Then, $\phi^{\sup} (\mathcal{S}^\prime)$ exists and
    $
        \phi^{\sup} (\mathcal{S}^\prime) 
        = \bigwedge_{s \in \mathcal{S}^\prime} \phi^{\sup} (\{ s \})
        .
    $
\end{proposition}

For the proof of \Cref{prop:constraint_learning_lattice_map} and an example, see \Cref{sec:apx_constraint_learning}.

\subsection{Learning Objective Function}

Fix $\phi \in \Phi$ such that $\hat{x}^{*} (s) \in \mathcal{X} (\phi , s)$ for any $s \in \mathcal{S}$. When the state set $\mathcal{S}$ is finite and $\mathbb{P}_{\mathcal{S}}$ is uniform, the suboptimality loss $\ell^{\mathrm{sub},0}$ is convex and Lipschitz continuous in $\theta$, with a subgradient at $\theta \in \Theta$ given by $F (\theta) := \mathbb{E} \left( f ( x^* (\theta  , \phi, S ) ,S)- f (\hat{x}^* (S) , S) \right)$ (\Cref{prop:SL_is_Lipschitz_convex}, \Cref{sec:apx_gamma_lowerbound}).

Let $\{ \alpha_k \}_k$ be a nonnegative real sequence, referred to as the learning rate. Let $\Proj_{\Theta}$ denote the projection onto the weight set $\Theta$.
As an inverse optimization method for estimating the objective function of the mixed-integer linear program \Cref{eq:FOP_linear}, one may employ the projected subgradient method for the suboptimality loss, shown in \Cref{alg:psgd_suboptimality_loss}.
The learning-rate sequence $\{ \alpha_k \}_k$ is said to be the square-root step size (SRSS) rule if there exists a constant $\alpha>0$ such that
$
    \alpha_k = \alpha k^{-\frac{1}{2}}
$.
The learning-rate sequence $\{ \alpha_k \}_k$ is said to be the square-root step length (SRSL) rule if there exists a constant $\alpha>0$ such that, when $F(\theta^k)\neq 0$,
$
    \alpha_k = \alpha k^{-\frac{1}{2}}\| F (\theta^k ) \|^{-1} 
$,
and otherwise $\alpha_k=0$.

\begin{algorithm}[t]
        \caption{Minimization of the suboptimality loss}\label{alg:psgd_suboptimality_loss}
        \begin{algorithmic}[1]
            \STATE Initialize $\theta^1 \in \Theta$
            \FOR{ $k = 1 , \ldots, K-1 $}
                \STATE For each $s \in \mathcal{S}$, solve $x^* (\theta^k ,\phi, s )$
                \STATE $\theta^{k+1} \leftarrow \Proj_{\Theta} \left( \theta^k - \alpha_{k} F ( \theta^k) \right)  $ 
            \ENDFOR
            \RETURN $ \displaystyle \theta_\mathrm{best}^{K} \in \argmin_{\theta \in \{ \theta^{k} \}_{k=1}^K } \mathbb{E}\ell^{\mathrm{sub},0} (\theta, \phi ,S)$ 
        \end{algorithmic}
    \end{algorithm}

Let $\mathrm{diam}(\Theta) := \sup_{\theta , \theta^\prime \in \Theta} \| \theta - \theta^\prime \|$.
With the SRSS or SRSL learning rate, \Cref{alg:psgd_suboptimality_loss} reaches the minimum value $0$ of the suboptimality loss in \emph{finitely} many iterations, and the number of iterations admits an explicit upper bound of order $O ( \gamma(\ell^{\mathrm{sub}})^{-2} )$, where $\gamma(\ell^{\mathrm{sub}}) > 0$ is a geometric constant \citep[Corollaries~4.6--4.9]{kitaoka2024fast}; a lower bound on $\gamma(\ell^{\mathrm{sub}})$ for ILPs is given by \citet{kitaoka2026explicit}.
The precise statements (\Cref{def:gamma_L_sub,prop:intention_learning_complete_Polyhedral,prop:iteration_bound_objective,prop:gamma_lowerbound_ILP_simple,prop:gamma_lowerbound_ILP}) are collected in \Cref{sec:apx_gamma_lowerbound}.

\section{Proposed method}
\label{sec:proposed_method}

\begin{algorithm}[t]
        \caption{Maximizing the feasible set and then minimizing the suboptimality loss}\label{alg:Solve_IOP}
        \begin{algorithmic}[1] 
            \STATE Set $\varepsilon \geq 0$
            \STATE Compute $\phi^{\sup} (\mathcal{S})$
            \STATE Compute $\theta^{\sup}$ satisfying\\ $
    \mathbb{E}_{S} \ell^{\mathrm{sub}, 0} ( \hat{x}^{*} (S), \theta^{\sup} , \phi^{\sup} ,S ) \leq \varepsilon$
            \STATE Return $\theta^{\sup} \in \Theta, \phi^{\sup} \in \Phi$
        \end{algorithmic}
    \end{algorithm}
Under \Cref{assu:WIRL,assu:WIRL-uniqueness,assu:IOP_std}, we define in \Cref{alg:Solve_IOP} an algorithm for solving \Cref{eq:IOP_linear}.
For example, line~2 can be implemented by \Cref{prop:constraint_learning_lattice_map}, which corresponds to constraint learning studied in, e.g., \citet{kolb2017learning,kumar2019automating,suenaga2024attempt}; line~3 can be implemented by \Cref{alg:psgd_suboptimality_loss}.
Setting $\varepsilon = 0$ in line~1 requires the observed solutions to be reproduced exactly on the sample (\Cref{theo:alg:Solve_IOP_is_solve_IOP,theo:alg:Solve_IOP_with_min_suboptloss_is_solve_IOP_LP}); a positive $\varepsilon$ allows an approximate solution, e.g., when \Cref{alg:psgd_suboptimality_loss} is stopped by an iteration limit. Line~3 uses the loss with $\lambda = 0$: by the definition of $\phi^{\sup}$, the observed solutions are feasible on the sample, so the penalty term vanishes there. The penalty with $\lambda > 0$ measures constraint violations at unseen states and is used in the generalization analysis (\Cref{sec:statistical_learning_theory}) and in the held-out evaluation (\Cref{sec:heldout}).

\section{Theory of Imitativeness}
\label{sec:imitation_thery}

In this section, we show \Cref{alg:Solve_IOP} can solve \cref{eq:IOP_linear}.
A key reason why \Cref{alg:Solve_IOP}---that is, learning the constraints followed by learning the objective---can solve \Cref{eq:IOP_linear} is the validity of the following proposition.
\begin{proposition}
    \label{prop:true_para_include_sup_para}
    Suppose \Cref{assu:WIRL,assu:WIRL-uniqueness,assu:IOP_std}, and that $\phi^{\sup}$ exists (e.g., $\mathcal{S}$ is finite; \Cref{prop:constraint_learning_lattice_map}).
    Then, for $s \in \mathcal{S}$, if
    $\hat{x}^{*}(s) \in \mathbf{FOP}(\theta,\phi^{\mathrm{true}},s)$,
    then
    $\hat{x}^{*}(s) \in \mathbf{FOP}(\theta,\phi^{\sup},s)$.
    Moreover, for any $s \in \mathcal{S}$, the vector $f\bigl(x^*(\theta^{\mathrm{true}},\phi^{\sup},s),s\bigr)$ is uniquely determined.
\end{proposition}

By \Cref{prop:true_para_include_sup_para}, using the parameter $\phi^{\sup}$ obtained via constraint learning, the given observation $\hat{x}^*$ satisfies
$\hat{x}^*(s) \in \mathbf{FOP}(\theta^{\mathrm{true}},\phi^{\sup},s)$.
Combining \Cref{prop:intention_learning_complete_Polyhedral,prop:true_para_include_sup_para}, we conclude that, for MILPs, one can solve \Cref{eq:IOP_linear} by learning the constraints and subsequently learning the objective. In particular, the following theorem holds.

\begin{theorem}
    \label{theo:alg:Solve_IOP_is_solve_IOP}
    Suppose \Cref{assu:WIRL,assu:WIRL-uniqueness,assu:IOP_std}.
    Let the state set $\mathcal{S}$ be finite, and let $\mathbb{P}_{\mathcal{S}}$ be the uniform distribution.
    Let $\varepsilon=0$.
    Then, for a sufficiently large number of iterations $K$, the output $(\theta^{\sup},\phi^{\sup})$ of \Cref{alg:Solve_IOP}---where line~3 is implemented by \Cref{alg:psgd_suboptimality_loss} with the learning rate chosen as either SRSS or SRSL---satisfies
    $\hat{x}^{*}(s) \in \mathbf{FOP}(\theta^{\sup},\phi^{\sup},s)$, i.e., \Cref{eq:IOP_linear}.
\end{theorem}

From \Cref{lem:Psi_set_is_almost_Phi} and \Cref{theo:alg:Solve_IOP_is_solve_IOP}, we obtain the following corollary.
\begin{theorem}
    \label{theo:alg:Solve_IOP_with_min_suboptloss_is_solve_IOP_LP}
    Suppose \Cref{assu:WIRL,assu:IOP_std}.
    Let the state set $\mathcal{S}$ be finite, let $\mathbb{P}_{\mathcal{S}}$ be the uniform distribution, and let the weight space be $\Theta = \Delta^{D-1}$.
    Let $\varepsilon=0$.
    Then, for almost every $\theta^{\mathrm{true}} \in \Delta^{D-1}$, and for a sufficiently large number of iterations $K$, the output $(\theta^{\sup},\phi^{\sup})$ of \Cref{alg:Solve_IOP}---where line~3 is implemented by \Cref{alg:psgd_suboptimality_loss} with the learning rate chosen as either SRSS or SRSL---satisfies
    $\hat{x}^{*}(s) \in \mathbf{FOP}(\theta^{\sup},\phi^{\sup},s)$, i.e., \Cref{eq:IOP_linear}.
\end{theorem}

The proofs of \Cref{prop:true_para_include_sup_para} and \Cref{theo:alg:Solve_IOP_is_solve_IOP,theo:alg:Solve_IOP_with_min_suboptloss_is_solve_IOP_LP} are provided in \Cref{sec:apx_imitation_theory}.

Moreover, combining \Cref{theo:alg:Solve_IOP_with_min_suboptloss_is_solve_IOP_LP} with existing iteration-count bounds and lower bounds on $\gamma(\ell^{\mathrm{sub}})$, the number of iterations needed to solve \Cref{eq:IOP_linear} for ILP forward problems is bounded above by explicit functions of the sample size $N$, the feature dimension $D$ and range, and the learning-rate constant (\Cref{prop:iteration_count_solve_IOP}, \Cref{sec:apx_gamma_lowerbound}).

\section{Statistical Learning Theory}
\label{sec:statistical_learning_theory}

In this section, we evaluate the generalization error of \Cref{alg:Solve_IOP}.
Standard generalization bounds for losses that are Lipschitz continuous on the parameter space are stated for a metric, whose covering numbers control the supremum of the empirical process (cf.~\citealt[Theorem 5]{shalev2009stochastic}; \citealt[Problem 5.12]{van2014probability}; \citealt{liao2020notes}; \citealt[Theorem 8.3.23]{vershynin2020high}).
In inverse optimization, however, the quantity controlling the increments of the suboptimality loss (the right-hand side of \Cref{prop:suboptimality_Lipschitz_short}) is only a pseudo-metric: distinct constraint parameters defining the same feasible region are at distance $0$, and two weights are not separated at a state whose feasible region is the singleton $\{ \hat{x}^{*} (s) \}$.
We therefore extend the covering-number analysis under a sub-Gaussian assumption to pseudo-metric spaces (\Cref{theo:learning_theory_expectation,theo:learning_theory_prob}, \Cref{sec:apx_learning_theory}) and apply it to inverse optimization (\Cref{theo:subpoptimality_learning_E_pl_short,theo:subpoptimality_learning_prob_pl_short}).

Before presenting a generalization-error analysis of \Cref{alg:Solve_IOP} under \Cref{assu:WIRL,assu:WIRL-uniqueness,assu:IOP_std}, we examine the relationship between $\phi^{\sup}$ and $\phi^{\mathrm{true}}$.
\begin{proposition}
    \label{prop:true_constraints_para_dominated}
    Suppose \Cref{assu:WIRL,assu:IOP_std}. Then for any nonempty subset $\mathcal{S}^\prime \subset \mathcal{S}$ for which $\phi^{\sup} (\mathcal{S}^\prime)$ exists (in particular, for any nonempty finite subset by \Cref{prop:constraint_learning_lattice_map}), we have $\phi^{\sup} (\mathcal{S}^\prime) \geq \phi^{\mathrm{true}}$. In particular, if $\phi^{\sup}$ exists, then
    $\phi^{\sup} \geq \phi^{\mathrm{true}}$.
\end{proposition}
\Cref{prop:true_constraints_para_dominated} implies that the parameter $\phi^{\sup}$ obtained via constraint learning satisfies $\phi^{\sup} \geq \phi^{\mathrm{true}}$.
In what follows, let $\Phi \subset \mathbb{R}^J$.
For $\delta \in \mathbb{R}_{\geq 0}^J$, define
$
    \Phi (\delta) := \Phi \cap \prod_{j=1}^J [ \phi^{\mathrm{true}}_j , \phi^{\mathrm{true}}_j + \delta_j ]  $
.
For $\phi , \phi^\prime \in \Phi$, define an equivalence relation $\phi \sim \phi^\prime$ by requiring that, for $\mathbb{P}_{\mathcal{S}}$-a.e.\ $s \in \mathcal{S}$, $\mathcal{X} (\phi , s) = \mathcal{X} (\phi^\prime , s)$.
Let $[\phi]$ denote the equivalence class of $\phi \in \Phi$.
Moreover, we write $\mathbf{1}_{A}$ for the indicator function of an event $A$.
Then the suboptimality loss admits the following Lipschitz property.

\begin{proposition}[{See \Cref{prop:suboptimality_Lipschitz} for details.}]
    \label{prop:suboptimality_Lipschitz_short}
    Suppose the setting of \Cref{assu:WIRL,assu:IOP_std}.
    Assume that $\sup_{\theta \in \Theta} \| \theta \|_2 \leq 1$.
    Then, for any $s \in \mathcal{S}$, $\theta , \theta^\prime \in \Theta$, and any $\phi , \phi^\prime \geq \phi^{\mathrm{true}}$ satisfying $g (\hat{x}^{*} (s) , \phi , s) \leq 0$ and $g (\hat{x}^{*} (s) , \phi^\prime , s) \leq 0$, we have
    \begin{align*}
        & \left|
            \ell^{\mathrm{sub},\lambda} ( \theta , \phi , s )
            -
            \ell^{\mathrm{sub},\lambda} ( \theta^\prime , \phi^\prime , s )
        \right|
        \notag\\
        & \leq L_f d^H (\mathcal{X} (\phi , s ) , \mathcal{X} (\phi^\prime , s ) )
        +
        L_f d^H \left( \mathcal{X} (\phi^{\mathrm{true}} , s), \{ \hat{x}^{*} (s) \}\right)
        \| \theta - \theta^\prime \|
        .
    \end{align*}
    Here, $d^H$ denotes the Hausdorff distance (cf.\ \Cref{sec:Hausdorff_metric}).
\end{proposition}

When the constraint parameter output by \Cref{alg:Solve_IOP} is equivalent to $\phi^{\mathrm{true}}$ under $\sim$, the first term in \Cref{prop:suboptimality_Lipschitz_short} vanishes.
In what follows, we write $C_{\mathcal{X}} := \| d^H ( \mathcal{X} (\phi^{\mathrm{true}} , S), \{ \hat{x}^{*} (S) \} ) \|_{\psi_2}$ (here $\| \bullet \|_{\psi_2}$ denotes the sub-Gaussian norm; see \Cref{sec:apx_learning_theory}).
For an i.i.d.\ sample $(S^{(1)} , \ldots , S^{(N)})$, we also write $\mathcal{S}_N := \{ S^{(1)} , \ldots , S^{(N)} \}$.
Combining this observation with \Cref{theo:learning_theory_expectation,theo:learning_theory_prob} yields \Cref{theo:subpoptimality_learning_E_pl_short,theo:subpoptimality_learning_prob_pl_short}.

\begin{theorem}[{See \Cref{theo:subpoptimality_learning_E_pl} for details.}]
    \label{theo:subpoptimality_learning_E_pl_short}
    Let $\delta \in \mathbb{R}_{\geq 0}^J$.
    Moreover, assume that either (i) $\mathcal{S}$ is finite, or (ii) $\Theta = \Delta^{D-1}$, $\mathcal{X} \subset \mathbb{Z}^{D_{\mathcal{X}}}$, and $f (\mathcal{X} \times \mathcal{S}) \subset \mathbb{Z}^{D}$ (i.e., the forward problem is an ILP).
    Suppose the setting of \Cref{prop:suboptimality_Lipschitz_short} and \Cref{assu:WIRL,assu:WIRL-uniqueness,assu:IOP_std}.
    Assume that $\Theta$ is contained in some $(D-1)$-dimensional affine subspace (e.g., $\Theta = \Delta^{D-1}$).
    For a sample $(S^{(1)} , \ldots , S^{(N)} )$, let $\theta^{*(N)}\in \Theta$ and $\phi^{*(N)}\in \Phi$ denote, respectively, the weight and the constraint parameter output by \Cref{alg:Solve_IOP} with $\varepsilon = 0$, where line~3 is implemented by \Cref{alg:psgd_suboptimality_loss} with a sufficiently large sample-independent $K$ and with the learning rate chosen as either SRSS or SRSL.
    Assume that
    $\Phi (\delta) / \sim  = \{ [ \phi^{\mathrm{true} } ] \}$
    .
    Assume that $\Phi (\delta)$ is at most countable, and that for any $\phi \in \Phi (\delta)$ and $s \in \mathcal{S}$, $\sup_{x \in \mathcal{X} (\phi , s)} \left\| f (x , s) - f (\hat{x}^{*} (s) , s) \right\| < \infty$.
    Assume that $\sum_{j=1}^J \mathbb{P} \left( \phi_j^{\sup} (\{ S \}) > \phi_j^{\mathrm{true}} + \delta_j \right)^N < 1$.
    Then
    \begin{align*}
        & \mathbb{E}_{S^{(1)} , \ldots , S^{(N)}} \left[
            \mathbb{E}_{S} \ell^{\mathrm{sub},\lambda} ( \theta^{*(N)} , \phi^{*(N)} , S )
            \, \middle| \, \phi^{*(N)} \in \Phi (\delta)
        \right]
        \\
        & \leq
        \frac{L_f}{\sqrt{N}}
            \left( 1 - \sum_{j=1}^J \mathbb{P} \left( \phi_j^{\sup} (\{ S \} ) >  \phi_j^{\mathrm{true}} + \delta_j \right)^N \right)^{-1}
            32 \sqrt{(1 + \log 2) \log 2}\, C_{\mathcal{X}}
            \sqrt{D-1}
            \notag
        .
    \end{align*}
\end{theorem}

\begin{theorem}[{See \Cref{theo:subpoptimality_learning_prob_pl} for details.}]
    \label{theo:subpoptimality_learning_prob_pl_short}
    Suppose the setting of \Cref{theo:subpoptimality_learning_E_pl_short}.
    Then, for any $u \geq 0$,
        \begin{align*}
        & \mathbb{P} \Bigl( \phi^{*(N)} \sim \phi^{\mathrm{true}} \\
        & \qquad \text{and } \mathbb{E} \ell^{\mathrm{sub},\lambda} ( \theta^{*(N)} , \phi^{*(N)} , S ) \\
        & \qquad \leq \frac{32 \sqrt{3 (1 + \log 2)}\, L_f C_{\mathcal{X}} }{\sqrt{N}} \times \left( 2 \sqrt{2 \log 2} \sqrt{D-1} + u \, \mathrm{diam} (\Theta) \right) \Bigr) \\
        & \geq 1 - \sum_{j=1}^J \mathbb{P} \left( \phi_j^{\sup} (\{ S \} ) >  \phi_j^{\mathrm{true}} + \delta_j \right)^N - 5 \exp \left( -3  u^2  \right) .
    \end{align*}
\end{theorem}
Assumptions (i) and (ii) in \Cref{theo:subpoptimality_learning_E_pl_short,theo:subpoptimality_learning_prob_pl} guarantee that, for fixed $N$, a sample-independent $K$ exists. Indeed, under (i) there are only finitely many samples of size $N$. Under (ii), since $\phi^{\sup} (\mathcal{S}_N) \geq \phi^{\mathrm{true}}$ (\Cref{prop:true_constraints_para_dominated}) and by \Cref{prop:land-morphism_monotone}, we have $\mathcal{X} (\phi^{\sup} (\mathcal{S}_N) , s) \subset \mathcal{X} (\phi^{\mathrm{true}} , s)$, so that, by the boundedness of the ranges in \Cref{assu:WIRL}, $\| m \|_2$ and $L (\ell^{\mathrm{sub}}) \leq \| m \|_2$ in \Cref{prop:iteration_count_solve_IOP} are bounded independently of the sample. Without (i) and (ii), the same conclusions hold if $K$ is chosen sufficiently large for each sample (\Cref{theo:suboptimality_learning_prob,theo:suboptimality_learning_E}).

Under the assumption $\Phi (\delta) / \sim = \{ [ \phi^{\mathrm{true}} ] \}$, we have $\phi^{*(N)} \sim \phi^{\mathrm{true}}$ on the event $\phi^{*(N)} \in \Phi (\delta)$. Hence \Cref{theo:subpoptimality_learning_E_pl_short,theo:subpoptimality_learning_prob_pl_short} can be read as a decomposition into the probability of recovering the constraints up to equivalence and the generalization error of the objective given this recovery.
Constraint learning is deliberately conservative: $\phi^{\sup}$ is the tightest parameter consistent with the observed solutions, so it admits no decision beyond the data, and this tightening preserves the optimality of the observed solutions (\Cref{prop:true_para_include_sup_para}). The cost of this conservativeness on unseen states is the probability term $\sum_{j=1}^J \mathbb{P} ( \phi_j^{\sup} (\{ S \}) > \phi_j^{\mathrm{true}} + \delta_j )^N$, which decays exponentially in $N$ whenever each constraint is active at the observed solution up to $\delta_j$ with positive probability, i.e., $\mathbb{P} ( \phi_j^{\sup} (\{ S \}) \leq \phi_j^{\mathrm{true}} + \delta_j ) > 0$ for each $j$.

\section{Numerical Experiments}
\label{sec:numerical_experiment}

We empirically demonstrate that \Cref{alg:Solve_IOP} can solve the DDIOP for MILPs (cf.\ \Cref{eq:IOP_linear}).
Implementation details and the computational environment are provided in \Cref{sec:devices}.
Let $\mathcal{S}_N$ be the set of $N$ states sampled independently from $\mathcal{S}$; \Cref{alg:Solve_IOP} is run with respect to the uniform distribution on $\mathcal{S}_N$.

\subsection{Single-Machine Scheduling with Release Dates and the Weighted Sum of Completion Times}
\label{sec:1-machine}

\begin{figure}[!htb]
    \vspace{-\intextsep}
    \centering
    \begin{minipage}[c]{0.44\textwidth}
        \centering
        \includegraphics[width=\linewidth]{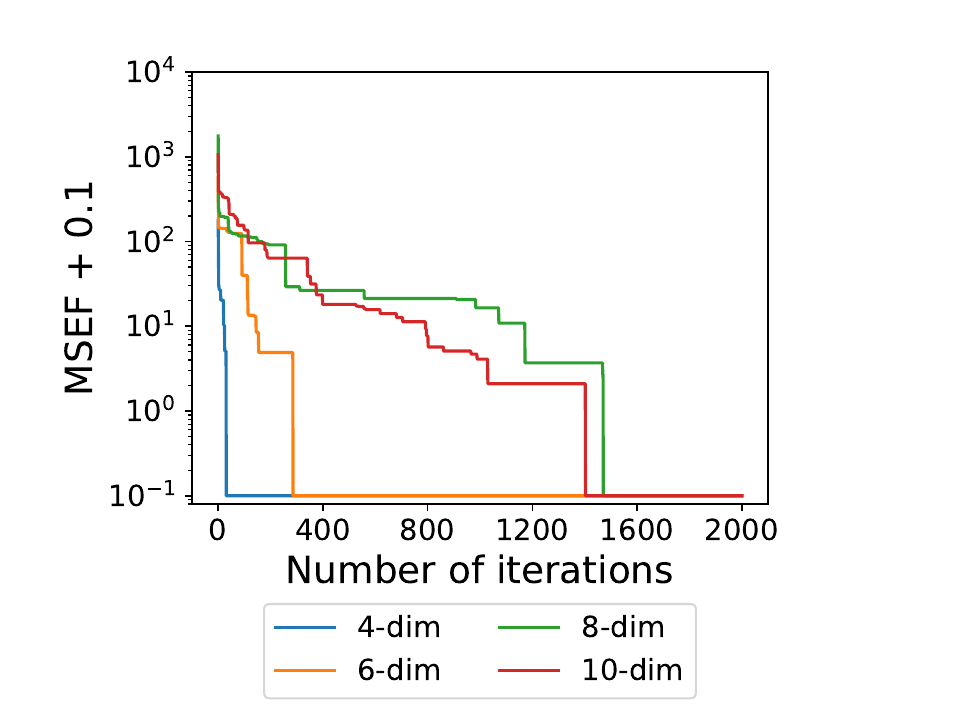}
        \vspace{-2\intextsep}
        \caption{Worst-case MSEF over all episodes at each iteration of \Cref{alg:psgd_suboptimality_loss} (\Cref{sec:1-machine}, $N = 5$, $10$ episodes, SRSL). The iteration at which a curve reaches $0$ is the worst-case number of iterations required to solve \Cref{eq:IOP_linear}. The vertical axis reports the mean squared error with respect to the feature (MSEF), plotted with an offset of $0.1$ to display $0$ on the logarithmic axis.}
        \label{fig:1-machine_N5_short}
    \end{minipage}
    \hfill
    \begin{minipage}[c]{0.55\textwidth}
        \centering
        \captionof{table}{Wall-time required to solve \Cref{eq:IOP_linear}}
        \label{tab:1-machine_N5_short}
        \setlength{\tabcolsep}{3pt}
        \begin{tabular}{lp{1.4cm}p{1.2cm}p{1.1cm}p{1.2cm}}
            \toprule
            $D$ & Decision variables & Con\-straints & Mean (s) & Median (s) \\
            \midrule
            4  & 16  & 40  & 0.13   & 0.09 \\
            6  & 36  & 96  & 2.24   & 1.63 \\
            8  & 64  & 176 & 62.76  & 51.45 \\
            \textbf{10} & \textbf{100} & \textbf{280} & \textbf{325.19} & \textbf{99.00} \\
            \bottomrule
        \end{tabular}
    \end{minipage}
        \vspace{-\intextsep}
\end{figure}

The single-machine scheduling problem of minimizing the weighted sum of completion times with release dates, denoted by $1|r_i|\sum \theta_i C_i$, is an ILP.
A detailed problem description is provided in \Cref{sec:1-machine_appx}.
We implemented \Cref{alg:Solve_IOP} using \Cref{alg:psgd_suboptimality_loss} with the SRSL learning-rate schedule.
Using this implementation, we conducted $10$ episodes with $N=5$ for $D=4,6,8,10$ (a subset of the settings $D = 4, \ldots, 10$ in \Cref{sec:1-machine_appx}), and report in \Cref{fig:1-machine_N5_short} the number of iterations in \Cref{alg:psgd_suboptimality_loss} required until \Cref{eq:IOP_linear} is solved.
As shown in \Cref{fig:1-machine_N5_short}, we empirically verified that \Cref{eq:IOP_linear} can be solved for all values of $D$ considered.
We also summarize in \Cref{tab:1-machine_N5_short} the wall-clock time required to solve \Cref{eq:IOP_linear}.
Even with $100$ decision variables, the problem is solved with a mean of $325$ seconds and a median of $99$ seconds.
For experimental results under other parameter settings, see \Cref{sec:1-machine_appx}.

\subsection{Single-Machine Scheduling with Release Dates and the Weighted Sum of Tardiness}
\label{sec:tardiness}
The single-machine scheduling problem with release dates of minimizing the weighted sum of tardiness \citep[\S~3]{PostekZoccaAMPL2024} can be formulated as an integer linear program (ILP).
A detailed problem description is provided in \Cref{sec:tardiness_appx}.
Using this implementation, we conducted experiments with $N=1$ and $D$ ranging from $3$ to $10$, performing $25$ episodes for each value of $D$.
The number of iterations required for \Cref{alg:Solve_IOP}—implemented with \Cref{alg:psgd_suboptimality_loss} using the SRSS learning-rate schedule—to solve \Cref{eq:IOP_linear} is shown in \Cref{fig:tardiness_SRSS_N1} (\Cref{sec:tardiness_appx}).
The number of iterations required for \Cref{alg:Solve_IOP}—implemented with \Cref{alg:psgd_suboptimality_loss} using the SRSL learning-rate schedule—to solve \Cref{eq:IOP_linear} is shown in \Cref{fig:tardiness_SRSL_N1} (\Cref{sec:tardiness_appx}).
From \Cref{fig:tardiness_SRSS_N1,fig:tardiness_SRSL_N1}, we confirmed that \Cref{eq:IOP_linear} can be solved for all values of $D$ considered, regardless of the learning rates.

For these figures and the behavior of the computational time and the suboptimality loss, see \Cref{sec:tardiness_appx}.

\subsection{Verifying the Generalization Error}
\label{sec:heldout}

On the same problem as in \Cref{sec:tardiness} ($D = 5$), we examine whether the loss on unseen states is consistent with the generalization error bounds (\Cref{theo:subpoptimality_learning_E_pl_short,theo:subpoptimality_learning_prob_pl_short}). In each episode, we sample $N + 100$ states, run \Cref{alg:Solve_IOP} on $N$ of them, and estimate the expected loss $\mathbb{E}_S \ell^{\mathrm{sub},1} ( \theta^{*(N)} , \phi^{*(N)} , S )$ on the remaining $100$ held-out states ($N \in \{ 1, 2, 4, 8, 16, 32 \}$, $25$ episodes each). Details are given in \Cref{sec:heldout_appx}.

\begin{figure}[!htb]
    \centering
    \begin{minipage}[c]{0.48\textwidth}
        \centering
        \includegraphics[width=\linewidth]{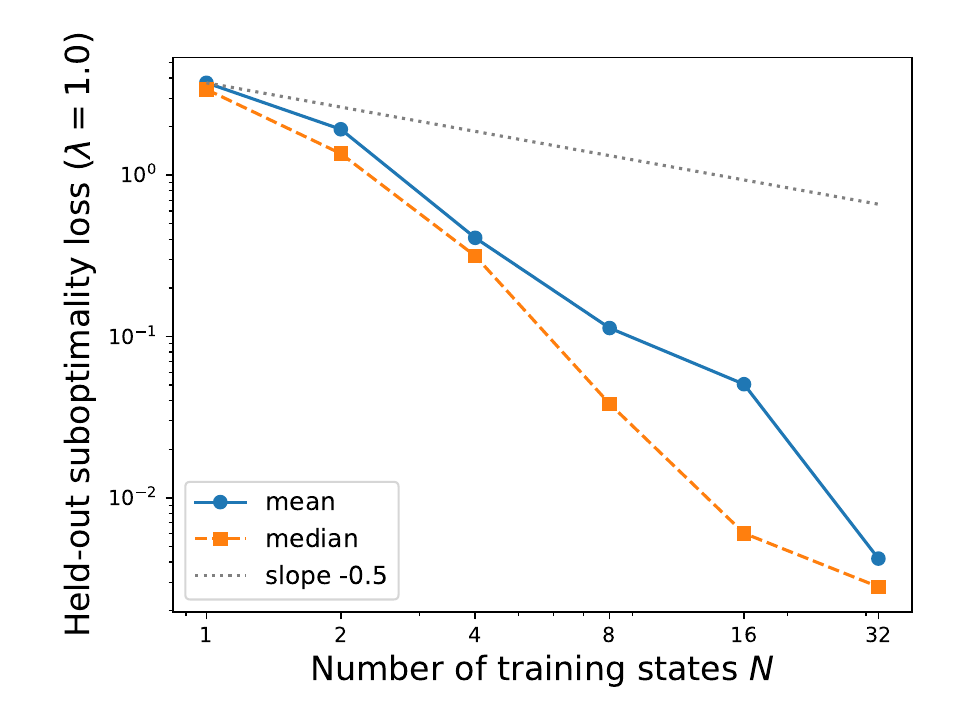}
        \caption{Mean and median of the suboptimality loss $\ell^{\mathrm{sub},1}$ on held-out states ($D = 5$, $25$ episodes for each $N$, log-log scale). The dotted line is a reference line of slope $-1/2$, corresponding to the rate of \Cref{theo:subpoptimality_learning_E_pl_short}.}
        \label{fig:heldout_SL}
    \end{minipage}
    \hfill
    \begin{minipage}[c]{0.48\textwidth}
        \centering
        \includegraphics[width=\linewidth]{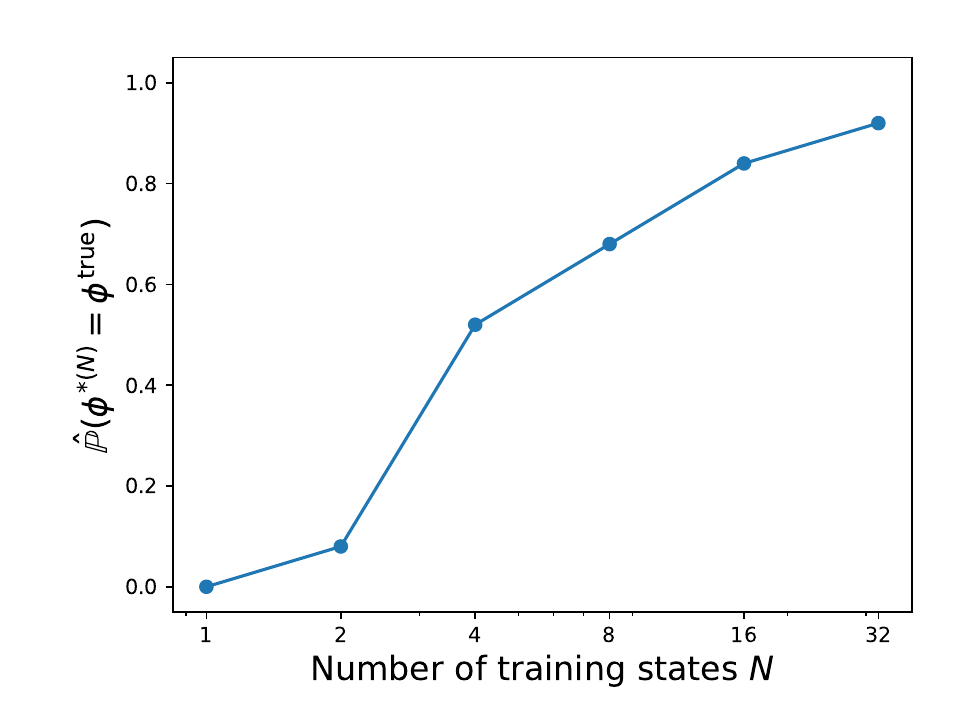}
        \caption{Exact recovery rate of the constraint parameter, $\hat{\mathbb{P}} ( \phi^{*(N)} = \phi^{\mathrm{true}} )$ ($D = 5$, $25$ episodes for each $N$).}
        \label{fig:phi_recovery}
    \end{minipage}
\end{figure}

The results are shown in \Cref{fig:heldout_SL,fig:phi_recovery}. The loss on the held-out states decays at least as fast as the theoretical rate $O(1/\sqrt{N})$ (the log-log slope of the mean is about $-1.9$). Moreover, the penalty term (the constraint side), dominant for small $N$, decays rapidly, and the exact recovery rate of the constraint parameters increases from $0\%$ at $N = 1$ to $92\%$ at $N = 32$. Already at $N = 8$, the penalty term is reduced by about $98\%$ (from $3.69$ to $0.085$; \Cref{tab:heldout_results}), so the conservativeness of constraint learning costs little with a moderate number of samples.

\FloatBarrier
\section{Conclusion}

We propose an efficient solution method for the inverse optimization problem of MILP. Specifically, we formulate a class of problems in which the objective function is represented as a linear combination of functions and the constraints are described by lattice homomorphisms, and propose a two-stage method that first learns the constraints and subsequently learns the objective function.

On the theoretical side, we show a theoretical guarantee of imitability under finite data distributions, develop statistical learning theory in pseudometric spaces and sub-Gaussian distributions, and apply it to derive generalization error bounds for inverse optimization. On the experimental side, we demonstrate that the proposed method exactly reproduces the observed optimal solutions of ILPs with up to $100$ decision variables, with a median computation time of $99$ seconds (mean $325$ seconds). This indicates that the proposed method is a practical solution for inverse optimization of scheduling problems of this size.

Finally, let us discuss future research directions. An important open question is whether there exist algorithms that solve the DDIOP for MILPs, learning both the constraints and the objective, faster than the proposed method. In the proposed method, constraint learning is given in closed form (\Cref{prop:constraint_learning_lattice_map}), and the number of iterations of objective learning is bounded by $O ( \gamma(\ell^{\mathrm{sub}})^{-2} )$ (\Cref{prop:iteration_count_solve_IOP}). Clarifying whether this iteration count can be improved, for example, through lower bounds on the number of calls to the forward problem required to solve the DDIOP, would lead to faster inverse optimization methods. Furthermore, investigating whether the generalization error bounds in inverse optimization, i.e.,~\Cref{theo:subpoptimality_learning_E_pl_short,theo:subpoptimality_learning_prob_pl_short}, are tight, as well as considering methods to obtain tighter bounds for inverse optimization, are important for designing faster inverse optimization algorithms.

\section*{AI Use Statement}
\label{sec:ai_statement}

In this work, we used generative AI tools for translation assistance (translating a Japanese draft of this paper into English) and for implementing the methods, i.e., writing and editing software code used in the numerical experiments.
The generative AI tools we used are Anthropic's Claude models (Claude Opus 4.6, 4.7, 4.8, 5, and 5.5, and Claude Fable 5 and 5.1).
We have not used generative AI tools for generating synthetic datasets, for developing the theoretical model or conceptual framework, for formulating the mathematical claims, for providing key elements of the proofs, for writing the proofs, for proposing or refining the hypotheses, for designing the research methods or experiments, or for interpreting the results, and curating or reformatting datasets and qualitative or thematic analysis are not applicable to this work.
Additionally, we used generative AI tools for editing to improve the readability of the paper and for formatting references.
We have reviewed all AI-assisted work.
All theoretical statements and their proofs were written and checked line by line by the author, and all AI-assisted code was verified and tested by the author against the experimental protocol reported in \Cref{sec:devices}.
We take responsibility for the final content of this work, including text, claims, or artifacts produced with the aid of generative AI.


\section*{Acknowledgement }
    We would like to thank Ryuta Matsuno for carefully reviewing this paper.

\bibliographystyle{apalike} 
\bibliography{suboptimality_loss.bib}

\newpage

\appendix

\crefalias{section}{appendix} 

\section{Related Works}
\label{sec:related_works_appendix}

\paragraph{Inverse optimization algorithms}

\citet{aswani2018inverse} proposed a method for learning both objective functions and constraints from states and optimization outcomes in convex programs (including linear programs), and showed its statistical consistency. However, the method presented in \citet{aswani2018inverse} encounters significant computational intractability issues.
As methods that estimate the constraints (the feasible region) via reduction to mathematical programming, there are \citet{chan2020inverse,ghobadi2020inferring}.
\citet{chan2020inverse} considered LPs and estimates the constraint matrix from a single observation. The cost vector is also included as a decision variable, but it serves only as an auxiliary variable that makes the observed solution optimal and is not estimated from data; an extension to multiple observations is left for future work.
\citet{ghobadi2020inferring} estimates the coefficient matrix and the right-hand side of the constraints (the feasible region) from multiple observations. However, neither of them learns the objective function.
Moreover, since \citet{chan2020inverse,ghobadi2020inferring} rely on LP duality, they are not directly applicable to integer or mixed-integer programs.
As inverse optimization that learns only the objective function (cost vector) in mixed-integer programs, there are \citet{bulut2021complexity,elhedhli2025inverse}.

\paragraph{Loss functions for inverse optimization}

The suboptimality loss for the objective function was introduced by \citet{mohajerin2018data}, who considered the setting where the constraints are known and only the objective function is learned.
\citet{ren2025inverse} extended this objective suboptimality loss to a form that also enables learning constraints (by adding a penalty term for constraint violation), and proposed the suboptimality loss, which evaluates whether the objective function and constraints have been correctly learned; in other words, whether the inverse optimization problem has been successfully solved.
The suboptimality loss is applicable not only to LP but also to MILP. The loss $\ell^{\mathrm{sub},\lambda}$ used in this work belongs to this line of work \citep{mohajerin2018data,ren2025inverse}.

\paragraph{Learning constraints}

\citet{ren2025inverse} proposed a method for learning constraints (a hypothesis class of feasible sets) from a dataset consisting of states and optimization results, under a known linear objective function.
The training problem is in general nonconvex and is handled by a smoothed local search (block coordinate descent), but it admits an exact reformulation as a mixed-integer linear program when the hypothesis class or the decision variables are restricted to be binary.
However, the objective function is assumed to be known; the method does not learn both the objective function and the constraints, and no generalization error analysis is given.

\citet{kolb2017learning,suenaga2024attempt} learn the constraint parameters with a pre-specified template for the constraints and a given two-dimensional (2-tensor) tabular dataset. \citet{kumar2019automating} uses a pre-defined constraint template to learn the constraint parameters from a 3-tensor dataset. These methods are superior in enabling constraint learning in MILP.

Our proposed method for learning constraints also uses constraint templates (cf.~\citet{kolb2017learning,kumar2019automating,suenaga2024attempt}) to learn constraints from the given dataset. The reason for adopting this method is that it enables learning constraints for both integer and mixed-integer cases.

\paragraph{Learning both constraints and objectives without learning guarantees}

Methods that learn both the constraints and the objective from optimal decisions without a learning guarantee also exist: \citet{tan2020learning_linear_programs} learn the parameters of LPs by solving the bilevel problem with gradient-based methods, which may be trapped in local minima, and \citet{paulus2021comboptnet} learn the constraints and costs of ILPs embedded as neural-network layers via surrogate gradients, and \citet{nandwani2022solver} learn them without calling an ILP solver during training, using a soft-margin loss that separates the observed solutions from generated negative samples. Since none of them guarantees that the parameters can be learned, we include them neither in \Cref{tab:compare_inverse_optimization} nor in our experiments.

\paragraph{Learning objective functions}

Inverse optimization methods for learning objective functions of MILP include methods based on suboptimality loss in the offline setting~\citep{suzuki2019TV,Kitaoka-2023-convergence-IRL,Kitaoka-2023-imitation-WIRL,kitaoka2024fast,zattoni2024learning} and the online setting~\citep{barmann2017emulating,barmann2018online,besbes2021online,besbes2025contextual,gollapudi2021contextual,sakaue2025online}.
In particular, \citet{kitaoka2024fast} showed that applying gradient-based optimization methods to the suboptimality loss reaches its minimum value $0$ (exact consistency with the observed data) within finitely many iterations, and bounded the number of iterations from above via a problem-dependent geometric constant $\gamma(\ell^{\mathrm{sub}})$.
Furthermore, \citet{kitaoka2026explicit} gave an explicit lower bound on $\gamma(\ell^{\mathrm{sub}})$ when the forward problem is an integer linear program, and thereby evaluated this iteration upper bound as an explicit function of the sample size, the dimension and ranges of the features, and the structure of the constraint coefficient matrix.
However, both consider the setting where the constraints are known and only the objective function is learned.
In this work, we use the iteration upper bound of \citet{kitaoka2024fast} in \Cref{prop:iteration_count_solve_IOP}, and the lower bound of \citet{kitaoka2026explicit} in \Cref{sec:apx_gamma_lowerbound}.%

\paragraph{Generalization error analysis for inverse optimization}
Several generalization error analyses are known for inverse optimization and predict-then-optimize.
\citet{elbalghiti2019generalization} established, in the predict-then-optimize framework, an $\tilde{O}(1/\sqrt{N})$ generalization bound for the SPO loss based on a margin and the Natarajan dimension.
\citet{li2025fenchel} provided generalization bounds and calibration bounds based on the Rademacher complexity for inverse optimization with the Fenchel--Young loss.
\citet{sun2023maximum} gave an excess-risk analysis for inverse linear programming.
\citet{mohajerin2018data} provided out-of-sample guarantees via distributionally robust optimization with the Wasserstein distance.
Moreover, \citet{aswani2018inverse} showed the statistical consistency (an asymptotic property) of estimators that learn both the objective function and the constraints in convex programs.
However, all of these either target learning only the objective function (or predict-then-optimize) or remain at asymptotic consistency; no finite-sample generalization error analysis has been given for inverse optimization that learns \emph{both} the constraints and the objective function in MILP.
This work provides the first generalization error analysis for this setting.%

\paragraph{Statistical learning theory}

As approaches for generalization error analysis, the use of Rademacher complexity (cf.~\citet{liao2020notes}) as well as results such as \citet[Theorem 8.3.23]{vershynin2020high}, \citet[Theorem 5]{shalev2009stochastic}, and \citet[Problem 5.12]{van2014probability} are known.
\citet[Theorem 8.3.23]{vershynin2020high} establishes a generalization bound under the assumption that the class generated by the parameters is a class of Boolean functions.
\citet[Theorem 5]{shalev2009stochastic} showed that, in a $D$-dimensional Euclidean space, if the loss function is Lipschitz continuous with respect to the parameters, the generalization error is $O_{\mathbb{P}} \left( \sqrt{\frac{D \log N}{N}} \right)$.
\citet[Problem 5.12]{van2014probability} demonstrated that, in a metric space, when the loss function is $L$-Lipschitz with respect to the parameters, the generalization error is $O_{\mathbb{P}} (L N^{-1/2})$.

One approach to generalization error analysis is to use Dudley's inequality~\citep{dudley1967sizes} (cf.~\citet[Theorem~8.3.23]{vershynin2020high}). Dudley's inequality bounds the expected supremum of a stochastic process by the covering number of the parameter space, where this covering number is defined with respect to a metric on the parameter space.
Among Dudley-integral-type bounds, one with small constants is \citet[Theorem~10.1]{lifshits2012lectures}. Note that the generic chaining bound via Talagrand's $\gamma_2$ functional (\citealp[Theorems 8.5.5 and 8.6.1]{vershynin2020high}) is in general sharper than Dudley-integral-type bounds.
On the other hand, as probabilistic inequalities (tail bounds) that bound the supremum of a stochastic process by a constant multiple of the integral of the covering numbers, we mention \citet[Theorem 8.1.6]{vershynin2020high} and \citet[Theorem 5.29]{van2014probability}.
Using such probabilistic inequalities, one can also perform generalization error analysis~\citep[Problem~5.12]{van2014probability}.

All of these results use covering numbers with respect to a \emph{metric} on the parameter space. In contrast, a framework that introduces a data-dependent \emph{pseudo-metric} on the hypothesis set and evaluates the generalization error via its covering number (a data-dependent fractal dimension) without assuming Lipschitz continuity has been studied by \citet{dupuis2023generalization} (\citet{simsekli2020hausdorff} give a bound based on the fractal dimension under the assumption of uniform Lipschitz continuity of the loss). Our statistical learning theory (\Cref{theo:learning_theory_expectation,theo:learning_theory_prob}) differs from these in that (i) it uses, as the distance, the pseudo-metric defined by the sub-Gaussian norm of the loss (rather than the loss values themselves), (ii) it carries out the direction of the chaining method suggested by \citet{dupuis2023generalization} in a setting different from theirs, namely with the pseudo-metric induced by the sub-Gaussian norm, and gives explicit constants, and (iii) it applies the resulting generalization bound to inverse optimization (suboptimality-loss minimization).%

\paragraph{Characterization of lattice homomorphisms on Euclidean space}

    \citet[Example~3.3]{dantas2021norm} provides a characterization of bounded linear lattice homomorphism functionals for several Banach lattice spaces; in particular, they give such a characterization on $\ell^p$ spaces.
    \Cref{theo:characterization_lattice_homomorphism} removes the bounded-linearity assumption and provides a characterization of lattice homomorphism functionals on standard Euclidean spaces.

\section{Lattice in Euclidean Space}
\label{sec:lattice}

For $\phi^{(1)}, \phi^{(2)} \in \mathbb{R}^{d}$, we write
\begin{align*}
    \phi^{(1)} \lor \phi^{(2)} & := ( \max (\phi^{(1)}_1, \phi^{(2)}_1 ), \ldots, \max (\phi^{(1)}_d, \phi^{(2)}_{d})), \\
    \phi^{(1)} \land \phi^{(2)} & := ( \min (\phi^{(1)}_1, \phi^{(2)}_1 ), \ldots, \min (\phi^{(1)}_d, \phi^{(2)}_{d})).
\end{align*}
For $\phi^{(1)}, \phi^{(2)} \in \mathbb{R}^{d}$, we write
\[
    \phi^{(1)} \leq \phi^{(2)} 
    \; :\Leftrightarrow\;
    \forall i = 1, \ldots , d,\,
    \phi_i^{(1)} \leq \phi_i^{(2)}.
\]
Let $\Phi \subset \mathbb{R}^d$ be a subset. A tuple $(\Phi , \land , \lor )$ is called a lattice if and only if for all $\phi^{(1)}, \phi^{(2)} \in \Phi$, $\phi^{(1)} \land \phi^{(2)} \in \Phi $,
     and $\phi^{(1)} \lor \phi^{(2)} \in \Phi $.
If $(\Phi, \land, \lor )$ is a lattice, a map $g \colon \Phi \to \mathbb{R}^{J}$ is a lattice homomorphism if and only if for all $\phi^{(1)}, \phi^{(2)} \in \Phi$,
\[
        g (\phi^{(1)} \land \phi^{(2)} ) = g (\phi^{(1)}) \land g(\phi^{(2)} ), 
        \quad
        g (\phi^{(1)} \lor \phi^{(2)} ) = g (\phi^{(1)}) \lor g(\phi^{(2)} ) .
\]

\section{Characterization of Lattice Homomorphisms}
\label{sec:characterization_lattice_homomorphism}

The following characterization of lattice homomorphisms underlies the examples satisfying \Cref{assu:IOP_std} (\Cref{exa:lattice_homomorphism_affine,exa:lattice_homomorphism_affine_check}).

\begin{theorem}
    \label{theo:characterization_lattice_homomorphism}
    Let $I_1,\ldots,I_d \subset \mathbb{R}$ be nonempty sets, and let $\Phi = \prod_{i=1}^d I_i$.\footnote{The triple $(\Phi,\land,\lor)$ forms a lattice.}
    Let $g \colon \Phi \to \mathbb{R}^J$.
    Then, $g$ is a lattice homomorphism if and only if, for each $j=1,\ldots,J$, there exist $i_j \in \{1,\ldots,d\}$ and a univariate monotone nondecreasing function $h_j$ such that $g_j(\phi)=h_j(\phi_{i_j})$ for any $\phi=(\phi_1,\ldots,\phi_d)\in\Phi$.
\end{theorem}

The proof is given later in this section.

\subsection{Order-Preserving Property of Lattice Homomorphisms}

\begin{definition}
    Let $\Phi \subset \mathbb{R}^{d}$ be a subset.
    Assume that $( \Phi , \land , \lor )$ forms a lattice.
    A map $g \colon \Phi \to \mathbb{R}^{J}$ is called a $\land$-homomorphism (resp. a $\lor$-homomorphism) if, for any $\phi^{(1)},\phi^{(2)} \in \Phi$, it holds that
    \[
        g (\phi^{(1)} \land \phi^{(2)} ) = g (\phi^{(1)}) \land g(\phi^{(2)})
        ,
        \quad
        \left( \text{resp. }
        g (\phi^{(1)} \lor \phi^{(2)} ) = g (\phi^{(1)}) \lor g(\phi^{(2)})
        \right)
        .
    \]
    The map $g \colon \Phi \to \mathbb{R}^{J}$ is called a lattice homomorphism if it is both a $\land$-homomorphism and a $\lor$-homomorphism.
\end{definition}

\begin{proposition}
    \label{prop:1-dim_set_is_and_closed}
    Let $\Phi \subset \mathbb{R}$.
    Then, $(\Phi , \lor , \land)$ is a lattice.
\end{proposition}

\begin{proof}
    Let $\phi^{(1)},\phi^{(2)} \in \Phi$ be arbitrary.
    By symmetry, we may assume $\phi^{(1)} \leq \phi^{(2)}$.
    Then, we have
    \[
        \phi^{(1)} \land \phi^{(2)}
        = \phi^{(1)} \in \Phi
        .
    \]
    Similarly, we can show that $\phi^{(1)} \lor \phi^{(2)} \in \Phi$.
\end{proof}

\begin{proposition}
    \label{prop:land-morphism_monotone}
    Let $\Phi \subset \mathbb{R}^{d}$ be a subset.
    Assume that $( \Phi , \land , \lor )$ forms a lattice.
    Then, for any $\land$-homomorphism or $\lor$-homomorphism $g \colon \Phi \to \mathbb{R}^J$, we have
    \[
        \phi^{(1)} \leq \phi^{(2)} 
        \Rightarrow
        g (\phi^{(1)})
        \leq
        g (\phi^{(2)})
        .
    \]
\end{proposition}

\begin{proof}
    Let $\phi^{(1)},\phi^{(2)} \in \Phi$ with $\phi^{(1)} \leq \phi^{(2)}$.
    Then $\phi^{(1)} \land \phi^{(2)} = \phi^{(1)}$ and $\phi^{(1)} \lor \phi^{(2)} = \phi^{(2)}$.
    If $g$ is a $\land$-homomorphism, we have
    \begin{align*}
        g ( \phi^{(1)})
        = g (\phi^{(1)} \land \phi^{(2)})
        = g (\phi^{(1)}) \land g(\phi^{(2)})
        \leq g (\phi^{(2)})
        .
    \end{align*}
    If $g$ is a $\lor$-homomorphism, we have
    \begin{align*}
        g ( \phi^{(2)})
        = g (\phi^{(1)} \lor \phi^{(2)})
        = g (\phi^{(1)}) \lor g(\phi^{(2)})
        \geq g (\phi^{(1)})
        .
    \end{align*}
\end{proof}

\subsection{Lattices on the Real Line}

\begin{proposition}
    \label{prop:one_variable_land-morphism_iff_monotone}
    Let $\Phi \subset \mathbb{R}$ be a nonempty set.
    Let $g \colon \Phi \to \mathbb{R}$.
    Then, the following statements are equivalent.
    \begin{description}
        \item[(1)] $g$ is a $\land$-homomorphism,
        \item[(2)] $g$ is a $\lor$-homomorphism,
        \item[(3)] $g$ is a lattice homomorphism,
        \item[(4)] $g$ is monotone increasing.
    \end{description}
\end{proposition}

\begin{proof}
    First, by Proposition~\ref{prop:1-dim_set_is_and_closed}, $( \Phi , \land , \lor )$ is a lattice.

    (1) $\Rightarrow$ (4): 
    If $g$ is a $\land$-homomorphism, then by Proposition~\ref{prop:land-morphism_monotone}, $g$ is monotone increasing.
    
    (4) $\Rightarrow$ (1): 
    Suppose that $g$ is monotone increasing.
    Let $\phi^{(1)},\phi^{(2)} \in \Phi$ be arbitrary.
    By symmetry, we may suppose $\phi^{(1)} \leq \phi^{(2)}$.
    Then, we have
    \begin{align*}
        g (\phi^{(1)})
        =
        g (\phi^{(1)} \land \phi^{(2)})
        \leq g ( \phi^{(1)}) \land g ( \phi^{(2)})
        \leq g ( \phi^{(1)})
    \end{align*}
    Therefore, it follows that $g$ is a $\land$-homomorphism.

    (2) $\Leftrightarrow$ (4) follows by an argument similar to (1) $\Leftrightarrow$ (4).

    (3) $\Rightarrow$ (1) is clear from the definition.
    
    (1) $\Rightarrow$ (3): 
    By (1) $\Rightarrow$ (4) and (4) $\Rightarrow$ (2), (3) holds.
\end{proof}

\subsection{Proof of Theorem~\ref{theo:characterization_lattice_homomorphism}}

\begin{proposition}
    \label{prop:characterization_lattice_map_bdd_dom}
    Let $I_1 , \ldots , I_d \subset \mathbb{R}$ be nonempty sets, each having minimum and maximum elements.
    Let $\Phi = \prod_{i=1}^d I_i$, and assume that $(\Phi , \land , \lor )$ forms a lattice.
    Let $g \colon \Phi \to \mathbb{R}$.
    Then, $g$ is a lattice homomorphism if and only if there exist $j \in \{1, \ldots , d\}$ and a monotone increasing univariate function $h$ such that $g(\phi) = h (\phi_j)$ for every $\phi = (\phi_1 , \ldots , \phi_d ) \in \Phi$.
\end{proposition}

\begin{proof}
    The sufficiency follows from \Cref{prop:one_variable_land-morphism_iff_monotone}.
    We show the necessity.

    For each $i=1 , \ldots , d$, let $s_i = \min I_i$ and $t_i = \max I_i$.
    By translation, we may assume $s_i = 0$.
    Let $\delta_{ij}$ denote the Kronecker delta, and let $\mathbf{e}_i = ( \delta_{ij })$ be the $i$-th standard basis vector in $\mathbb{R}^d$.
    Define $\hat{g} (\phi ) = g (\phi) - g(0)$, which is also a lattice map and $\hat{g} (0) = 0$.
    For any $i, i' \in \{ 1, \ldots , d \}$, $i \neq i'$, we have
    \[
        \hat{g} (t_i \mathbf{e}_i)  \land \hat{g} ( t_{i'} \mathbf{e}_{i'} ) 
        =
        \hat{g} (t_i \mathbf{e}_i \land t_{i'} \mathbf{e}_{i'} ) 
        = \hat{g} (0) = 0
    \]
    holds.
    Since $\hat{g} (t_i \mathbf{e}_i) \geq \hat{g} (0) = 0$ for any $i$ by \Cref{prop:land-morphism_monotone}, there exists $i^* \in \{ 1, \ldots , d \}$ such that, for any $i \in \{ 1, \ldots , d \} \setminus \{ i^* \}$, $\hat{g}(t_i \mathbf{e}_i) = 0 $.
    The map $t \mapsto \hat{g}(t \mathbf{e}_i)$ is a univariate lattice homomorphism, so by \Cref{prop:one_variable_land-morphism_iff_monotone}, this map is monotone increasing.
    Since $\hat{g}(0 \mathbf{e}_i) = \hat{g}(t_i \mathbf{e}_i) =0$, it follows that for any $\phi_i \in I_i $, $\hat{g}(\phi_i \mathbf{e}_i) =0$.

    Moreover, $t \mapsto \hat{g}(t \mathbf{e}_{i^*})$ is monotone increasing and $\hat{g}(0 \mathbf{e}_{i^*}) = 0 $, so for any $\phi_{i^*} \in I_{i^*} $, $\hat{g}(\phi_{i^*} \mathbf{e}_{i^*}) \geq 0 = \hat{g}(\phi_i \mathbf{e}_i)$.
    Therefore, for any $\phi \in \Phi$, we have
    \begin{align*}
        \hat{g} (\phi) 
        = \hat{g} \left( \bigvee_{i=1}^d \phi_i \mathbf{e}_i \right)
        = \bigvee_{i=1}^d \hat{g} \left( \phi_i \mathbf{e}_i \right)
        = \hat{g} (\phi_{i^*} \mathbf{e}_{i^*} )
        .
    \end{align*}
    Since $g (\phi ) = \hat{g} (\phi) + g(0)$, we have
    \begin{equation*}
        g (\phi ) = g (\phi_{i^*} \mathbf{e}_{i^*} )
        .
    \end{equation*}
    Since the map $I_{i^*} \ni \phi_{i^*} \mapsto g (\phi_{i^*} \mathbf{e}_{i^*})$ is monotone increasing, the proposition follows.
\end{proof}

\begin{proposition}
    \label{prop:characterization_lattice_map}
    Let $I_1 , \ldots , I_d \subset \mathbb{R}$ be nonempty sets.
    Let $\Phi = \prod_{i=1}^d I_i$.\footnote{The triple $(\Phi ,\land , \lor ) $ forms a lattice.}
    Let $g \colon \Phi \to \mathbb{R}$.
    Then, $g$ is a lattice homomorphism if and only if there exist $i \in \{1 , \ldots , d\}$ and a monotone increasing univariate function $h$ such that $g (\phi) = h (\phi_i)$ for any $\phi = (\phi_1 , \ldots , \phi_d )\in \Phi$.
\end{proposition}

\begin{proof}
    The sufficiency follows from \Cref{prop:one_variable_land-morphism_iff_monotone}.
    We show the necessity.

    Let $\{ s_i^m \}_{m \in \mathbb{Z}_{\geq 1 }} \subset I_i$ be a sequence converging monotonically decreasing to $\inf I_i$ and $\{ t_i^m \}_{m \in \mathbb{Z}_{\geq 1 }}  \subset I_i$ be a sequence converging monotonically increasing to $\sup I_i$.
    By \Cref{prop:characterization_lattice_map_bdd_dom}, the statement holds for each $\prod_{i=1}^d ( I_i \cap [s_i^m , t_i^m ] )$, i.e.,
    for each $m$, there exists a monotone increasing function $h^m $ and $i_m \in \{1, \ldots, d\}$ such that for any $\phi \in \prod_{i=1}^d ( I_i \cap [s_i^m , t_i^m ] )$, $g (\phi ) = h^m (\phi_{i_m})$.
    By assumption, for any $\phi \in \prod_{i=1}^d ( I_i \cap [s_i^m , t_i^m ] )$, $h^{m+1} (\phi_{i_{m+1}}) = h^m (\phi_{i_m})$.
    If $g$ is not constant on $\prod_{i=1}^d ( I_i \cap [s_i^m , t_i^m ] )$, then $i_{m+1} = i_m$; indeed, if $i_{m+1} \neq i_m$, then on this product set $g$ does not depend on $\phi_{i_{m+1}}$ since $g = h^m (\phi_{i_m})$, and does not depend on the other coordinates since $g = h^{m+1} (\phi_{i_{m+1}})$, so $g$ is constant there. If $g$ is constant on $\Phi$, the proposition holds with any $i$ and a constant function $h$. Otherwise, let $m_0$ be the smallest $m$ such that $g$ is not constant on $\prod_{i=1}^d ( I_i \cap [s_i^m , t_i^m ] )$. Then $i_m = i_{m_0}$ for any $m \geq m_0$, and for $m < m_0$ we may replace $i_m$ by $i_{m_0}$ since $h^m$ is constant. Hence we may assume $i_m = i_1$ for all $m$.
    Define the function $h \colon I_{i_1} \to \mathbb{R}$ by $h(\phi_{i_1}) := h^m (\phi_{i_1})$ for $\phi_{i_1} \in I_{i_1} \cap [s_{i_1}^m , t_{i_1}^m ]$.
    Then $h$ is well-defined and monotone increasing.
\end{proof}

\begin{proposition}
    \label{prop:lattice_map_iff_each_comp_lattice_map}
    Suppose that $(\Phi , \land , \lor )$ forms a lattice.
    Then, a map $g = (g_1 , \ldots , g_J ) \colon \Phi \to \mathbb{R}^J $ is a lattice homomorphism if and only if each $g_j \colon \Phi \to \mathbb{R}$ is a lattice homomorphism.
\end{proposition}

\begin{proof}
    The necessity is clear from the definition. Indeed,
    \begin{align*}
        \begin{pmatrix}
            g_1  (\phi \land \phi^\prime ) \\
            \vdots \\
            g_J  (\phi \land \phi^\prime )
        \end{pmatrix}
        = 
        g (\phi \land \phi^\prime ) 
        =
        g (\phi ) \land g (\phi^\prime )
        =
        \begin{pmatrix}
            g_1 (\phi) \land g_1 (\phi^\prime ) \\
            \vdots \\
            g_J (\phi) \land g_J (\phi^\prime )
        \end{pmatrix}
    \end{align*}
    holds.

    The sufficiency follows from
    \begin{align*}
        g (\phi \land \phi^\prime )
        =
        \begin{pmatrix}
            g_1  (\phi \land \phi^\prime ) \\
            \vdots \\
            g_J  (\phi \land \phi^\prime )
        \end{pmatrix}
        = 
        \begin{pmatrix}
            g_1 (\phi) \land g_1 (\phi^\prime ) \\
            \vdots \\
            g_J (\phi) \land g_J (\phi^\prime )
        \end{pmatrix}
        =
        g (\phi ) \land g (\phi^\prime )
    \end{align*}
    as required.
\end{proof}

{
\renewcommand{\proofname}{Proof of \Cref{theo:characterization_lattice_homomorphism}}
\begin{proof}
    The theorem follows from \Cref{prop:characterization_lattice_map,prop:lattice_map_iff_each_comp_lattice_map}.
\end{proof}
}

\section{Learning Constraints}
\label{sec:apx_constraint_learning}

{
\renewcommand{\proofname}{Proof of \Cref{prop:constraint_learning_lattice_map}}
\begin{proof}
    For a nonempty $\mathcal{S}^{\prime\prime} \subset \mathcal{S}$, let
    \[
        F (\mathcal{S}^{\prime\prime}) := \{ \phi \in \Phi \mid g ( \hat{x}^{*} (s) , \phi , s) \leq 0 \text{ for } s \in \mathcal{S}^{\prime\prime} \}
        .
    \]
    By \Cref{assu:IOP_std}, for any $s \in \mathcal{S}$, the set $F (\{ s \})$ has a maximal element, which is the maximum $\phi^{\sup} (\{ s \})$ since $F (\{ s \})$ is closed under $\lor$. Set $\phi^\prime := \bigwedge_{s \in \mathcal{S}^\prime} \phi^{\sup}(\{s\})$, which belongs to $\Phi$ since $(\Phi , \land , \lor)$ is a lattice and $\mathcal{S}^\prime$ is finite.
    Since for any $s \in \mathcal{S}^\prime$,
    \begin{align}
        g \left( \hat{x}^*(s),\, \bigwedge_{s^\prime \in \mathcal{S}^\prime} \phi^{\sup}(\{s^\prime\}),\, s \right)
        & =
        \bigwedge_{s^\prime \in \mathcal{S}^\prime}
        g \left( \hat{x}^*(s),\, \phi^{\sup}(\{s^\prime\}),\, s \right)
        \notag \\
        & \leq g \left( \hat{x}^*(s),\, \phi^{\sup}(\{s\}),\, s \right)
        \leq 0
        \notag
    \end{align}
    holds, we have $\phi^\prime \in F (\mathcal{S}^\prime)$. On the other hand, for any $\phi \in F (\mathcal{S}^\prime)$ and $s \in \mathcal{S}^\prime$, we have $\phi \in F (\{ s \})$ and hence $\phi \leq \phi^{\sup}(\{s\})$; thus $\phi \leq \phi^\prime$. Therefore, $\phi^\prime$ is the maximum of $F (\mathcal{S}^\prime)$, that is, $\phi^{\sup}(\mathcal{S}^\prime)$ exists and $\phi^{\sup}(\mathcal{S}^\prime) = \phi^\prime$.
\end{proof}
}

\begin{example}
    In the setting of \Cref{exa:lattice_homomorphism_affine_check} (with $\Phi^+ = \mathbb{R}^{J^+}$ and $\Phi^- = \mathbb{R}^{J^-}$, using the constraint map $\check{g}$ and the parameter space $\check{\Phi}$), for $s \in \mathcal{S}$,
    \[
        \phi^{\sup} (\{ s \}) = \left( 0 , - h^+ \bigl(\hat{x}^* (s) , s \bigr) , h^- \bigl(\hat{x}^* (s) , s \bigr) \right)
    \]
    holds. Applying \Cref{prop:constraint_learning_lattice_map}, we obtain
    \begin{align*}
        \phi^{\sup} ( \mathcal{S}^\prime )
        & = \Bigl( 0 , - \bigvee_{s \in \mathcal{S}^\prime } h^+ \bigl(\hat{x}^* (s) , s \bigr) , \bigwedge_{s \in \mathcal{S}^\prime } h^- \bigl(\hat{x}^* (s) , s \bigr) \Bigr)
    \end{align*}
    as claimed. That is, the learned feasible region is given by
    \[
        h^+ (x , s) \leq \bigvee_{s^\prime \in \mathcal{S}^\prime } h^+ (\hat{x}^* (s^\prime) , s^\prime )
        \quad \text{and} \quad
        h^- (x , s) \geq \bigwedge_{s^\prime \in \mathcal{S}^\prime } h^- (\hat{x}^* (s^\prime) , s^\prime )
        .
    \]
\end{example}

\section{Theory of Imitativeness}
\label{sec:apx_imitation_theory}

\begin{lemma}[{\citealp[Lemma~3.3]{kitaoka2024fast}}]
    \label{lem:Psi_set_is_almost_Phi}
    Assume \Cref{assu:WIRL} and let $\Theta = \Delta^{D-1} \subset \mathbb{R}^D$.
    Fix $\phi \in \Phi$.
    Let the state set $\mathcal{S} = \{ s^{(1)} , \ldots , s^{(N)} \}$ be finite.
    Then, for (almost) every $\theta \in \Delta^{D-1}$ (with respect to the measure on $\Delta^{D-1}$ induced by the Lebesgue measure), and for any $n=1,\ldots,N$, the vector $f(x^*(\theta,\phi,s^{(n)}), s^{(n)})$ is uniquely determined.
\end{lemma}

In this section, we prove that when $\mathcal{S}$ is a finite set and both $f$ and $g$ are piecewise linear maps, \Cref{alg:Solve_IOP} can be used to solve \cref{eq:IOP_linear}.

{
\renewcommand{\proofname}{Proof of \Cref{prop:true_constraints_para_dominated}}
\begin{proof}
    Suppose, to the contrary, that $\phi^{\sup} (\mathcal{S}^\prime) \not \geq \phi^{\mathrm{true}}$.
    Define $\phi^{\prime} := \phi^{\mathrm{true}} \lor \phi^{\sup} (\mathcal{S}^\prime)$. Then $\phi^{\prime} \in \Phi$ and $\phi^{\prime} > \phi^{\sup} (\mathcal{S}^\prime)$.
    For any $s \in \mathcal{S}^\prime$, we have $g (\hat{x}^{*} (s) , \phi^{\mathrm{true}} , s ) \leq 0$ since $\hat{x}^{*} (s) \in \mathcal{X} (\phi^{\mathrm{true}} , s)$, and since the map $g(\hat{x}^{*}(s),\bullet,s)$ is a lattice homomorphism, we have
    \begin{align}
        g (\hat{x}^{*} (s) , \phi^{\prime} , s )
        & = g (\hat{x}^{*} (s) , \phi^{\mathrm{true}} \lor \phi^{\sup} (\mathcal{S}^\prime) , s )
        \notag \\
        & = g (\hat{x}^{*} (s) , \phi^{\mathrm{true}} , s ) \lor g (\hat{x}^{*} (s) , \phi^{\sup} (\mathcal{S}^\prime) , s ) \leq 0
        \notag
    \end{align}
    This contradicts the choice of $\phi^{\sup} (\mathcal{S}^\prime)$.
\end{proof}
}

{
\renewcommand{\proofname}{Proof of \Cref{prop:true_para_include_sup_para}}
\begin{proof}
    By the definition of $\phi^{\sup}$,
    \[
    \hat{x}^{*}(s) \in 
        \left\{ 
            x^\prime \in \mathcal{X}
        \, \middle| \,
            g (x^\prime,\phi^{\sup},s) \leq 0
    \right\}
    .
    \]
    By \Cref{prop:true_constraints_para_dominated}, we have $\phi^{\mathrm{true}} \leq \phi^{\sup}$.
    Hence, by \Cref{prop:land-morphism_monotone},
    \begin{equation}
        \left\{ 
            x^\prime \in \mathcal{X}
        \, \middle| \,
            g (x^\prime,\phi^{\sup},s) \leq 0
        \right\}
        \subset
        \left\{ 
            x^\prime \in \mathcal{X}
            \, \middle| \,
            g (x^\prime,\phi^{\mathrm{true}},s) \leq 0
        \right\}
        \label{eq:true_para_include_sup_para_01}
    \end{equation}
    holds.
    Moreover, for any
    \[
        x^\prime \in 
        \left\{ 
            x^\prime \in \mathcal{X}
        \, \middle| \,
            g (x^\prime,\phi^{\sup},s) \leq 0
        \right\}
    \]
    we have, by \Cref{eq:true_para_include_sup_para_01} and $\hat{x}^{*}(s) \in \mathbf{FOP} (\theta, \phi^{\mathrm{true}}, s)$,
    \[
        \theta^{\top} f (x^\prime,s) \leq \theta^{\top} f (\hat{x}^{*} (s) ,s)
    \]
    . Therefore, $\hat{x}^{*}(s) \in \mathbf{FOP} (\theta, \phi^{\sup}, s)$.

    Finally, applying the first claim with $\theta = \theta^{\mathrm{true}}$, the objective value of any element of $\mathbf{FOP} (\theta^{\mathrm{true}}, \phi^{\sup}, s)$ equals the optimal value of $\mathbf{FOP} (\theta^{\mathrm{true}}, \phi^{\mathrm{true}}, s)$; hence by \Cref{eq:true_para_include_sup_para_01}, $\mathbf{FOP} (\theta^{\mathrm{true}}, \phi^{\sup}, s) \subset \mathbf{FOP} (\theta^{\mathrm{true}}, \phi^{\mathrm{true}}, s)$. Thus, by \Cref{assu:WIRL-uniqueness}, $f (x^*(\theta^{\mathrm{true}},\phi^{\sup},s) , s)$ is uniquely determined.
\end{proof}
}
\begin{theorem}
    \label{theo:alg:Solve_IOP_is_solve_IOP_classic}
    Suppose \Cref{assu:WIRL,assu:WIRL-uniqueness,assu:IOP_std}. Let the state set $\mathcal{S}$ be finite, and let $\mathbb{P}_{\mathcal{S}}$ be the uniform distribution.
    Take $\theta^{\sup} \in \Theta$ such that
    \begin{equation}
        \theta^{\sup} \in \argmin_{\theta \in \Theta }
        \mathbb{E}_{S} \ell^{\mathrm{sub},0} ( \theta , \phi^{\sup} (\mathcal{S}) ,S )
        .
        \label{eq:IOP_E_suboptimality_loss_sup}
    \end{equation}
    Then, for any $s \in \mathcal{S}$, we have
    $\hat{x}^{*}(s) \in \mathbf{FOP} (\theta^{\sup} , \phi^{\sup} (\mathcal{S}), s)$.
\end{theorem}

\begin{proof}
    By \Cref{prop:true_para_include_sup_para}, we have
    $\hat{x}^{*}(s) \in \mathbf{FOP} (\theta^{\mathrm{true}}, \phi^{\sup}, s)$.
    Therefore,
    \[
        0=
         \mathbb{E}_{S} \ell^{\mathrm{sub},0} ( \theta^{\mathrm{true} }, \phi^{\sup} (\mathcal{S}) ,S )
         =
         \min_{\theta \in \Theta }
          \mathbb{E}_{S} \ell^{\mathrm{sub},0} ( \theta , \phi^{\sup} (\mathcal{S}) ,S )
        .
    \]
    Since the state space $\mathcal{S}$ is finite, $\mathbb{P}_{\mathcal{S}}$ is the uniform distribution, and the suboptimality loss $\ell^{\mathrm{sub},0}$ takes nonnegative values, we have, for any $s \in \mathcal{S}$,
    \[
        \ell^{\mathrm{sub},0} ( \theta^{\sup} , \phi^{\sup} (\mathcal{S}) ,s ) = 0
        .
    \]
    By the definition of the suboptimality loss, this implies $\hat{x}^{*} (s ) \in \mathbf{FOP} (\theta^{\sup}, \phi^{\sup}, s)$.
\end{proof}

{
\renewcommand{\proofname}{Proof of \Cref{theo:alg:Solve_IOP_is_solve_IOP}}
\begin{proof}
    By \Cref{prop:true_para_include_sup_para}, we have
    $\hat{x}^{*}(s) \in \mathbf{FOP} (\theta^{\mathrm{true}}, \phi^{\sup}, s)$.
    By \Cref{prop:intention_learning_complete_Polyhedral}, the minimum value in line~3 of \Cref{alg:Solve_IOP} is $0$, and the minimum is attained within $k = O( \gamma ( \ell^{\mathrm{sub}} )^{-2} )$ iterations.
    The theorem follows from \Cref{theo:alg:Solve_IOP_is_solve_IOP_classic}.
\end{proof}
}

{
\renewcommand{\proofname}{Proof of \Cref{theo:alg:Solve_IOP_with_min_suboptloss_is_solve_IOP_LP}}
\begin{proof}
    By \Cref{lem:Psi_set_is_almost_Phi}, \Cref{assu:WIRL-uniqueness} holds for almost every $\theta^{\mathrm{true}}\in\Delta^{D-1}$ (with respect to the Lebesgue measure on $\Delta^{D-1}$).
    Fix such a $\theta^{\mathrm{true}}$. Then \Cref{theo:alg:Solve_IOP_is_solve_IOP} applies, which yields the claim.
\end{proof}
}

\section{Iteration-count bounds and a lower bound on the constant $\gamma(\ell^{\mathrm{sub}})$}
\label{sec:apx_gamma_lowerbound}

The suboptimality loss is known to satisfy the following properties.
\begin{proposition}[{\citealt[Proposition~3.1]{barmann2018online}; \citealt[Lemma~4.8]{Kitaoka-2023-convergence-IRL}}]
    \label{prop:SL_is_Lipschitz_convex}
    Suppose that \Cref{assu:WIRL} holds. Let the state set $\mathcal{S}$ be finite, and let the probability distribution $\mathbb{P}_{\mathcal{S}}$ be the uniform distribution.
    Fix $\phi \in \Phi$ such that $\hat{x}^{*} (s) \in \mathcal{X} (\phi , s)$ for any $s \in \mathcal{S}$, and consider $\theta \mapsto \mathbb{E} \ell^{\mathrm{sub},0} ( \hat{x}^{*} (S) , \theta , \phi , S )$.
    Then the following statements hold:
    (A) the suboptimality loss $\ell^{\mathrm{sub},0}$ is convex;
    (B) the suboptimality loss $\ell^{\mathrm{sub},0}$ is Lipschitz continuous; and
    (C) a subgradient of the suboptimality loss $\ell^{\mathrm{sub},0}$ at $\theta \in \Theta$ is given by
        $
       F (\theta) = 
            \mathbb{E} \left( f ( x^* (\theta  , \phi, S ) ,S)- f (\hat{x}^* (S) , S) \right) 
        $
        .
\end{proposition}

Let the diameter of $\mathcal{X}$ be $\mathrm{diam} (\mathcal{X}) = \sup_{x , x^\prime \in \mathcal{X}} \| x - x^\prime \|$, and $\mathrm{diam}(\Theta)$ is as defined in the main text. Following \citet{kitaoka2024fast}, we define the constants that appear in the iteration-count bound (\Cref{prop:iteration_bound_objective}) as follows.
\begin{proposition}[{\citealp[Proposition~10.14]{kitaoka2024fast}}]
    \label{prop:image_PL_union_polytopes}
    Let $X \subset \mathbb{R}^{n_1}$ be a finite union of polytopes (the union of finitely many convex polytopes). Let $T \colon \mathbb{R}^{n_1} \to \mathbb{R}^{n_2}$ be piecewise linear, i.e., there exist closed convex sets $Q_1 , \ldots , Q_L$, each defined by finitely many linear inequalities, such that $\mathbb{R}^{n_1} = \bigcup_{l=1}^L Q_l$ and $T$ is affine on each $Q_l$. Then $T(X)$ is a finite union of polytopes, and $\mathrm{Conv}\, T(X)$ is a convex polytope.
\end{proposition}

\begin{definition}[{\citealp[Eqs.~(4.2) and (4.3)]{kitaoka2024fast}}]
    \label{def:gamma_L_sub}
    Suppose \Cref{assu:WIRL,assu:WIRL-uniqueness}. Let the state set $\mathcal{S}$ be finite, and fix $\phi \in \Phi$. For any $s \in \mathcal{S}$, the convex hull $\mathrm{Conv} f ( \mathcal{X} (\phi, s ) ,s )$ is a convex polytope (by \Cref{assu:WIRL} and \Cref{prop:image_PL_union_polytopes}); let $Y(\phi , s)$ denote its vertex set. Define the Lipschitz constant $L(\ell^{\mathrm{sub}})$ of the suboptimality loss $\ell^{\mathrm{sub},0}$ (under $\phi$) and the positive geometric constant $\gamma(\ell^{\mathrm{sub}})$ by
    \begin{align*}
        L(\ell^{\mathrm{sub}})
        &:= \sup_{\substack{a^s \in f(\mathcal{X}(\phi,s),s),\\ s \in \mathcal{S}}} \left\| \frac{1}{|\mathcal{S}|} \sum_{s \in \mathcal{S}} \left( a^s - f(\hat{x}^*(s),s) \right) \right\|
        \quad \left( \leq L_f \, \mathrm{diam}(\mathcal{X}) \right), \\
        \gamma(\ell^{\mathrm{sub}})
        &:= \max_{\theta \in \Theta} \min_{( a^{s} )_s \in \prod_{s \in \mathcal{S}} Y(\phi , s) \setminus \{ f (\hat{x}^* (s), s) \}_s } \frac{1}{|\mathcal{S}|} \sum_{s \in \mathcal{S}} \theta^\top \left( f (\hat{x}^* (s),s) - a^s \right)
    \end{align*}
    where $L_f$ is the Lipschitz constant of $f$ (\Cref{assu:WIRL}).
\end{definition}

In this setting, \Cref{alg:psgd_suboptimality_loss} can attain the minimum value $0$ of the suboptimality loss \citep[Corollaries~4.6--4.9]{kitaoka2024fast}.

\begin{proposition}[{\citealp[Corollaries~4.7, 4.9]{kitaoka2024fast}}]
    \label{prop:intention_learning_complete_Polyhedral}
    Suppose \Cref{assu:WIRL,assu:WIRL-uniqueness}.
    Let the state set $\mathcal{S}$ be finite, and let the probability distribution $\mathbb{P}_{\mathcal{S}}$ be the uniform distribution.
    Fix $\phi \in \Phi$, and suppose that the map $\hat{x}^{*} \colon \mathcal{S} \to \mathcal{X}$ satisfies $\hat{x}^{*}(s) \in \mathbf{FOP}(\theta^{\mathrm{true}},\phi,s)$ for any $s \in \mathcal{S}$ and that $f (\mathbf{FOP} (\theta^{\mathrm{true}} , \phi , s ) , s)$ is a singleton.
    Let $\{ \theta^k \}_k$ be the sequence generated by \Cref{alg:psgd_suboptimality_loss} with the learning rate chosen as either SRSS or SRSL.
    Then, after finitely many iterations $k =  O (\gamma(\ell^{\mathrm{sub}})^{-2})$, we have $\mathbb{E} \ell^{\mathrm{sub}, 0} ( \hat{x}^{*} (S) , \theta^k , \phi , S ) = 0$.
\end{proposition}

\begin{remark}
    In \Cref{prop:intention_learning_complete_Polyhedral}, for each of the above learning-rate rules, one can explicitly bound the number of iterations required to reach the minimum value $0$ of $\mathbb{E} \ell^{\mathrm{sub}, 0}$ \citep[Corollaries~4.7, 4.9]{kitaoka2024fast}. We quote the explicit form, in our notation, in \Cref{prop:iteration_bound_objective} (\Cref{sec:apx_gamma_lowerbound}).
\end{remark}

For the geometric constant $\gamma(\ell^{\mathrm{sub}})$ (\Cref{def:gamma_L_sub}) appearing in this iteration-count bound, when the forward problem is an ILP, \citet{kitaoka2026explicit} provides an explicit lower bound in terms of the problem size. We quote the statement for general ILPs.
\begin{proposition}[\citealp{kitaoka2026explicit}]
    \label{prop:gamma_lowerbound_ILP_simple}
    Suppose \Cref{assu:WIRL,assu:WIRL-uniqueness}.
    Let the state set $\mathcal{S}$ be finite, and let $N = | \mathcal{S} |$.
    Fix $\phi \in \Phi$, and let the weight space be $\Theta = \Delta^{D-1}$ ($D \geq 2$).
    Suppose that the decision space satisfies $\mathcal{X} \subset \mathbb{Z}^{D_{\mathcal{X}}}$ and that the features satisfy $f (\mathcal{X} \times \mathcal{S}) \subset \mathbb{Z}^{D}$ (that is, the forward problem is an ILP).
    Define $m = (m_1 , \ldots , m_D)$ by
    \[
        m_i := \sup_{s \in \mathcal{S}} \left( \max_{x \in \mathcal{X} (\phi , s)} f_i (x , s) - \min_{x \in \mathcal{X} (\phi , s)} f_i (x , s) \right)
        .
    \]
    Then the constant $\gamma(\ell^{\mathrm{sub}})$ of \Cref{def:gamma_L_sub} satisfies
    \[
        \gamma(\ell^{\mathrm{sub}}) \ge \frac{1}{N^{D}\,\max(D-1,\sqrt{2})\,\|m\|_2^{D-1}}
        .
    \]
\end{proposition}
Sharper lower bounds that exploit the structure of the constraints (total unimodularity, Graver bases, and M-convexity) are quoted in \Cref{prop:gamma_lowerbound_ILP} (\Cref{sec:apx_gamma_lowerbound}).

We first quote from \citet{kitaoka2024fast}, in our notation, the explicit form of the iteration-count bound of \Cref{prop:intention_learning_complete_Polyhedral}.

\begin{proposition}[{\citealp[Corollaries~4.7, 4.9]{kitaoka2024fast}}]
    \label{prop:iteration_bound_objective}
    Suppose \Cref{assu:WIRL,assu:WIRL-uniqueness}.
    Let the state set $\mathcal{S}$ be finite, and let $\mathbb{P}_{\mathcal{S}}$ be the uniform distribution.
    Fix $\phi \in \Phi$, and suppose that the map $\hat{x}^{*} \colon \mathcal{S} \to \mathcal{X}$ satisfies $\hat{x}^{*}(s) \in \mathbf{FOP}(\theta^{\mathrm{true}},\phi,s)$ for any $s \in \mathcal{S}$ and that $f (\mathbf{FOP} (\theta^{\mathrm{true}} , \phi , s ) , s)$ is a singleton.
    Let $\alpha>0$ be the learning-rate constant, and let $L(\ell^{\mathrm{sub}})$ and $\gamma(\ell^{\mathrm{sub}})$ be the constants of \Cref{def:gamma_L_sub}.
    Then \Cref{alg:psgd_suboptimality_loss} with the SRSS learning rate attains $\mathbb{E} \ell^{\mathrm{sub},0} ( \hat{x}^{*} (S) , \theta^{k} , \phi , S ) = 0$ within at most
    \[
        \max\left(1,\;\left(\frac{\mathrm{diam}(\Theta)^{2}+(1+\log 2)\,\alpha^{2}\,L(\ell^{\mathrm{sub}})^{2}}{\alpha\,\gamma(\ell^{\mathrm{sub}})}\right)^{2}-2\right)
    \]
    iterations, and \Cref{alg:psgd_suboptimality_loss} with the SRSL learning rate within at most
    \[
        \max\left(1,\;\left(\frac{L(\ell^{\mathrm{sub}})\left(\mathrm{diam}(\Theta)^{2}+(1+\log 2)\,\alpha^{2}\right)}{\alpha\,\gamma(\ell^{\mathrm{sub}})}\right)^{2}-2\right)
    \]
    iterations.
\end{proposition}

Next, in the setting of \Cref{prop:gamma_lowerbound_ILP_simple}, we quote from \citet{kitaoka2026explicit} sharper lower bounds on $\gamma(\ell^{\mathrm{sub}})$ that exploit the structure of the constraints.

\begin{proposition}[\citealp{kitaoka2026explicit}]
    \label{prop:gamma_lowerbound_ILP}
    Suppose the assumptions of \Cref{prop:gamma_lowerbound_ILP_simple}, and let $\gamma(\ell^{\mathrm{sub}})$ be the constant of \Cref{def:gamma_L_sub}.
    Then the following hold.
    (1) Assume that the constraints are expressed as $\mathcal{X}(\phi,s)=\{x\in\mathbb{Z}^{D}\mid Ax\le b(s)\}$ ($A\in\mathbb{Z}^{p\times D}$) and that the features are the identity map $f_i(x,s)=x_i$. Let $C:=g_\infty(\widetilde{A})$ be the $\ell_\infty$ norm of the Graver basis of the slack-augmented matrix $\widetilde{A}=[A\mid I_p]$ (an integer $\ge1$ determined by $A$ alone). Then
    \[
        \gamma(\ell^{\mathrm{sub}})\ge\frac{1}{N\,\max(D-1,\sqrt{2})\,(2C\sqrt{D})^{D-1}}=\Omega\!\left(\frac{1}{N\,D^{(D+1)/2}\,(2C)^{D-1}}\right).
    \]
    In particular, if $A$ is totally unimodular (every square submatrix has determinant $0,\pm1$), then $C=1$ and $\gamma(\ell^{\mathrm{sub}})=\Omega\!\left(1/(N\,D^{(D+1)/2}\,2^{D-1})\right)$.
    (2) Assume that, for each $s \in \mathcal{S}$, the feasible region $\mathcal{X}(\phi,s)$ is an M-convex or M$^\natural$-convex set, and that the features are the identity map $f_i(x,s)=x_i$. Then
    \[
        \gamma(\ell^{\mathrm{sub}}) = \Omega\!\left( \frac{1}{N D^{2}} \right)
        ,
    \]
    and this lower bound does not depend on the ranges $m$ of the features.
\end{proposition}

Substituting the lower bound of \Cref{prop:gamma_lowerbound_ILP} into the SRSS iteration-count bound of \Cref{prop:iteration_bound_objective}, and since $\gamma(\ell^{\mathrm{sub}})$ enters the denominator squared, the totally unimodular case yields $O\!\left(N^{2}\,D^{D+1}\,4^{D-1}\,\mathrm{poly}(\mathrm{diam}(\Theta),\alpha,L(\ell^{\mathrm{sub}}))\right)$ iterations (substituting $C=1$ into $(2C)^{2(D-1)}$ gives $4^{D-1}$). See \citet{kitaoka2026explicit} for the proof and the details of $g_\infty(\widetilde{A})$.

\begin{theorem}[Number of iterations to solve the inverse optimization problem that learns both the constraints and the objective]
    \label{prop:iteration_count_solve_IOP}
    Suppose \Cref{assu:WIRL,assu:IOP_std}.
    Let the state set $\mathcal{S}$ be finite, let $N = | \mathcal{S} |$, and let $\mathbb{P}_{\mathcal{S}}$ be the uniform distribution.
    Let the weight space be $\Theta = \Delta^{D-1}$ ($D \geq 2$), and let $\varepsilon = 0$.
    Suppose that the decision space satisfies $\mathcal{X} \subset \mathbb{Z}^{D_{\mathcal{X}}}$ and that the features satisfy $f (\mathcal{X} \times \mathcal{S}) \subset \mathbb{Z}^{D}$ (that is, the forward problem is an ILP).
    Let $\alpha > 0$ be the learning-rate constant, and let $L(\ell^{\mathrm{sub}})$ and $\gamma(\ell^{\mathrm{sub}})$ be the constants of \Cref{def:gamma_L_sub} determined under the constraint parameter $\phi = \phi^{\sup}$.
    Moreover, define the vector $m$ as in \Cref{prop:gamma_lowerbound_ILP_simple} under $\phi = \phi^{\sup}$.
    Then, for almost every $\theta^{\mathrm{true}} \in \Delta^{D-1}$, the following hold.
    (1) \Cref{alg:Solve_IOP}, where line~3 is implemented by \Cref{alg:psgd_suboptimality_loss} with the SRSS learning rate, solves \Cref{eq:IOP_linear} within at most
    \[
        \max\left(1,\;\left(\frac{\mathrm{diam}(\Theta)^{2}+(1+\log 2)\,\alpha^{2}\,L(\ell^{\mathrm{sub}})^{2}}{\alpha\,\gamma(\ell^{\mathrm{sub}})}\right)^{2}-2\right)
    \]
    iterations.
    (2) \Cref{alg:Solve_IOP}, where line~3 is implemented by \Cref{alg:psgd_suboptimality_loss} with the SRSL learning rate, solves \Cref{eq:IOP_linear} within at most
    \[
        \max\left(1,\;\left(\frac{L(\ell^{\mathrm{sub}})\left(\mathrm{diam}(\Theta)^{2}+(1+\log 2)\,\alpha^{2}\right)}{\alpha\,\gamma(\ell^{\mathrm{sub}})}\right)^{2}-2\right)
    \]
    iterations.
    (3) The constant $\gamma(\ell^{\mathrm{sub}})$ appearing in (1) and (2) satisfies $\gamma(\ell^{\mathrm{sub}}) \ge 1 / (N^{D}\,\max(D-1,\sqrt{2})\,\|m\|_2^{D-1})$.
    In particular, the numbers of iterations in (1) and (2) are bounded from above by explicit functions of $N$, $D$, $\|m\|_2$, $\mathrm{diam}(\Theta)$, $\alpha$, and $L(\ell^{\mathrm{sub}})$.
\end{theorem}

{
\renewcommand{\proofname}{Proof of \Cref{prop:iteration_count_solve_IOP}}
\begin{proof}
    In \Cref{alg:Solve_IOP}, constraint learning (computing $\phi^{\sup}$ in line~2) is given in closed form by \Cref{prop:constraint_learning_lattice_map} and requires no iteration.
    Hence the number of iterations required to solve \Cref{eq:IOP_linear} equals that of \Cref{alg:psgd_suboptimality_loss} in line~3, which learns the objective.
    Here, by \Cref{prop:true_para_include_sup_para}, the observation satisfies $\hat{x}^{*}(s)\in\mathbf{FOP}(\theta^{\mathrm{true}},\phi^{\sup},s)$, so once $\phi^{\sup}$ is fixed, solving \Cref{eq:IOP_linear} reduces to inverse optimization of the objective alone with the constraints fixed at $\phi=\phi^{\sup}$.
    Moreover, by \Cref{lem:Psi_set_is_almost_Phi}, \Cref{assu:WIRL-uniqueness} holds for almost every $\theta^{\mathrm{true}}\in\Delta^{D-1}$.
    Therefore, (1) and (2) follow from \Cref{theo:alg:Solve_IOP_with_min_suboptloss_is_solve_IOP_LP} and \Cref{prop:iteration_bound_objective} applied with $\phi = \phi^{\sup}$ (the explicit form of the iteration bound of \Cref{prop:intention_learning_complete_Polyhedral}, quoted in \Cref{sec:apx_gamma_lowerbound}).
    Statement (3) follows by applying \Cref{prop:gamma_lowerbound_ILP_simple} with $\phi = \phi^{\sup}$.
\end{proof}
}

\begin{remark}
    For example, when the constraint coefficient matrix is totally unimodular, substituting the lower bound of \Cref{prop:gamma_lowerbound_ILP} into (1) of \Cref{prop:iteration_count_solve_IOP}, \Cref{alg:Solve_IOP} with the SRSS learning rate solves \Cref{eq:IOP_linear} within at most $O\!\left(N^{2}\,D^{D+1}\,4^{D-1}\,\mathrm{poly}\!\left(\mathrm{diam}(\Theta),\alpha,L(\ell^{\mathrm{sub}})\right)\right)$ iterations.
\end{remark}

\section{Statistical Learning Theory}
\label{sec:apx_learning_theory}

A random variable $S$ is said to be sub-Gaussian if there exists $t > 0$ such that
$
    \mathbb{E} \exp( S^2 /t^2 ) \leq 2 .
$
The sub-Gaussian norm of a random variable $S$ is defined as
\[
    \| S \|_{\psi_2} := \inf \left\{ t > 0 \;\middle|\; \mathbb{E} \exp\left( S^2 / t^2 \right) \leq 2 \right\}.
\]
Let $\ell \colon \Theta \times \mathcal{S} \to \mathbb{R}$ be a loss function.
Let $S$ be a random variable taking values in $\mathcal{S}$.
We assume $S^{(n)} \sim S$ are independent and identically distributed random variables.
Define $\theta^{*(N)}$ and $\theta^{*}$ to be any element of
\[
    \argmin_{\theta \in \Theta}
        \frac{1}{N} \sum_{n=1}^N \ell \left(\theta,\, S^{(n)} \right),
    \quad 
    \argmin_{\theta \in \Theta}
        \mathbb{E} \ell(\theta,\, S) ,
\]
respectively.

For any $s \in \mathcal{S}$, let $d_s$ be a pseudometric\footnote{A function $d \colon \Theta \times \Theta \to \mathbb{R}_{\geq 0 }$ is called pseudometric if for every $\theta, \theta^\prime \theta^{\prime\prime} \in \Theta$, (1) $d (\theta , \theta ) = 0 $,
(2) $d (\theta , \theta^\prime ) = d (\theta^\prime , \theta ) $,
(3) $d (\theta , \theta^{\prime \prime } ) \leq d (\theta , \theta^{\prime } ) + d (\theta^\prime , \theta^{\prime \prime } )$.
} 
on $\Theta$. For any $\theta, \theta' \in \Theta$, define
\[
    d_{\mathcal{S}} (\theta,\, \theta')
    := \| d_{S}(\theta,\, \theta') \|_{\psi_2} < \infty,
\]
where $\| \bullet \|_{\psi_2}$ denotes the sub-Gaussian norm.
We also define $\mathrm{diam} (\Theta , d_{\mathcal{S}}) := \sup_{\theta , \theta^\prime \in \Theta} d_{\mathcal{S}} (\theta , \theta^\prime )$.

Consider the situation where, for any $\theta, \theta' \in \Theta$ and $\mathbb{P}_{\mathcal{S}}$-almost every $s \in \mathcal{S}$,
\begin{equation}
    | \ell (\theta,\, s ) - \ell (\theta',\, s ) |
    \leq d_s (\theta,\, \theta')
    .
    \label{eq:assu:loss_estimated_pseudometric}
\end{equation}
In what follows, we assume that $\mathbb{E}_S | \ell (\theta , S) | < \infty$ for any $\theta \in \Theta$ and that $\theta^{*(N)}$ can be chosen as a measurable function of the sample $( S^{(1)} , \ldots , S^{(N)} )$. We also assume that the suprema over the sample and the covering integrals appearing in the following theorems and their proofs are measurable.
Define the centered empirical process by $X_{\theta} := \frac{1}{N} \sum_{n=1}^N \ell ( \theta , S^{(n)}) - \mathbb{E}_S \ell ( \theta , S)$ for $\theta \in \Theta$.

\subsection{Sub-Gaussian Random Variables}

\begin{proposition}[{Cf. \citet[Proposition 2.5.2]{vershynin2020high}}]
    \label{prop:subgaussian_1-->4}
    If
    \[
        \mathbb{P} (|S|\geq t) \leq 2 \exp (-t^2/K^2)
    \]
    then
    
    \[
        \| S \|_{\psi_2} \leq \sqrt{3} K
        .
    \]
    Conversely, if $\| S \|_{\psi_2} \leq K$, then for any $t \geq 0$, we have
    $\mathbb{P} (|S|\geq t) \leq 2 \exp (-t^2/K^2)$.
    
\end{proposition}

\begin{proof}
    
    For a nonnegative random variable $Y$, we have
    \[
        \mathbb{E} \exp (Y)
        = \int_0^{\infty} \mathbb{P} \left( \exp (Y) \geq s \right) ds
        = 1 + \int_0^{\infty} \mathbb{P} (Y \geq t) \exp (t) dt
    \]
    (the second equality follows from the fact that the integrand equals $1$ for $s \in [0,1]$, together with the change of variables $s = \exp(t)$).
    Applying the assumption with $Y = S^2 / (3K^2)$, we obtain
    \begin{align*}
        \mathbb{E} \exp \left( \frac{S^2}{3K^2} \right)
        & = 1 + \int_0^{\infty} \mathbb{P} \left( |S| \geq \sqrt{3} K \sqrt{t} \right) \exp (t) dt \\
        & \leq 1 + \int_0^{\infty} 2 \exp ( - 3t ) \exp (t) dt
        = 1 + \int_0^{\infty} 2 \exp ( - 2t ) dt
        = 2
    \end{align*}
    By the definition of the sub-Gaussian norm, $\| S \|_{\psi_2} \leq \sqrt{3} K$ follows.

    The converse follows from the fact that, by the monotone convergence theorem, $\mathbb{E} \exp ( S^2 / \| S \|_{\psi_2}^2 ) \leq 2$ and hence $\mathbb{E} \exp ( S^2 / K^2 ) \leq 2$ for $K \geq \| S \|_{\psi_2}$, together with Markov's inequality:
    \begin{align*}
        \mathbb{P} (|S| \geq t )
        & = \mathbb{P} \left( \exp \left( \frac{S^2}{K^2} \right) \geq \exp \left( \frac{t^2}{K^2} \right) \right)
        \\
        & \leq \mathbb{E} \exp \left( \frac{S^2}{K^2} \right) \exp \left( - \frac{t^2}{K^2} \right)
        \leq 2 \exp \left( - \frac{t^2}{K^2} \right)
    \end{align*}
    (cf. \citet[Proposition 2.5.2]{vershynin2020high}).
    
\end{proof}

\begin{proposition}[{\citealp[Cf.][Proposition 2.6.1]{vershynin2026}}]
    \label{prop:subgaussian_4-->5}
    Let $S$ be a random variable with mean $0$.
    Then, if $\| S \|_{\psi_2} = K$, for any $\lambda \in \mathbb{R}$,
    \[
        \mathbb{E} \exp( \lambda S) \leq \exp \left( \frac{(1 + \log 2)\lambda^2}{2}K^2 \right)
        .
    \]
\end{proposition}

\begin{proof}
    If $K = 0$, then $S = 0$ almost surely and the claim is trivial; hence assume $K > 0$. By the monotone convergence theorem, $\mathbb{E} \exp ( S^2 / K^2 ) \leq 2$. Moreover, $e^x \leq x + e^{x^2}$ for any $x \in \mathbb{R}$.
    If $| \lambda | K \leq 1$, then by $\mathbb{E} S = 0$ and Jensen's inequality for the concave function $y \mapsto y^{\lambda^2 K^2}$,
    \[
        \mathbb{E} e^{\lambda S}
        \leq \mathbb{E} e^{\lambda^2 S^2}
        = \mathbb{E} \left( e^{S^2 / K^2} \right)^{\lambda^2 K^2}
        \leq \left( \mathbb{E} e^{S^2 / K^2} \right)^{\lambda^2 K^2}
        \leq \exp ( \lambda^2 K^2 \log 2 )
        .
    \]
    If $| \lambda | K > 1$, then by $\lambda S \leq \lambda^2 K^2 / 2 + S^2 / ( 2 K^2 )$ and Jensen's inequality,
    \begin{align*}
        \mathbb{E} e^{\lambda S}
        & \leq e^{\lambda^2 K^2 / 2} \, \mathbb{E} e^{S^2 / ( 2 K^2 )}
        \leq e^{\lambda^2 K^2 / 2} \left( \mathbb{E} e^{S^2 / K^2} \right)^{1/2}
        \\
        & \leq \exp \left( \frac{\lambda^2 K^2}{2} + \frac{\log 2}{2} \right)
        \leq \exp \left( \frac{(1 + \log 2) \lambda^2 K^2}{2} \right)
        .
    \end{align*}
    Since $\log 2 \leq (1 + \log 2) / 2$, the proposition follows in both cases.
\end{proof}

\begin{proposition}[{\citealp[Cf.][Proposition 2.5.2]{vershynin2020high}}]
    \label{prop:subgaussian_5-->1}
    Let $S$ be a random variable with mean $0$.
    If, for any $\lambda \in \mathbb{R}$,
    \[
        \mathbb{E} \exp( \lambda S) \leq \exp (\lambda^2 K^2)
    \]
    then for any $t \geq 0$,
    \[
        \mathbb{P} (|S|\geq t) \leq 2 \exp \left( - \frac{t^2}{4 K^2} \right).
    \]
\end{proposition}

\begin{proposition}[{\citealp[Cf.][Proposition 2.6.1]{vershynin2020high}}]
    \label{prop:meanzero_subGaussian_sum}
    Let $S_i$ be independent sub-Gaussian random variables with mean $0$.
    Then,
    \[
        \left\| \sum_{i=1}^N S_i \right\|_{\psi_2}^2
        \leq 6 (1 + \log 2)
        \sum_{i=1}^N \left\| S_i \right\|_{\psi_2}^2
    \]
\end{proposition}

\begin{proof}
    
    Set $K_i := \| S_i \|_{\psi_2}$.
    By \Cref{prop:subgaussian_4-->5}, for each $i$ and any $\lambda \in \mathbb{R}$, we have
    $\mathbb{E} \exp (\lambda S_i ) \leq \exp \left( \frac{(1 + \log 2)\lambda^2}{2} K_i^2 \right)$.
    By independence,
    \[
        \mathbb{E} \exp \left( \lambda \sum_{i=1}^N S_i \right)
        \leq \exp \left( \frac{(1 + \log 2)\lambda^2}{2} \sum_{i=1}^N K_i^2 \right)
    \]
    Applying \Cref{prop:subgaussian_5-->1} with $K^2 = \frac{1 + \log 2}{2}\sum_{i=1}^N K_i^2$, for any $t \geq 0$, we have
    \[
        \mathbb{P} \left( \left| \sum_{i=1}^N S_i \right| \geq t \right)
        \leq 2 \exp \left( - \frac{t^2}{2 (1 + \log 2) \sum_{i=1}^N K_i^2 } \right)
    \]
    Applying \Cref{prop:subgaussian_1-->4} with $K^2 = 2 (1 + \log 2)\sum_{i=1}^N K_i^2$, we obtain
    \[
        \left\| \sum_{i=1}^N S_i \right\|_{\psi_2}^2
        \leq 3 \cdot 2 (1 + \log 2) \sum_{i=1}^N K_i^2
        = 6 (1 + \log 2) \sum_{i=1}^N K_i^2
    \]
    (cf. the proof of \citet[Proposition 2.6.1]{vershynin2020high}).
    
\end{proof}

\begin{proposition}[{\citealp[Cf.][Proposition 2.5.2]{vershynin2020high}}]
    Let $S$ be a sub-Gaussian random variable. Then
    \[
        \mathbb{E} |S| \leq \sqrt{\log 2} \| S \|_{\psi_2}
        .
    \]
\end{proposition}

\begin{proof}
    For any $t > \| S \|_{\psi_2}$, Jensen's inequality for the convex function $x \mapsto \exp (x^2 / t^2)$ yields $\exp ( ( \mathbb{E} |S| )^2 / t^2 ) \leq \mathbb{E} \exp ( S^2 / t^2 ) \leq 2$, i.e., $\mathbb{E} |S| \leq \sqrt{\log 2}\, t$. Letting $t \downarrow \| S \|_{\psi_2}$ yields the proposition.
\end{proof}

\begin{proposition}[{\citealp[Cf.][Lemma 2.6.8]{vershynin2020high}}]
    \label{prop:centering_subGaussian}
    Let $S$ be a sub-Gaussian random variable.
    Then $S - \mathbb{E}S$ is also sub-Gaussian and
    \[
        \| S - \mathbb{E} S \|_{\psi_2} 
        \leq 2 
        \| S \|_{\psi_2}
        .
    \]
\end{proposition}

\begin{proof}
    For any $t > \| S \|_{\psi_2}$, since $\exp (x) \geq 1 + x$, we have $1 + \mathbb{E} S^2 / t^2 \leq \mathbb{E} \exp (S^2 / t^2) \leq 2$; hence $\mathbb{E} S^2 < \infty$, and in particular $\mathbb{E} S$ is well defined.
    Let $S^\prime$ be an independent copy of $S$, and take any $t > \| S - S^\prime \|_{\psi_2}$.
    Since $x \mapsto \exp (x^2 / t^2)$ is convex and $S - \mathbb{E} S = \mathbb{E}_{S^\prime} [ S - S^\prime \mid S ]$, Jensen's inequality for conditional expectations yields
    \[
        \mathbb{E} \exp \left( \frac{(S - \mathbb{E} S)^2}{t^2} \right)
        \leq \mathbb{E} \exp \left( \frac{(S - S^\prime)^2}{t^2} \right)
        \leq 2 .
    \]
    Therefore, $\| S - \mathbb{E} S \|_{\psi_2} \leq \| S - S^\prime \|_{\psi_2}$.
    Since the sub-Gaussian norm satisfies the triangle inequality (cf.\ \citet[Exercise 2.5.7]{vershynin2020high}), we have $\| S - S^\prime \|_{\psi_2} \leq \| S \|_{\psi_2} + \| S^\prime \|_{\psi_2} = 2 \| S \|_{\psi_2}$, and the proposition follows.
\end{proof}

\subsection{Pseudometric Spaces and Metric Spaces}

Let $(\Theta, d)$ be a pseudometric space.
A family of subsets $\mathcal{N}$ of $\Theta$ is called an $\varepsilon$-cover of the pseudometric space $(\Theta, d)$ if for every $\mathcal{N}_i \in \mathcal{N}$, there exists $\theta^i \in \Theta$ such that $\mathcal{N}_i = \{ \theta \in \Theta \mid d(\theta, \theta^i) < \varepsilon \}$, and $\Theta = \bigcup_i \mathcal{N}_i$ holds.
The $\varepsilon$-covering number $N(\Theta , d , \varepsilon)$ of the pseudometric space $(\Theta, d)$ is defined as the minimal cardinality of such an $\varepsilon$-cover.
A subset $A$ of $\Theta$ is said to be \emph{$\varepsilon$-separated} in the pseudometric space $(\Theta, d)$ if $d (a , b) \geq \varepsilon$ holds for any distinct $a , b \in A$.
The $\varepsilon$-packing number $P(\Theta, d, \varepsilon)$ of the pseudometric space $(\Theta, d)$ is the maximal cardinality of an $\varepsilon$-separated subset.

\begin{proposition}[{\citealp[Cf.][Lemma 4.2.8]{vershynin2020high}}]
    \label{prop:covereing_leq_packing}
    Let $(\Theta, d)$ be a pseudometric space.
    Then,
    \[
        N (\Theta , d , \varepsilon)
        \leq P (\Theta , d , \varepsilon)
        \leq N (\Theta , d , \varepsilon /2 )
        .
    \]
\end{proposition}

\begin{proposition}
    \label{prop:packing_ineq}
    Let $(\Theta, d)$ be a pseudometric space.
    Let $\Theta^\prime \subset \Theta$.
    Then,
    \[
        P (\Theta^\prime , d , \varepsilon)
        \leq P (\Theta , d , \varepsilon)
        .
    \]
\end{proposition}

\begin{proof}
    Let $A$ be an $\varepsilon$-separated subset of $\Theta^\prime$.
    Since separation depends only on the values of $d$, $A \subset \Theta^\prime \subset \Theta$ is also an $\varepsilon$-separated subset of $\Theta$.
    The proposition follows.
\end{proof}

Let $(\Omega, \mathcal{B}, \mathbb{P})$ be a probability space.
Let $(\Theta,d)$ be a pseudometric space.
A process $(S_{\theta})_{\theta \in \Theta}$ is said to be a sub-Gaussian process if for each $\theta \in \Theta$, $S_{\theta}$ is a random variable on $(\Omega, \mathcal{B}, \mathbb{P})$, and there exists a constant $L \geq 0$ such that for any $\theta, \theta^\prime \in \Theta$,
\begin{equation}
    \| S_{\theta} - S_{\theta^\prime} \|_{\psi_2} \leq L d( \theta , \theta^\prime )
    \label{eq:subgaussian_process}
\end{equation}
is satisfied.

\begin{definition}
    \label{def:separable_process}
    A stochastic process $(S_{\theta})_{\theta \in \Theta}$ on a pseudometric space $(\Theta , d )$ is said to be \emph{separable} if there exists an at most countable subset $\Theta^\prime \subset \Theta$ that is dense in $(\Theta , d )$ such that, for any $\theta^{\prime\prime} \in \Theta$, almost surely
    \[
        \sup_{\theta \in \Theta} \left( S_{\theta} - S_{\theta^{\prime\prime}} \right)
        = \sup_{\theta \in \Theta^\prime} \left( S_{\theta} - S_{\theta^{\prime\prime}} \right)
        .
    \]
    In particular, if $\Theta$ is at most countable, or if $(\Theta , d )$ is separable and $\theta \mapsto S_{\theta}$ is almost surely $d$-continuous, then $(S_{\theta})_{\theta \in \Theta}$ is separable.
\end{definition}%

In a pseudometric space $(\Theta, d)$, define the equivalence relation $\theta \sim \theta^\prime$ if $d(\theta, \theta^\prime) = 0$.
Set $\Theta^* = \Theta / \sim$.
Let $[ \theta ]$ denote the equivalence class of $\theta \in \Theta$.
A metric on $\Theta^*$ is given by
$d^* ( [ \theta ] , [\theta^\prime]) = d (\theta , \theta^\prime )$.
Then, the space $(\Theta^*, d^*)$ forms a metric space (cf.~\cite[p.~58]{howes1995modern}).

For a sub-Gaussian process $(S_{\theta})_{\theta \in \Theta}$ on a pseudometric space $(\Theta, d)$, by defining $S_{[\theta]} := S_{\theta}$, the collection $(S_{\theta})_{\theta \in \Theta}$ induces a sub-Gaussian process on the metric space $(\Theta^*, d^*)$.
Moreover, if $(S_{\theta})_{\theta \in \Theta}$ is separable, then $(S_{[\theta]})_{[\theta] \in \Theta^*}$ is also separable (it suffices to take the image of the same countable set).

\subsection{Dudley-type Integral Inequalities}

\begin{proposition}[{Cf. \citealp{dudley1967sizes}; \citealp[Theorem 10.1]{lifshits2012lectures}}]
    \label{prop:dudley_integral_ineq_E}
    Let $(S_{\theta})_{\theta \in \Theta}$ be a separable stochastic process (\Cref{def:separable_process}) on a metric space $(\Theta, d)$ with $\mathbb{E} S_{\theta} = 0$ for any $\theta \in \Theta$.
    Assume further that there exists a constant $L > 0$ such that, for any $\theta , \theta^\prime \in \Theta$ and any $\lambda \in \mathbb{R}$,
    \begin{equation}
        \mathbb{E} \exp \left( \lambda ( S_{\theta} - S_{\theta^\prime} ) \right)
        \leq \exp \left( \frac{\lambda^2 L^2 d (\theta , \theta^\prime )^2}{2} \right)
        \label{eq:subgaussian_process_MGF}
    \end{equation}
    holds.
    Then,
    \[
        \mathbb{E} \sup_{\theta \in \Theta } S_{\theta}
        \leq 4\sqrt{2} L
             \int_0^{\infty} \sqrt{\log N( \Theta, d, \varepsilon ) } \, d \varepsilon
         .
    \]
\end{proposition}

\begin{proposition}
    \label{prop:dudley_integral_ineq_E_pseudo}
    \Cref{prop:dudley_integral_ineq_E} still holds for a process satisfying the same assumptions on a pseudometric space $(\Theta,d)$.
\end{proposition}

\begin{proof}
    Since $\sup_{\theta \in \Theta } S_{\theta} = \sup_{ [ \theta ] \in \Theta^* } S_{ [ \theta ] }$ and $N( \Theta, d, \varepsilon ) = N( \Theta^*, d^*, \varepsilon )$, and since separability and \cref{eq:subgaussian_process_MGF} are preserved under passing to the quotient metric space, the proposition follows from \Cref{prop:dudley_integral_ineq_E}.
\end{proof}

\begin{proposition}[{\citealp[Theorem 5.29]{van2014probability}}]
    \label{prop:dudley_integral_ineq_prob_old}
    Let $(S_{\theta})_{\theta \in \Theta}$ be a separable sub-Gaussian process (\Cref{def:separable_process}) on a metric space $(\Theta, d)$.
    Let $L > 0$ be the constant appearing in \cref{eq:subgaussian_process}.
    Then, for any $\theta^\prime \in \Theta$ and any $u \geq 0$,
    \[
        \mathbb{P}
        \left[
        \sup_{\theta \in \Theta } ( S_{\theta} - S_{\theta^\prime} ) \geq C L \left( 
            \int_0^{\infty} \sqrt{\log N( \Theta, d, \varepsilon ) } \, d \varepsilon + u\, \mathrm{diam} (\Theta , d) 
         \right)
         \right] 
         \leq 2 \exp \left( -u^2 \right)
         ,
    \]
    where $C= 6 (1+ 2 / \log 2)$.
\end{proposition}

\begin{remark}
    In \Cref{prop:dudley_integral_ineq_prob_old}, the value $C= 6 (1+ 2 / \log 2)$ follows from the proof of \citet[Theorem 5.29]{van2014probability}.
\end{remark}

\begin{proposition}
    \label{prop:dudley_integral_ineq_prob}
    Let $(S_{\theta})_{\theta \in \Theta}$ be a separable sub-Gaussian process (\Cref{def:separable_process}) on a metric space $(\Theta, d)$.
    Let $L$ be the constant appearing in \cref{eq:subgaussian_process}.
    Let $C > 6\sqrt{3}$.
    Then, for any $\theta^\prime \in \Theta$ and any $u \geq 0$,
    \begin{align*}
        & \mathbb{P}
        \left[
        \sup_{\theta \in \Theta } ( S_{\theta} - S_{\theta^\prime} ) > C L \left(
            \int_0^{\infty} \sqrt{\log N( \Theta, d, \varepsilon ) } d \varepsilon + u\, \mathrm{diam} (\Theta , d)
        \right)
        \right]
        \notag\\
        & \leq 2 \zeta \left( \frac{C^2}{36} -2 \right) \exp \left( -\frac{C^2}{36} u^2  \right)
         ,
    \end{align*}
    where $\zeta$ denotes the Riemann zeta function,
    \[
        \zeta (u) := \sum_{j=1}^{\infty} j^{-u}
        .
    \]
    
    In particular, for $C = 8\sqrt{2}$,
    \begin{align*}
        & \mathbb{P}
        \left[
        \sup_{\theta \in \Theta } ( S_{\theta} - S_{\theta^\prime} ) > 8 \sqrt{2} L \left(
            \int_0^{\infty} \sqrt{\log N( \Theta, d, \varepsilon ) } d \varepsilon + u\, \mathrm{diam} (\Theta , d)
        \right)
        \right]
        \\
        & \leq 5 \exp \left( - 3 u^2  \right)
         ,
    \end{align*}
    and for $C = 12$,
    \begin{align*}
        & \mathbb{P}
        \left[
        \sup_{\theta \in \Theta } ( S_{\theta} - S_{\theta^\prime} ) > 12 L \left(
            \int_0^{\infty} \sqrt{\log N( \Theta, d, \varepsilon ) } d \varepsilon + u\, \mathrm{diam} (\Theta , d)
        \right)
        \right]
        \\
        & \leq 3.3 \exp \left( - 4 u^2  \right)
         .
    \end{align*}
    
\end{proposition}

\begin{remark}
    The proof of \Cref{prop:dudley_integral_ineq_prob} was inspired by \citet{kadmos2025}.
\end{remark}

{
\renewcommand{\proofname}{Proof of \Cref{prop:dudley_integral_ineq_prob}}
\begin{proof}
    First, if
    \[
    \int_0^{\infty} \sqrt{\log N( \Theta, d, \varepsilon ) } d \varepsilon = \infty,
    \]
    then, since the threshold is $+ \infty$, the event
    \begin{equation*}
        \Bigl\{
        \sup_{\theta \in \Theta } ( S_{\theta} - S_{\theta^\prime} ) > C L \left( 
            \int_0^{\infty} \sqrt{\log N( \Theta, d, \varepsilon ) } d \varepsilon + u\, \mathrm{diam} (\Theta , d) 
        \right)
        \Bigr\}
    \end{equation*}
    is empty, so its probability is $0$, and thus the proposition follows trivially.
    Therefore, we may assume that
    \[
        \int_0^{\infty} \sqrt{\log N( \Theta, d, \varepsilon ) } d \varepsilon < \infty.
    \]
    Moreover, if $L = 0$, then for any $\theta \in \Theta$ we have $\| S_{\theta} - S_{\theta^\prime} \|_{\psi_2} = 0$, i.e., $S_{\theta} = S_{\theta^\prime}$ almost surely. Hence, when $\Theta$ is finite or countable, $\sup_{\theta \in \Theta} ( S_{\theta} - S_{\theta^\prime} ) = 0$ almost surely (the general case is reduced to a countable set by separability, as in the last paragraph of this proof), so the event with the strict inequality is contained in a null set and the proposition holds trivially. Hence we assume $L > 0$.
    We assume $\Theta$ is finite. 
    Fix $\theta^\prime \in \Theta$, and set $\varepsilon_k = 2^{-k} \mathrm{diam} (\Theta , d)$.
    
    A set $\mathcal{N} \subset \Theta$ is called an $\varepsilon$-net of $(\Theta , d)$ if the open balls of radius $\varepsilon$ centered at the points of $\mathcal{N}$ cover $\Theta$ (i.e., $\mathcal{N}$ is the set of centers of an $\varepsilon$-cover).
    Let $\kappa$ be the minimal $k \in \mathbb{Z}_{\geq 1}$ such that a minimal $\varepsilon_k$-net of $\Theta$ is $\Theta$ itself. Since $\Theta$ is a metric space with finitely many points, a minimal $\varepsilon_k$-net is $\Theta$ itself once $\varepsilon_k$ is smaller than the minimal distance between two distinct points, so such $k$ exists.
    
    Define the set $\mathcal{N}_0 := \{ \theta^\prime \}$.
    Here we do not require $\mathcal{N}_0$ to be an $\varepsilon_0$-cover; in what follows we only use the fact that, by the definition of $\mathrm{diam} (\Theta , d)$, for any $\theta \in \Theta$ we have $d (\theta , \theta^\prime ) \leq \varepsilon_0 = \mathrm{diam} (\Theta , d)$.
    Let $\{ \mathcal{N}_k \}_{1 \leq k \leq \kappa}$ be a sequence of sets such that each $\mathcal{N}_k$ is a minimal $\varepsilon_k$-net of $\Theta$ and
    \[
        \left| \{
            1 \leq k < \kappa \mid \mathcal{N}_k = \mathcal{N}_{k+1}
        \} \right|
    \]
    is maximized.

    First, we show that there exists no $1 \leq k < \kappa$ such that $\mathcal{N}_k \neq \mathcal{N}_{k+1}$ and $| \mathcal{N}_k | = | \mathcal{N}_{k+1} |$. Suppose such a $k$ existed. Take one such $k$, and take the minimal $1 \leq j \leq k$ such that $\mathcal{N}_j = \mathcal{N}_{j+1} = \cdots = \mathcal{N}_k$ (by the minimality of $j$, either $j = 1$ or $\mathcal{N}_{j-1} \neq \mathcal{N}_j$). Since $\mathcal{N}_{k+1}$ is an $\varepsilon_{k+1}$-net, it is also an $\varepsilon_{i}$-net for each $i = j , \ldots , k$, and since its cardinality $|\mathcal{N}_{k+1}| = |\mathcal{N}_k| = |\mathcal{N}_i|$ is the minimal cardinality, $\mathcal{N}_{k+1}$ is a minimal $\varepsilon_i$-net. Redefining $\mathcal{N}_i := \mathcal{N}_{k+1}$ for $i = j , \ldots , k$, the equalities for the pairs $(i , i+1)$ ($j \leq i < k$) are preserved, the equality for the pair $(k , k+1)$ newly holds, and for the pair $(j-1 , j)$ (when $j \geq 2$) the number of equalities does not decrease because $\mathcal{N}_{j-1} \neq \mathcal{N}_{j}$ before the redefinition. Hence the number of equalities strictly increases, contradicting the maximality of $\{ \mathcal{N}_k \}_{1 \leq k \leq \kappa }$.
    In particular, in the range $1 \leq k < \kappa$, all nets of cardinality $1$ (if any appear) can be taken to be equal.

    For any $\theta \in \Theta$, for $k \geq 1$ let $\pi_k (\theta)$ denote a point of $\mathcal{N}_k$ closest to $\theta$, and set $\pi_0 (\theta) := \theta^\prime$. When $\mathcal{N}_k = \mathcal{N}_{k+1}$, we set $\pi_{k+1} = \pi_k$.
    Then, for any $k \geq 0$, we have $d (\theta , \pi_k (\theta) ) \leq \varepsilon_k$.

    If $\mathcal{N}_k = \mathcal{N}_{k+1}$, then
    \begin{equation}
        \mathbb{P} \left(
            \sup_{\theta \in \Theta} | S_{\pi_k(\theta) } -  S_{\pi_{k+1} (\theta) } | > 0 
            \right)
            = 0 
            .
            \label{eq:pi_k_if_0}
    \end{equation}
    
    Moreover,
    \begin{equation}
        d (\pi_{k} (\theta), \pi_{k+1} (\theta))
        \leq
        d (\pi_{k} (\theta), \theta)
        + d ( \theta, \pi_{k+1} (\theta))
        \leq 3 \varepsilon_{k+1}
        \label{eq:pi_k_pi_k}
    \end{equation}
    holds.
    Since $(S_{\theta})_{\theta \in \Theta}$ is a sub-Gaussian process, by \Cref{prop:subgaussian_1-->4},
    \begin{align}
        & \mathbb{P} \left( \sup_{\theta \in \Theta }  | S_{\pi_{k} (\theta)} - S_{\pi_{k+1} (\theta)} | \geq u \right)
        \notag\\
        & \leq \sum_{ \{ ( \pi_{k} (\theta), \pi_{k+1} (\theta) ) \mid \theta \in \Theta \} } \mathbb{P} \left( | S_{\pi_{k} (\theta)} - S_{\pi_{k+1} (\theta)} | \geq u \right)
        \notag \\
        & \leq \sum_{ \{ ( \pi_{k} (\theta), \pi_{k+1} (\theta) ) \mid \theta \in \Theta \} } 2 \exp \left( - \frac{u^2}{ \| S_{\pi_{k} (\theta)} - S_{\pi_{k+1} (\theta)} \|_{\psi_2}^2} \right) 
        . 
        \notag
    \end{align}
    Since $(S_{\theta})_{\theta \in \Theta}$ is a sub-Gaussian process,
    \begin{align}
        & \sum_{ \{ ( \pi_{k} (\theta), \pi_{k+1} (\theta) ) \mid \theta \in \Theta \} } 2 \exp \left( - \frac{u^2}{\| S_{\pi_{k} (\theta)} - S_{\pi_{k+1} (\theta)} \|_{\psi_2}^2} \right)
        \notag \\
        & \leq \sum_{ \{ ( \pi_{k} (\theta), \pi_{k+1} (\theta) ) \mid \theta \in \Theta \} } 2 \exp \left( - \frac{u^2}{d (\pi_{k} (\theta), \pi_{k+1} (\theta) )^2 L^2} \right)
        .
        \notag
    \end{align}
    By \cref{eq:pi_k_pi_k},
    \begin{align}
        & \sum_{ \{ ( \pi_{k} (\theta), \pi_{k+1} (\theta) ) \mid \theta \in \Theta \} } 2 \exp \left( - \frac{u^2}{d (\pi_{k} (\theta), \pi_{k+1} (\theta) )^2 L^2} \right)
        \notag \\ 
        & \leq \sum_{ \{ ( \pi_{k} (\theta), \pi_{k+1} (\theta) ) \mid \theta \in \Theta \} } 2 \exp \left( - \frac{u^2}{ 9 \varepsilon_{k+1}^2 L^2} \right)
        \notag \\
        & \leq 2 \left| \mathcal{N}_k \right| \left| \mathcal{N}_{k+1} \right| \exp \left( - \frac{u^2}{ 9 \varepsilon_{k+1}^2 L^2} \right)
        \leq
        2 \left| \mathcal{N}_{k+1} \right|^2 \exp \left( - \frac{u^2}{ 9 \varepsilon_{k+1}^2 L^2} \right)
        .
        \notag
    \end{align}
    Define the indicator function
    \begin{align}
        \delta_{\mathcal{N} } (k ) = 
        \begin{cases}
            0, & \mathcal{N}_k = \mathcal{N}_{k+1}, \\
            1, & \mathcal{N}_k \not = \mathcal{N}_{k+1}.
        \end{cases}
    \end{align}
    Then,
    \begin{equation}
        \mathbb{P} \left( \sup_{\theta \in \Theta }  | S_{\pi_{k} (\theta)} - S_{\pi_{k+1} (\theta)} | \geq u \right)
        \leq
        2 \left| \mathcal{N}_{k+1} \right|^2 \exp \left( - \frac{u^2}{ 9 \varepsilon_{k+1}^2 L^2} \right) 
        .
        \label{eq:pi_k_subgaussian_estimate}
    \end{equation}

    On the other hand, since $\pi_0 (\theta ) = \theta^\prime$ and $\pi_{\kappa} (\theta ) = \theta$, we have
    \begin{align}
        \sup_{\theta \in \Theta} | S_{\theta} - S_{\theta^{\prime}} |
        \leq
        \sum_{k=0}^{\kappa-1} \sup_{\theta \in \Theta} | S_{\pi_{k+1} (\theta)} - S_{\pi_{k} (\theta)} |
    \end{align}

    Regarding the inequality, since the covering number $N (\Theta , d , \varepsilon )$ is nonincreasing in $\varepsilon$, for $\varepsilon \leq \varepsilon_{k+1}$ we have $N (\Theta , d , \varepsilon ) \geq N (\Theta , d , \varepsilon_{k+1} ) = |\mathcal{N}_{k+1}|$.
    Moreover, the intervals $[\varepsilon_{k+2} , \varepsilon_{k+1}]$ ($k = 0 , \ldots , \kappa -1$) are pairwise disjoint, and the total of their lengths satisfies $\sum_{k=0}^{\kappa-1} \varepsilon_{k+2} \leq \varepsilon_1 \leq \mathrm{diam} (\Theta , d)$.
    Therefore,
    \begin{align}
        & CL \left( \int_0^\infty \sqrt{\log N (\Theta , d ,\varepsilon ) } d\varepsilon + u \mathrm{diam} (\Theta , d) \right)
        \notag \\
        & \geq CL\sum_{k=0}^{\kappa-1} \int_{\varepsilon_{k+2} }^{\varepsilon_{k+1}} \left( \sqrt{\log N(\Theta ,d ,\varepsilon)}+ u \right) d \varepsilon
        \notag \\
        & \geq CL \sum_{k=0}^{\kappa-1} \varepsilon_{k+2}
        \left(
            \sqrt{\log (|\mathcal{N}_{k+1}|)}+u
        \right)
        \notag \\
        & =
        \frac{C}{2} L \sum_{k=0}^{\kappa-1} \varepsilon_{k+1}
        \left(
            \sqrt{\log (|\mathcal{N}_{k+1}|)}+u
        \right)
    \end{align}
    and hence
    
    \begin{align}
        & \mathbb{P} \left( \sup_{\theta \in \Theta} | S_{\theta} - S_{\theta^\prime} | >
        CL \left( \int_0^\infty \sqrt{\log N (\Theta , d ,\varepsilon ) } d\varepsilon + u \mathrm{diam} (\Theta , d) \right)
        \right)
        \notag \\
        & \leq
        \mathbb{P} \left( \sum_{k=0}^{\kappa-1} \sup_{\theta \in \Theta} | S_{\pi_{k+1} (\theta)} - S_{\pi_{k} (\theta)} |
        >
         \frac{C}{2} L \sum_{k=0}^{\kappa-1} \varepsilon_{k+1}
        \left(
            \sqrt{\log (|\mathcal{N}_{k+1}|)}+u
        \right)
        \right)
        .
    \end{align}
    Here, if $\sum_k a_k > \sum_k b_k$, then there exists $k$ such that $a_k > b_k$. We have
    \begin{align}
        &
        \mathbb{P} \left( \sum_{k=0}^{\kappa-1} \sup_{\theta \in \Theta} | S_{\pi_{k+1} (\theta)} - S_{\pi_{k} (\theta)} |
        >
         \frac{C}{2} L \sum_{k=0}^{\kappa-1} \varepsilon_{k+1}
        \left(
            \sqrt{\log (|\mathcal{N}_{k+1}|)}+u
        \right)
        \right)
        \notag \\
        & \leq
        \mathbb{P} \left( \bigcup_{k=0}^{\kappa-1} \left\{
            \sup_{\theta \in \Theta} | S_{\pi_{k+1} (\theta)} - S_{\pi_{k} (\theta)} |
            >
            \frac{C}{2} L \varepsilon_{k+1}
            \left(
                \sqrt{\log (|\mathcal{N}_{k+1}|)}+u
            \right)
        \right\}
        \right)
        \notag \\
        & \leq
        \sum_{k=0}^{\kappa-1}
        \mathbb{P} \left(
            \sup_{\theta \in \Theta} | S_{\pi_{k+1} (\theta)} - S_{\pi_{k} (\theta)} |
            >
            \frac{C}{2} L \varepsilon_{k+1}
            \left(
                \sqrt{\log (|\mathcal{N}_{k+1}|)}+u
            \right)
        \right)
        .
    \end{align}
    In the last sum, the terms with $\mathcal{N}_k = \mathcal{N}_{k+1}$ (i.e., $\delta_{\mathcal{N}} (k) = 0$) vanish by \cref{eq:pi_k_if_0} (since the threshold is nonnegative, the event with the strict inequality is contained in $\{ \sup_{\theta} | S_{\pi_{k+1} (\theta)} - S_{\pi_{k} (\theta)} | > 0 \}$). The terms with $\delta_{\mathcal{N}} (k) = 1$ can be bounded by the inclusion $\{ \cdot > \cdot \} \subset \{ \cdot \geq \cdot \}$ and \cref{eq:pi_k_subgaussian_estimate}. Therefore,
    \begin{align}
        & \sum_{k=0}^{\kappa-1}
        \mathbb{P} \left(
            \sup_{\theta \in \Theta} | S_{\pi_{k+1} (\theta)} - S_{\pi_{k} (\theta)} |
            >
            \frac{C}{2} L \varepsilon_{k+1}
            \left(
                \sqrt{\log (|\mathcal{N}_{k+1}|)}+u
            \right)
        \right)
        \notag \\
        & \leq
        \sum_{k=0}^{\kappa-1}
                2 \left| \mathcal{N}_{k+1} \right|^2 \exp \left( - \frac{\left(
                    \frac{C}{2} L \varepsilon_{k+1}
                    \left( \sqrt{\log (|\mathcal{N}_{k+1}|)}+u \right)
                \right)^2}{
                    9 \varepsilon_{k+1}^2 L^2
                }
                \right)
            \delta_{\mathcal{N}} (k)
        \notag \\
        & \leq
        \sum_{k=0}^{\kappa-1}
                2 \left| \mathcal{N}_{k+1} \right|^2 \exp \left( - \frac{C^2}{36} \left(
                    \sqrt{\log (|\mathcal{N}_{k+1}|)}+u
                \right)^2
                \right) \delta_{\mathcal{N}} (k)
        \notag \\
        & \leq
        \sum_{k=0}^{\kappa-1}
            2 \left| \mathcal{N}_{k+1} \right|^2 \exp \left( - \frac{C^2}{36} \left(
                \log (|\mathcal{N}_{k+1}|) + u^2
            \right)
                \right)
                \delta_{\mathcal{N}} (k)
        \notag \\
        & \leq
        \sum_{k=0}^{\kappa-1}
            2 \left| \mathcal{N}_{k+1} \right|^{2- C^2 / 36} \delta_{\mathcal{N}} (k) \exp \left( - \frac{C^2}{36} u^2
            \right)
            .
    \end{align}
    
    By the choice of $\{\mathcal{N}_{k} \}_k$, for the indices $k$ with $\delta_{\mathcal{N}} (k) =1$ the cardinalities $\left| \mathcal{N}_{k+1} \right|$ are strictly increasing in $k$, and hence pairwise distinct.
    Moreover, in the range $1 \leq k < \kappa$ there is no change between nets of equal cardinality, so $\delta_{\mathcal{N}} (k) =1$ and $\left| \mathcal{N}_{k+1} \right| = 1$ can occur only for $k=0$.
    Therefore, there exists an injection from the set $\{ \left| \mathcal{N}_{k+1} \right| \, | \, k \in \mathbb{Z}_{\geq 0 }, \delta_{\mathcal{N}} (k) =1 , \left| \mathcal{N}_{k+1} \right| \geq 2 \} $ into $\mathbb{Z}_{\geq 2}$, and
    \begin{align}
        &
        \sum_{k=0}^{\kappa-1}
            2 \left| \mathcal{N}_{k+1} \right|^{2- C^2 / 36} \delta_{\mathcal{N}} (k) \exp \left( - \frac{C^2}{36} u^2
            \right)
        \notag \\
        & \leq
            2 \left( 1 + \sum_{m=2}^{\infty} m^{2- C^2 / 36} \right) \exp \left( - \frac{C^2}{36} u^2
            \right)
        \notag \\
        & =
            2 \zeta \left( \frac{C^2}{36} -2 \right) \exp \left( - \frac{C^2}{36} u^2
            \right) .
    \end{align}
    
    Summing up,
    \begin{align}
        & \mathbb{P} \left( \sup_{\theta \in \Theta} | S_{\theta} - S_{\theta^\prime} | >
        CL \left( \int_0^\infty \sqrt{\log N (\Theta , d ,\varepsilon ) } d\varepsilon + u \mathrm{diam} (\Theta , d) \right)
        \right)
        \notag\\
        & \leq 2 \zeta \left( \frac{C^2}{36} -2 \right) \exp \left( - \frac{C^2}{36} u^2
            \right) .
    \end{align}

    Next, we consider the case where $\Theta$ is a countably infinite set. Suppose $\Theta = \{ \theta^j \mid j \in \mathbb{Z}_{\geq 1} \}$ (by reordering the enumeration, we may assume without loss of generality that $\theta^\prime = \theta^1$).
    For $J \in \mathbb{Z}_{\geq 1}$, define $\Theta^J := \{ \theta^n \mid n = 1,\ldots, J \}$ (note that $\theta^\prime \in \Theta^J$).
    By applying the proposition to $(\Theta^J, d)$, we have
    \begin{align}
        & \mathbb{P} \left( \sup_{\theta \in \Theta^J} | S_{\theta} - S_{\theta^\prime} | >
        CL \left( \int_0^\infty \sqrt{\log N (\Theta^J , d ,\varepsilon ) } d\varepsilon + u\, \mathrm{diam} (\Theta^J , d) \right)
        \right)
        \notag\\
        & \leq 2 \zeta \left( \frac{C^2}{36} -2 \right) \exp \left( - \frac{C^2}{36} u^2
            \right).
    \end{align}
    Since $\mathrm{diam} (\Theta^J , d) \leq \mathrm{diam} (\Theta , d)$,
    we have
    \begin{align}
        & \mathbb{P} \left( \sup_{\theta \in \Theta^J} | S_{\theta} - S_{\theta^\prime} | >
        CL \left( \int_0^\infty \sqrt{\log N (\Theta^J , d ,\varepsilon ) } d\varepsilon + u \mathrm{diam} (\Theta , d) \right)
        \right)
        \notag \\
        & \leq \mathbb{P} \left( \sup_{\theta \in \Theta^J} | S_{\theta} - S_{\theta^\prime} |
        -
        CL \int_0^\infty \sqrt{\log N (\Theta^J , d ,\varepsilon ) } d\varepsilon
        >
        CL
        u \mathrm{diam} (\Theta , d) \right)
        \notag\\
        & \leq 2 \zeta \left( \frac{C^2}{36} -2 \right) \exp \left( - \frac{C^2}{36} u^2
            \right)
        .
    \end{align}
    When $u = 0$, the assertion is trivial because $2 \zeta \left( \frac{C^2}{36} - 2 \right) > 1$; hence we may assume $u > 0$. Moreover, if $\mathrm{diam} (\Theta , d) = 0$, then $\Theta$ would be a singleton since $(\Theta , d)$ is a metric space, contradicting that $\Theta$ is countably infinite; hence $\mathrm{diam} (\Theta , d) > 0$.
    Take an arbitrary constant $\delta \in \left( 0 , CL u \, \mathrm{diam} (\Theta , d) \right)$, and define the sequence of events
    \[
        A_J (\delta) := \left\{
            \omega \in \Omega
        \middle|
        \begin{array}{l}
            \displaystyle
            \sup_{\theta \in \Theta^J} | S_{\theta} - S_{\theta^\prime} | -
            CL \int_0^\infty \sqrt{\log N (\Theta^J , d ,\varepsilon ) } d\varepsilon
            \\
            \displaystyle
            > CL u \, \mathrm{diam} (\Theta , d) - \delta
        \end{array}
        \right\} .
    \]
    On the event
    \[
        A_\infty := \left\{
            \omega \in \Omega
        \middle|
        \begin{array}{l}
            \displaystyle
            \liminf_{J \to \infty} \sup_{\theta \in \Theta^J} | S_{\theta} - S_{\theta^\prime} |
            -
            \limsup_{J \to \infty}
            CL \int_0^\infty \sqrt{\log N (\Theta^J , d ,\varepsilon ) } d\varepsilon
            \\
            \displaystyle
            > CL u \, \mathrm{diam} (\Theta , d)
        \end{array}
        \right\} ,
    \]
    since $\sup_{\theta \in \Theta^J} | S_{\theta} - S_{\theta^\prime} |$ is nondecreasing in $J$, for sufficiently large $J$ this quantity is at least $\liminf_{J^\prime \to \infty} \sup_{\theta \in \Theta^{J^\prime}} | S_{\theta} - S_{\theta^\prime} | - \delta / 2$, and, for sufficiently large $J$, we also have $CL \int_0^\infty \sqrt{\log N (\Theta^J , d ,\varepsilon ) } d\varepsilon \leq \limsup_{J^\prime \to \infty} CL \int_0^\infty \sqrt{\log N (\Theta^{J^\prime} , d ,\varepsilon ) } d\varepsilon + \delta / 2$.
    Therefore, $A_\infty \subset \liminf_{J \to \infty} A_J (\delta)$.
    By Fatou's lemma (in the form for events) and the above bound for each $J$, applied with $u^\prime := u - \delta / (CL \, \mathrm{diam} (\Theta , d)) \geq 0$ in place of $u$, we obtain
    \begin{align}
        & \mathbb{P} \left( A_\infty \right)
        \leq
        \mathbb{P} \left( \liminf_{J \to \infty} A_J (\delta) \right)
        \leq
        \liminf_{J \to \infty} \mathbb{P} \left( A_J (\delta) \right)
        \notag \\
        & \leq
        2 \zeta \left( \frac{C^2}{36} -2 \right) \exp \left( - \frac{C^2}{36} \left( u - \frac{\delta}{CL \, \mathrm{diam} (\Theta , d)} \right)^2 \right)
        .
        \notag
    \end{align}
    Since the right-hand side is continuous in $\delta$, letting $\delta \to 0$ yields
    \begin{align}
        \mathbb{P} \left( A_\infty \right)
        \leq
        2 \zeta \left( \frac{C^2}{36} -2 \right) \exp \left( - \frac{C^2}{36} u^2 \right)
        .
        \label{eq:dudley_integral_ineq_prob_030}
    \end{align}%
    
    The sequence of random variables
    \[
        \sup_{\theta \in \Theta^J} | S_{\theta} - S_{\theta^\prime} |
    \]
    is nondecreasing in $J$ and converges, as $J \to \infty$, to
    \[
        \sup_{\theta \in \Theta} | S_{\theta} - S_{\theta^\prime} |.
    \]
    Therefore,
    \begin{equation}
        \liminf_{J \to \infty} \sup_{\theta \in \Theta^J} | S_{\theta} - S_{\theta^\prime} | 
        =
        \sup_{\theta \in \Theta} | S_{\theta} - S_{\theta^\prime} |
        .
        \label{eq:convergence_sup_S}
    \end{equation}

    From \Cref{prop:covereing_leq_packing,prop:packing_ineq}
    \begin{align}
        \sqrt{\log N (\Theta^J , d ,\varepsilon ) }
        \leq
        \sqrt{\log P (\Theta^J , d ,\varepsilon ) }
        \leq
        \sqrt{\log P (\Theta , d ,\varepsilon ) }
        \leq
        \sqrt{\log N (\Theta , d ,\varepsilon /2 ) }.
    \end{align}
    Since $ \sqrt{\log N (\Theta , d ,\varepsilon /2 ) }$ is integrable over $(0, \infty)$ with respect to $\varepsilon$, it follows from the reverse Fatou's inequality that
    \begin{align}
        \limsup_{J \to \infty}
        \int_0^\infty \sqrt{\log N (\Theta^J , d ,\varepsilon ) } d\varepsilon
        & \leq \int_0^\infty \sqrt{\log \left( \limsup_{J \to \infty} N (\Theta^J , d ,\varepsilon ) \right)} d\varepsilon
        .
    \end{align}
    Here, since for sufficiently large $J$, $\Theta^J$ contains an $\varepsilon$-net of $(\Theta, d)$, we have
    \[
        \limsup_{J \to \infty} N (\Theta^J , d ,\varepsilon )
        \leq N (\Theta , d ,\varepsilon ).
    \]
    Summarizing the above, we obtain
    \begin{equation}
        \limsup_{J \to \infty}
        \int_0^\infty \sqrt{\log N (\Theta^J , d ,\varepsilon ) }\, d\varepsilon
        \leq 
        \int_0^\infty \sqrt{\log N (\Theta , d ,\varepsilon ) }\, d\varepsilon
        \label{eq:covering_Theta_J_leq_Theta}
    \end{equation}

    From \cref{eq:dudley_integral_ineq_prob_030,eq:convergence_sup_S,eq:covering_Theta_J_leq_Theta}
    \begin{align}
        & 2 \zeta \left( \frac{C^2}{36} -2 \right) \exp \left( - \frac{C^2}{36} u^2
        \right)
        \notag \\
        & \geq \mathbb{P} \left( \left\{
            \omega \in \Omega
        \middle|
            \sup_{\theta \in \Theta} | S_{\theta} - S_{\theta^\prime} |
            >
            CL \left( \int_0^\infty \sqrt{\log N (\Theta , d ,\varepsilon ) } d\varepsilon
            + u \mathrm{diam} (\Theta , d) \right)
        \right\}
        \right)
        .
    \end{align}

    Finally, we consider the general case where $(\Theta, d)$ is an arbitrary metric space.
    By the separability of the process $(S_{\theta})_{\theta \in \Theta}$ (\Cref{def:separable_process}), there exists an at most countable dense subset $\Theta' \subset \Theta$ such that, almost surely,
    $\sup_{\theta \in \Theta} (S_{\theta} - S_{\theta^\prime}) = \sup_{\theta \in \Theta'} (S_{\theta} - S_{\theta^\prime})$;
    moreover, for any $\varepsilon_0 > 0$,
    $N(\Theta', d, \varepsilon + \varepsilon_0) \leq N(\Theta, d, \varepsilon)$,
    and $\mathrm{diam} (\Theta' , d) = \mathrm{diam} (\Theta , d)$.
    Replacing $\Theta'$ by $\Theta' \cup \{ \theta^\prime \}$ if necessary (which preserves countability, density, and the three properties above), we may assume $\theta^\prime \in \Theta'$.
    Since $N(\Theta, d, \varepsilon)$ is nonincreasing in $\varepsilon$ and takes values in $\mathbb{Z}_{\geq 1}$,
    for almost every $\varepsilon$, we have
    \begin{align}
        N(\Theta', d, \varepsilon)
        = \lim_{\varepsilon_0 \to 0,\, \varepsilon_0 > 0} N(\Theta', d, \varepsilon + \varepsilon_0)
        \leq N(\Theta, d, \varepsilon).
    \end{align}
    Summarizing the above, we conclude that the proposition also holds for general $(\Theta, d)$.
\end{proof}
}

\begin{proposition}
    \label{prop:dudley_integral_ineq_prob_pseudo}
    \Cref{prop:dudley_integral_ineq_prob} also holds for a separable sub-Gaussian process on a pseudometric space $(\Theta, d)$.
\end{proposition}

\begin{proof}
    We have
    \[
        \sup_{\theta \in \Theta } S_{\theta} = \sup_{ [ \theta ] \in \Theta^* } S_{ [ \theta ] } ,
        \quad
        N( \Theta, d, \varepsilon ) = N( \Theta^*, d^*, \varepsilon ) ,
        \quad
        \mathrm{diam} (\Theta , d) = \mathrm{diam} (\Theta^* , d^*)
        .
    \]
    Therefore, the proposition follows from \Cref{prop:dudley_integral_ineq_prob}.
\end{proof}

\subsection{Statistical Learning Theory for Sub-Gaussian Random Variables}

\begin{theorem}
    \label{theo:learning_theory_expectation}
    We assume that the loss function $\ell \colon \Theta \times \mathcal{S} \to \mathbb{R}$ satisfies \cref{eq:assu:loss_estimated_pseudometric}.
    Assume that the process $( - X_{\theta} )_{\theta \in \Theta}$ is separable on the pseudometric space $(\Theta , d_{\mathcal{S}})$ (\Cref{def:separable_process}).
    Then,
    \begin{align}
        & \mathbb{E}_{S^{(1)}, \ldots , S^{(N)} } \mathbb{E}_{S} \ell (\theta^{*(N)} , S )
        - \mathbb{E}_{S} \ell (\theta^{*} , S )
        \notag \\
        & \leq
        \frac{8 \sqrt{2 (1 + \log 2)} }{\sqrt{N}}
            \int_0^{\infty} \sqrt{\log N( \Theta, d_{\mathcal{S}}, \varepsilon ) } d \varepsilon
         .
         \notag
    \end{align}
\end{theorem}

\begin{proof}
    Let $X_{\theta}$ be the centered empirical process defined at the beginning of \Cref{sec:apx_learning_theory}.
    First,
    \begin{align}
        & \mathbb{E}_S \ell (\theta^{*(N)} , S  ) - \mathbb{E}_S \ell (\theta^{*} , S  )
        \notag \\
        & \leq
        \left( \mathbb{E}_S \ell (\theta^{*(N)} , S  )
        - \frac{1}{N} \sum_{n=1}^N \ell (\theta^{*(N)} , S^{(n)}  ) \right)
        \notag \\
        & \quad + \left(
        \frac{1}{N} \sum_{n=1}^N \ell (\theta^{*(N)} , S^{(n)}  )
        - \frac{1}{N} \sum_{n=1}^N \ell (\theta^{*} , S^{(n)}  )
        \right)
        \notag \\
        & \quad +
        \left(
            \frac{1}{N} \sum_{n=1}^N \ell (\theta^{*} , S^{(n)}  )
            - \mathbb{E}_S \ell (\theta^{*} , S  )
        \right) ,
        \notag
    \end{align}
    and the second term is $\leq 0$ by the definition of $\theta^{*(N)}$. Hence
    \begin{align}
        \mathbb{E}_S \ell (\theta^{*(N)} , S  ) - \mathbb{E}_S \ell (\theta^{*} , S  )
        & \leq - X_{\theta^{*(N)}} + X_{\theta^{*}}
        \notag \\
        & \leq \sup_{\theta \in \Theta } ( X_{\theta^{*}} - X_{\theta} )
        \leq X_{\theta^{*}} + \sup_{\theta \in \Theta } ( - X_{\theta} ) .
        \label{eq:lemma_emperical_estimate}
    \end{align}
    Since $\mathbb{E} X_{\theta^{*}} = 0$, taking expectations yields
    \[
        \mathbb{E}_{S^{(1)},\ldots,S^{(N)}} \mathbb{E}_S \ell(\theta^{*(N)}, S) - \mathbb{E}_S \ell(\theta^*, S)
        \leq \mathbb{E}_{S^{(1)},\ldots,S^{(N)}} \sup_{\theta \in \Theta } ( - X_{\theta} ) .
    \]
    We now bound the right-hand side.

    For any $\theta , \theta^\prime \in \Theta$, write
    \[
        X_{\theta} - X_{\theta^\prime}
        = \frac{1}{N} \sum_{n=1}^N Z^{(n)},
        \quad
        Z^{(n)} := \left( \ell (\theta , S^{(n)}) - \ell ( \theta^\prime , S^{(n)}) \right) - \mathbb{E}_S \left[ \ell ( \theta , S) - \ell ( \theta^\prime , S) \right] .
    \]
    The random variables $Z^{(1)} , \ldots , Z^{(N)}$ are independent and have mean $0$.
    By \Cref{prop:centering_subGaussian} and \cref{eq:assu:loss_estimated_pseudometric}, for each $n$,
    \begin{equation}
        \| Z^{(n)} \|_{\psi_2}
        \leq 2 \left\| \ell ( \theta , S^{(n)}) - \ell ( \theta^\prime , S^{(n)}) \right\|_{\psi_2}
        \leq 2 \left\| d_{S^{(n)}} ( \theta , \theta^\prime ) \right\|_{\psi_2}
        = 2 d_{\mathcal{S}} (\theta , \theta^\prime )
        \label{eq:checked_exp_subgaussian_process}
    \end{equation}
    (the last equality follows from $S^{(n)} \sim S$).
    By \Cref{prop:subgaussian_4-->5}, for each $n$ and any $\lambda \in \mathbb{R}$,
    $
        \mathbb{E} \exp ( \lambda Z^{(n)} / N )
        \leq \exp ( 2 (1 + \log 2) \lambda^2 d_{\mathcal{S}} (\theta , \theta^\prime )^2 / N^2 )
    $.
    By the independence of $Z^{(1)} , \ldots , Z^{(N)}$, for any $\lambda \in \mathbb{R}$,
    \begin{align*}
        \mathbb{E} \exp \left( \lambda ( X_{\theta} - X_{\theta^\prime} ) \right)
        & = \prod_{n=1}^N \mathbb{E} \exp \left( \frac{\lambda}{N} Z^{(n)} \right)
        \\
        & \leq \exp \left( \frac{\lambda^2}{2} \left( \frac{2 \sqrt{1 + \log 2}}{\sqrt{N}} \right)^2 d_{\mathcal{S}} (\theta , \theta^\prime )^2 \right) .
    \end{align*}
    Since this bound holds for any $\lambda \in \mathbb{R}$, the same bound holds for the sign-reversed increments $(-X_{\theta}) - (-X_{\theta^\prime})$. That is, the process $(-X_{\theta})_{\theta \in \Theta}$ satisfies \cref{eq:subgaussian_process_MGF} on $(\Theta , d_{\mathcal{S}})$ with $L = 2 \sqrt{1 + \log 2} / \sqrt{N}$.
    Since $\mathbb{E} ( - X_{\theta} ) = 0$ for any $\theta \in \Theta$ and the process $(-X_{\theta})_{\theta \in \Theta}$ is separable by assumption, \Cref{prop:dudley_integral_ineq_E_pseudo} yields
    \begin{align*}
        \mathbb{E}_{S^{(1)},\ldots,S^{(N)}} \sup_{\theta \in \Theta} ( - X_{\theta} )
        & \leq 4 \sqrt{2} \cdot \frac{2 \sqrt{1 + \log 2}}{\sqrt{N}} \int_0^\infty \sqrt{\log N(\Theta, d_{\mathcal{S}}, \varepsilon)}\, d\varepsilon
        \\
        & = \frac{8 \sqrt{2 (1 + \log 2)}}{\sqrt{N}} \int_0^\infty \sqrt{\log N(\Theta, d_{\mathcal{S}}, \varepsilon)}\, d\varepsilon .
    \end{align*}
\end{proof}

\begin{theorem}
    \label{theo:learning_theory_prob}
    We assume that the loss function $\ell \colon \Theta \times \mathcal{S} \to \mathbb{R}$ satisfies \cref{eq:assu:loss_estimated_pseudometric}.
    Assume that the process $( - X_{\theta} )_{\theta \in \Theta}$ is separable on the pseudometric space $(\Theta , d_{\mathcal{S}})$ (\Cref{def:separable_process}).
    Let $C > 6 \sqrt{3}$.
    Then, for any $u \geq 0$,
    \begin{align}
        &
        \mathbb{P}
        \left(
        \begin{array}{l}
            \displaystyle \mathbb{E}_S \ell (\theta^{*(N)} , S )
            - \mathbb{E}_S \ell (\theta^{*} , S )    \\
            \displaystyle
            >
            \frac{ 2 \sqrt{6 (1 + \log 2)}\, C}{\sqrt{N}}
            \left(
                \int_0^{\infty} \sqrt{\log N( \Theta , d_{\mathcal{S}}, \varepsilon ) } d \varepsilon + u \mathrm{diam} (\Theta , d_{\mathcal{S}})
            \right)
        \end{array}
        \right)
        \notag \\
        & \leq 2 \zeta \left( \frac{C^2}{36} -2 \right) \exp \left( -\frac{C^2}{36} u^2  \right)
        .
        \notag
    \end{align}
    In particular, for $C = 8 \sqrt{2}$,
        \begin{align}
        &
        \mathbb{P}
        \left(
        \begin{array}{l}
            \displaystyle \mathbb{E}_S \ell (\theta^{*(N)} , S )
            - \mathbb{E}_S \ell (\theta^{*} , S )
            \\
            \displaystyle >
            \frac{ 32 \sqrt{3 (1 + \log 2)} }{\sqrt{N}}
            \left(
                \int_0^{\infty} \sqrt{\log N( \Theta , d_{\mathcal{S}} , \varepsilon ) } d \varepsilon + u \mathrm{diam} (\Theta , d_{\mathcal{S}})
            \right)
        \end{array}
        \right)
        \notag \\
        & \leq 5 \exp \left( -3 u^2  \right)
        ,
        \notag
    \end{align}
    and for $C = 12$,
        \begin{align}
        &
        \mathbb{P}
        \left(
        \begin{array}{l}
            \displaystyle \mathbb{E}_S \ell (\theta^{*(N)} , S )
            - \mathbb{E}_S \ell (\theta^{*} , S )
            \\
            \displaystyle >
            \frac{ 24 \sqrt{6 (1 + \log 2)} }{\sqrt{N}}
            \left(
                \int_0^{\infty} \sqrt{\log N( \Theta , d_{\mathcal{S}}, \varepsilon ) } d \varepsilon + u \mathrm{diam} (\Theta , d_{\mathcal{S}})
            \right)
        \end{array}
        \right)
        \notag \\
        & \leq 3.3 \exp \left( -4 u^2  \right)
        .
        \notag
    \end{align}
\end{theorem}

\begin{proof}
    Define $X_{\theta}$ as in \Cref{theo:learning_theory_expectation}.
    By \cref{eq:lemma_emperical_estimate},
    $\mathbb{E}_S \ell (\theta^{*(N)} , S ) - \mathbb{E}_S \ell (\theta^{*} , S ) \leq \sup_{\theta \in \Theta } ( X_{\theta^{*}} - X_{\theta} )$
    holds.
    Moreover, with $Z^{(n)}$ as in the proof of \Cref{theo:learning_theory_expectation}, \Cref{prop:meanzero_subGaussian_sum} and \cref{eq:checked_exp_subgaussian_process} yield, for any $\theta , \theta^\prime \in \Theta$,
    \begin{align*}
        \| X_{\theta} - X_{\theta^\prime} \|_{\psi_2}
        & = \frac{1}{N} \left\| \sum_{n=1}^N Z^{(n)} \right\|_{\psi_2}
        \leq \frac{\sqrt{6 (1 + \log 2)}}{N} \left( \sum_{n=1}^N \| Z^{(n)} \|_{\psi_2}^2 \right)^{1/2}
        \\
        & \leq \frac{2 \sqrt{6 (1 + \log 2)}}{\sqrt{N}}\, d_{\mathcal{S}} ( \theta , \theta^\prime ) .
    \end{align*}
    That is, the process $(-X_{\theta})_{\theta \in \Theta}$ is a sub-Gaussian process on the pseudometric space $(\Theta , d_{\mathcal{S}})$, and the constant in \cref{eq:subgaussian_process} can be taken as $L = 2 \sqrt{6 (1 + \log 2)} / \sqrt{N}$.
    Since the process $(-X_{\theta})_{\theta \in \Theta}$ is separable by assumption,
    applying \Cref{prop:dudley_integral_ineq_prob_pseudo} to $( - X_{\theta} )_{\theta \in \Theta}$ with $\theta^\prime = \theta^*$ (note that $\sup_{\theta \in \Theta} ( (-X_{\theta}) - (-X_{\theta^*}) ) = \sup_{\theta \in \Theta} ( X_{\theta^*} - X_{\theta} )$), we obtain, for any $u \geq 0$,
    \begin{align*}
        &
        \mathbb{P}
        \Biggl[
        \sup_{\theta \in \Theta } ( X_{\theta^*} - X_{\theta} ) >
        \frac{ 2 \sqrt{6 (1 + \log 2)}\, C }{\sqrt{N}} \biggl(
            \int_0^{\infty} \sqrt{\log N( \Theta, d_{\mathcal{S}}, \varepsilon ) } d \varepsilon
        \\
        & \qquad \qquad
            + u \mathrm{diam} (\Theta , d_{\mathcal{S}})
        \biggr)
        \Biggr]
        \leq 2 \zeta \left( \frac{C^2}{36} -2 \right) \exp \left( -\frac{C^2}{36} u^2  \right)
        .
    \end{align*}
    Combining this with \cref{eq:lemma_emperical_estimate}, we obtain the claim for general $C$.
    For $C = 8\sqrt{2}$, we have $2 \sqrt{6 (1 + \log 2)}\, C = 32 \sqrt{3 (1 + \log 2)}$, and the tail is bounded by $5 \exp (-3u^2)$ by the case $C = 8\sqrt{2}$ of \Cref{prop:dudley_integral_ineq_prob}.
    For $C = 12$, we have $2 \sqrt{6 (1 + \log 2)}\, C = 24 \sqrt{6 (1 + \log 2)}$, and similarly the tail is bounded by $3.3 \exp (-4 u^2)$.
\end{proof}

\begin{remark}[Explicit constants]
    \Cref{theo:learning_theory_expectation,theo:learning_theory_prob} rely on a Dudley-type tail inequality for sub-Gaussian processes on pseudo-metric spaces (\Cref{prop:dudley_integral_ineq_prob}; for the extension to pseudo-metric spaces, see \Cref{prop:dudley_integral_ineq_prob_pseudo}). This inequality reduces the constant $C = 6(1 + 2/\log 2) \approx 23.3$ of \citet[Theorem 5.29]{van2014probability} to any constant larger than $6\sqrt{3} \approx 10.4$, and sharpens the tail decay from $2\exp(-u^2)$ to $2\zeta(C^2/36 - 2)\exp(-C^2 u^2 / 36)$ (e.g., $5\exp(-3u^2)$ for $C = 8\sqrt{2}$, and $3.3\exp(-4u^2)$ for $C = 12$). This inequality applies to general separable sub-Gaussian processes beyond inverse optimization (a tail inequality of the same type is itself known; \citealp[Theorem 8.1.6]{vershynin2020high}).
\end{remark}

\subsection{Covering Number and Diameter in Product Spaces}

\begin{proposition}
    \label{prop:covering_number_scaling}
    Let $(\Theta , d)$ be a pseudometric space, and let $L > 0$, $\varepsilon > 0$ be constants. Then,
    \[
        N (\Theta , L d ,\varepsilon )
        =
        N (\Theta , d ,\varepsilon / L )
        .
    \]
\end{proposition}

\begin{proof}
    Let $\mathcal{N}$ be an $\varepsilon$-cover of $\Theta$ with respect to the metric $L d$, realizing the covering number $N (\Theta , L d , \varepsilon )$. By definition, for any $\mathcal{N}_i \in \mathcal{N}$, there exists $\theta^i \in \Theta$ such that
    \[
        \mathcal{N}_i = \{ \theta \in \Theta \mid L d (\theta , \theta^i ) < \varepsilon \}
        = \{ \theta \in \Theta \mid d (\theta , \theta^i ) < \varepsilon / L \}
    \]
    holds. Therefore, $\mathcal{N}$ is an $\varepsilon / L$-cover of $(\Theta, d)$. This implies
    \[
        N (\Theta , L d ,\varepsilon )
        \geq
        N (\Theta , d ,\varepsilon / L )
    \]
    holds. Similarly, the reverse inequality can be shown.
\end{proof}

Let $(\Theta_1, d_1)$ and $(\Theta_2, d_2)$ be pseudometric spaces. Define the pseudometric $d_{12}$ on $\Theta_1 \times \Theta_2$ by
\[
    d_{12} ((\theta_{1}, \theta_{2}), (\theta_{1}^\prime, \theta_{2}^\prime))
    :=
    d_1 (\theta_1, \theta_1^\prime)
    +
    d_2 (\theta_2, \theta_2^\prime)
    .
\]

\begin{proposition}
    \label{prop:covering_integral_on_product_sp}
    Let $p_1, p_2 > 0$ be constants such that $p_1 + p_2 = 1$. Then,
    \[
        N (\Theta_1 \times \Theta_2 , d_{12} ,\varepsilon )
        \leq
        N (\Theta_1 , d_{1} , p_1 \varepsilon )
        N (\Theta_2 , d_{2} , p_2 \varepsilon )
        .
    \]
\end{proposition}

\begin{proof}
    Let $\mathcal{N}^j$ be a $p_j \varepsilon$-cover of $\Theta_j$ that realizes the covering number $N (\Theta_j , d_{j}, p_j \varepsilon )$. For any $\mathcal{N}_i^j \in \mathcal{N}^j$, there exists $\theta^{ij} \in \Theta_j$ such that $\mathcal{N}_i^j = \{ \theta^{j} \in \Theta_j \mid d_j (\theta^j , \theta^{ij} ) < p_j \varepsilon\}$.

    Consider the Cartesian product $\mathcal{N}_{i_1}^1 \times \mathcal{N}_{i_2}^2$, then
    \begin{align}
        \mathcal{N}_{i_1}^1 \times \mathcal{N}_{i_2}^2
        & = \{ \theta^{1} \in \Theta_1 \mid d_1 (\theta^1 , \theta^{i_1 1} ) < p_1 \varepsilon\} \times 
        \{ \theta^{2} \in \Theta_2 \mid d_2 (\theta^2 , \theta^{i_2 2} ) < p_2 \varepsilon\}
        \notag \\
        & \quad \subset
        \{ ( \theta^{1} , \theta^{2} ) \in \Theta_1 \times \Theta_2 \mid d_{12} ( ( \theta^1 , \theta^2 ) , ( \theta^{i_1 1} , \theta^{i_2 2} ) ) < \varepsilon \}
        \notag
    \end{align}
    holds. The collection $\{ \{ ( \theta^{1} , \theta^{2} ) \in \Theta_1 \times \Theta_2 \mid d_{12} (( \theta^1 ,  \theta^2 ) , ( \theta^{i_1 1} , \theta^{i_2 2} ) ) < \varepsilon \} \}_{i_1, i_2}$ forms an $\varepsilon$-cover of $\Theta_1 \times \Theta_2$. Therefore, the proposition holds.
\end{proof}

\begin{proposition}
    \label{prop:covering_integral_on_product_sp_2}
    Let $p_1, p_2 > 0$ be constants such that $p_1 + p_2 = 1$. Then,
    \begin{align*}
        & \int_0^{\infty} \sqrt{ \log N (\Theta_1 \times \Theta_2, d_{12}, \varepsilon ) } \, d \varepsilon
        \\
        & \leq
        \frac{1}{p_1}
        \int_0^{\infty} \sqrt{ \log N (\Theta_1, d_{1}, \varepsilon ) } \, d \varepsilon
        +
        \frac{1}{p_2}
        \int_0^{\infty} \sqrt{ \log N (\Theta_2, d_{2}, \varepsilon ) } \, d \varepsilon
        .
    \end{align*}
\end{proposition}

\begin{proof}
    Taking the logarithm of both sides of \Cref{prop:covering_integral_on_product_sp}, we have
    \[
        \log N (\Theta_1 \times \Theta_2 , d_{12} ,\varepsilon )
        \leq
        \log N (\Theta_1 , d_{1} , p_1 \varepsilon )
        + \log N (\Theta_2 , d_{2} , p_2 \varepsilon )
        .
    \]
    In general, since $\sqrt{\varepsilon_1 + \varepsilon_2} \leq \sqrt{\varepsilon_1} + \sqrt{\varepsilon_2}$, it follows that
    \[
        \sqrt{\log N (\Theta_1 \times \Theta_2 , d_{12} ,\varepsilon )}
        \leq
        \sqrt{
        \log N (\Theta_1 , d_{1} , p_1 \varepsilon )}
        + \sqrt{\log N (\Theta_2 , d_{2} , p_2 \varepsilon )}
        .
    \]
    Integrating both sides over $[0, \infty)$ and applying \Cref{prop:covering_number_scaling} yields the proposition.
\end{proof}

Let $B^D := \{ x \in \mathbb{R}^D \mid \| x \|_2 \leq 1 \}$ denote the $D$-dimensional unit ball.

\begin{proposition}[{Cf.\ \citet[Proposition 4.2.13]{vershynin2020high}}]
    \label{prop:covering_number_unit_ball}
    For a subset $\Theta \subset B^D$ and $\varepsilon > 0$,
    \[
        N ( \Theta, \| \bullet \|_2 , \varepsilon )
        \leq \left( \frac{2}{\varepsilon} + 1 \right)^D
        .
    \]
\end{proposition}

\begin{proof}
    Let $A \subset \Theta$ be an $\varepsilon$-separated subset. The open balls of radius $\varepsilon / 2$ centered at the points of $A$ are pairwise disjoint and contained in $(1 + \varepsilon / 2) B^D$. Comparing volumes, we obtain $| A | \leq ( (1 + \varepsilon / 2) / (\varepsilon / 2) )^D = ( 2 / \varepsilon + 1 )^D$. The proposition follows from \Cref{prop:covereing_leq_packing}.
\end{proof}

\begin{proposition}
    \label{prop:covering_integral_unit_ball}
    For a subset $\Theta \subset B^D$,
    \[
        \int_0^{\infty}
        \sqrt{ \log N ( \Theta , \| \bullet \|_2 , \varepsilon ) }
        d \varepsilon
        \leq \sqrt{D} \int_0^2 \sqrt{ \log \left( \frac{2}{\varepsilon} + 1 \right) } d \varepsilon
        \leq 2 \sqrt{2 \log 2} \sqrt{D}
        .
    \]
\end{proposition}

\begin{proof}
    If $\varepsilon > 2 \geq \mathrm{diam} (\Theta)$, the open ball of radius $\varepsilon$ centered at any point of $\Theta$ covers $\Theta$, so $N ( \Theta, \| \bullet \|_2 , \varepsilon ) = 1$. Hence
    \[
        \int_0^{\infty}
        \sqrt{\log N ( \Theta, \| \bullet \|_2 , \varepsilon )}
        d \varepsilon
        =
        \int_0^{2}
        \sqrt{\log N ( \Theta, \| \bullet \|_2 , \varepsilon )}
        d \varepsilon
        .
    \]
    By \Cref{prop:covering_number_unit_ball},
    \[
        \int_0^{2}
        \sqrt{\log N ( \Theta , \| \bullet \|_2 , \varepsilon )}
        d \varepsilon
        \leq
        \sqrt{D} \int_0^2 \sqrt{\log \left( \frac{2}{\varepsilon} + 1 \right)}
        d \varepsilon
        .
    \]
    Here, by the change of variables $\varepsilon = 2 t$, the Cauchy--Schwarz inequality, and $\int_0^1 \log ( 1 / t + 1 ) dt = 2 \log 2$,
    \begin{align*}
        \int_0^2 \sqrt{\log \left( \frac{2}{\varepsilon} + 1 \right)} d \varepsilon
        & = 2 \int_0^1 \sqrt{\log \left( \frac{1}{t} + 1 \right)} d t
        \\
        & \leq 2 \left( \int_0^1 \log \left( \frac{1}{t} + 1 \right) d t \right)^{1/2}
        = 2 \sqrt{2 \log 2}
        .
    \end{align*}
    Thus, the proposition follows.
\end{proof}

\begin{proposition}
    \label{prop:diam_of_scaling}
    Let $(\Theta, d)$ be a pseudometric space, and let $L > 0$ be a constant. Then,
    \[
        \mathrm{diam} (\Theta, L d)
        =
        L \, \mathrm{diam} (\Theta, d)
        .
    \]
\end{proposition}

\begin{proof}
    \begin{align*}
        \mathrm{diam} (\Theta, L d)
        &= \sup_{\theta^1, \theta^2 \in \Theta} L d (\theta^1, \theta^2) 
        = L \sup_{\theta^1, \theta^2 \in \Theta} d (\theta^1, \theta^2) 
        = L \, \mathrm{diam} (\Theta, d)
        .
    \end{align*}
\end{proof}

\begin{proposition}
    \label{prop:diam_on_product_sp}
    \[
        \mathrm{diam} (\Theta_1 \times \Theta_2, d_{12})
        =
        \mathrm{diam} (\Theta_1, d_{1})
        +
        \mathrm{diam} (\Theta_2, d_{2})
        .
    \]
\end{proposition}

\begin{proof}
    \begin{align*}
        & \mathrm{diam} (\Theta_1 \times \Theta_2, d_{12}) = \sup_{\substack{
            (\theta^1, \theta^2), \\
            (\theta^{1\prime}, \theta^{2\prime}) \in \Theta_1 \times \Theta_2
        }}
        \left(
            d_1 (\theta^1, \theta^{1\prime}) + d_2 (\theta^2, \theta^{2\prime})
        \right) \\
        &= \sup_{\theta^1, \theta^{1\prime} \in \Theta_1} d_1 (\theta^1, \theta^{1\prime})
        + \sup_{\theta^2, \theta^{2\prime} \in \Theta_2} d_2 (\theta^2, \theta^{2\prime})
        = \mathrm{diam} (\Theta_1, d_{1})
        + \mathrm{diam} (\Theta_2, d_{2})
        .
    \end{align*}
\end{proof}

\subsection{Hausdorff Distance}
\label{sec:Hausdorff_metric}

For nonempty compact sets $\mathcal{M}_1, \mathcal{M}_2$ in Euclidean space, the Hausdorff distance is defined as
\begin{align*}
    d^H (\mathcal{M}_1, \mathcal{M}_2 )
    = \max \left( \sup_{x^1 \in \mathcal{M}_1} \inf_{x^2 \in \mathcal{M}_2} \| x^1 - x^2 \|_2 ,
    \sup_{x^2 \in \mathcal{M}_2} \inf_{x^1 \in \mathcal{M}_1} \| x^1 - x^2 \|_2  
    \right)
    .
\end{align*}
The convex hull of a compact set $\mathcal{M}$ in Euclidean space is denoted by $\mathrm{Conv} \mathcal{M}$.

\begin{proposition}[{\citealp[Cf.][Lemma 1.8.14]{schneider2014convex}}]
    \label{prop:Hausdorff_distance_is_diff_LP}
    Let $\mathcal{M}_1, \mathcal{M}_2$ be compact convex sets in Euclidean space. Then,
    \[
        d^H (\mathcal{M}_1, \mathcal{M}_2 )
        = 
        \max_{\| \theta \|_2 =1} 
            \left| 
            \max_{x \in \mathcal{M}_1} \theta^\top x  - 
            \max_{x \in \mathcal{M}_2} \theta^\top x 
            \right|
        .
    \]
\end{proposition}

\begin{proposition}[{\citealp[Cf.][p.\ 64]{schneider2014convex}}]
    \label{prop:convex_hull_is_Lipschitz}
    For any compact sets $\mathcal{M}_1, \mathcal{M}_2$ in Euclidean space,
    \[
        d^H (\mathrm{Conv}\mathcal{M}_1, \mathrm{Conv}\mathcal{M}_2 )
        \leq
        d^H (\mathcal{M}_1, \mathcal{M}_2 )
        .
    \]
\end{proposition}

\begin{proposition}[{\citealp[Cf.][]{manfred2023lipshitz}}]
    \label{prop:Hausdorff_distance_and_Lipschitz_map}
    Let $L_f$ denote the Lipschitz constant of a map $f \colon \mathbb{R}^{d_1} \to \mathbb{R}^{d_2}$. For any compact sets $\mathcal{M}_1, \mathcal{M}_2$, we have
    \[
        d^H (f (\mathcal{M}_1), f (\mathcal{M}_2) )
        \leq 
        L_f
        d^H (\mathcal{M}_1, \mathcal{M}_2)
        .
    \]
\end{proposition}

\subsection{Lipschitz Property of Suboptimality Loss}

For a pseudometric $d_{\Phi}$ on the space $\Phi$, for any $\phi, \phi' \in \Phi$, define
\[
    d_{\Phi} (\phi, \phi') := \| d^H (\mathcal{X} (\phi, S), \mathcal{X} (\phi', S) ) \|_{\psi_2}
    .
\]
When each $\mathcal{X} (\phi , s)$ is nonempty, the equivalence relation of the pseudometric space $(\Phi , d_{\Phi})$ (i.e., $d_{\Phi} (\phi , \phi^\prime) = 0$) coincides with the equivalence relation $\sim$ on $\Phi$ defined in \Cref{sec:statistical_learning_theory}.

\begin{proposition}
    \label{prop:ReLU_Lipschitz}
    The Lipschitz constant of $\mathrm{ReLU}$ is $1$.
\end{proposition}

\begin{proof}
    For $u , v \in \mathbb{R}$, $\max (u , 0) \leq \max (v , 0) + \max (u - v , 0) \leq \max (v , 0) + | u - v |$; exchanging $u$ and $v$ gives $| \mathrm{ReLU} (u) - \mathrm{ReLU} (v) | \leq | u - v |$. Equality holds for $u = 1$, $v = 0$, so the Lipschitz constant is $1$.
\end{proof}

\begin{proposition}
    \label{prop:suboptimality_zero_Lipschitz}
    Assume that $\sup_{\theta \in \Theta} \| \theta \|_2 \leq 1$.
    Let $\hat{x}^{*}(s) = x^*(\theta^{\mathrm{true}}, \phi^{\mathrm{true}}, s)$.
    For any $s \in \mathcal{S}$ and any $x,x^\prime \in \mathcal{X}$, assume that
    \[
        \| f (x, s) - f (x^\prime, s) \|
        \leq L_f \| x - x^\prime \|
    \]
    holds.
    Moreover, assume that for any $x \in \mathcal{X}$ and $s \in \mathcal{S}$, the map $g(x,\bullet,s)\colon \Phi \to \mathbb{R}^J$ is a lattice homomorphism.
    Then, for any $s \in \mathcal{S}$, $\theta, \theta' \in \Theta$, and $\phi, \phi' \geq \phi^{\mathrm{true}}$ such that $\mathcal{X} (\phi , s)$ and $\mathcal{X} (\phi' , s)$ are nonempty, we have
    \begin{align*}
        & \left|
            \ell^{\mathrm{sub}, 0} (\theta, \phi, s)
            -
            \ell^{\mathrm{sub}, 0} (\theta', \phi', s)
        \right| 
        \notag\\
        & \leq
        L_f d^H ( \mathcal{X}(\phi, s), \mathcal{X}(\phi', s) )
        +
        L_f d^H ( \mathcal{X}(\phi^{\mathrm{true}}, s), \{ \hat{x}^* (s) \} )
        \| \theta - \theta' \|
        .
    \end{align*}
\end{proposition}

\begin{proof}
    By \Cref{prop:ReLU_Lipschitz},
    \begin{align}
        & \left|
            \ell^{\mathrm{sub},0} ( \theta , \phi , s )
            -
            \ell^{\mathrm{sub},0} ( \theta^\prime , \phi^\prime , s )  
        \right|
        \notag \\
        & \leq
        \left|
            \max_{x^{\star} \in \mathcal{X} (\phi , s)}
            \theta^\top \left( f ( x^{\star} , s) - f ( \hat{x}^{*}(s) , s) 
                 \right) 
            -
            \max_{x^{\star} \in \mathcal{X} (\phi^\prime , s)}
            \theta^{\prime \top} \left( f ( x^{\star} , s) - f ( \hat{x}^{*}(s) , s) 
                 \right)
        \right|
        \notag 
        \\
        & \leq 
        \left|
            \max_{x^{\star} \in \mathcal{X} (\phi , s)}
            \theta^\top \left( f ( x^{\star} , s) - f ( \hat{x}^{*}(s) , s) 
                 \right) 
            -
            \max_{x^{\star} \in \mathcal{X} (\phi^\prime , s)}
            \theta^{\top} \left( f ( x^{\star} , s) - f ( \hat{x}^{*}(s) , s) 
                 \right)
        \right|
        \notag \\
        & \quad 
        + 
        \left|
            \max_{x^{\star} \in \mathcal{X} (\phi^\prime , s)}
            \theta^\top \left( f ( x^{\star} , s) - f ( \hat{x}^{*}(s) , s) 
                 \right) 
            -
            \max_{x^{\star} \in \mathcal{X} (\phi^\prime , s)}
            \theta^{\prime \top} \left( f ( x^{\star} , s) - f ( \hat{x}^{*}(s) , s) 
                 \right)
        \right|
        \notag 
        \\
        & \leq 
        \left|
            \max_{x^{\star} \in \mathcal{X} (\phi , s)}
            \theta^\top \left( f ( x^{\star} , s)
                 \right) 
            -
            \max_{x^{\star} \in \mathcal{X} (\phi^\prime , s)}
            \theta^{\top} \left( f ( x^{\star} , s)
                 \right)
        \right|
        \notag \\
        & \quad 
        + 
        \left|
            \max_{x^{\star} \in \mathcal{X} (\phi^\prime , s)}
            \theta^\top \left( f ( x^{\star} , s) - f ( \hat{x}^{*}(s) , s) 
                 \right) 
            -
            \max_{x^{\star} \in \mathcal{X} (\phi^\prime , s)}
            \theta^{\prime \top} \left( f ( x^{\star} , s) - f ( \hat{x}^{*}(s) , s) 
                 \right)
        \right|
        .
        \label{eq:suboptimality_zero_Lipschitz_010}
    \end{align}
    By \Cref{prop:Hausdorff_distance_is_diff_LP,prop:convex_hull_is_Lipschitz,prop:Hausdorff_distance_and_Lipschitz_map}, the first term in \cref{eq:suboptimality_zero_Lipschitz_010} is
    \begin{align}
        & \left|
            \max_{x^{\star} \in \mathcal{X} (\phi , s)}
            \theta^\top \left( f ( x^{\star} , s)
                 \right) 
            -
            \max_{x^{\star} \in \mathcal{X} (\phi^\prime , s)}
            \theta^{\top} \left( f ( x^{\star} , s)
                 \right)
        \right|
        \notag \\
        & \leq  
        \left|
            \max_{a^{\star} \in \mathrm{Conv} f( \mathcal{X} (\phi , s) ,s )}
            \theta^\top a^\star 
            -
            \max_{a^{\star} \in \mathrm{Conv} f( \mathcal{X} (\phi^\prime , s) ,s )}
            \theta^{\top} a^\star
        \right|
        \notag \\
        & \leq 
        d^H \left(
            \mathrm{Conv} f ( \mathcal{X} (\phi , s) , s)
            ,
            \mathrm{Conv} f ( \mathcal{X} (\phi^\prime , s) , s)
        \right)
        \notag \\
        & \leq 
        d^H \left(
            f ( \mathcal{X} (\phi , s) , s)
            ,
            f ( \mathcal{X} (\phi^\prime , s) , s)
        \right) \leq 
        L_f
        d^H \left(
            \mathcal{X} (\phi , s)
            ,
            \mathcal{X} (\phi^\prime , s)
        \right)
        .
        \label{eq:suboptimality_zero_Lipschitz_020}
    \end{align}
    On the other hand, the second term of \cref{eq:suboptimality_zero_Lipschitz_010} is 
    \begin{align}
        &
        \left|
            \max_{x^{\star} \in \mathcal{X} (\phi^\prime , s)}
            \theta^\top \left( f ( x^{\star} , s) - f ( \hat{x}^{*}(s) , s) 
                 \right) 
            -
            \max_{x^{\star} \in \mathcal{X} (\phi^\prime , s)}
            \theta^{\prime \top} \left( f ( x^{\star} , s) - f ( \hat{x}^{*}(s) , s) 
                 \right)
        \right|
        \notag 
        \\
        & \leq
        \sup_{x^{\star} \in \mathcal{X} (\phi^\prime , s)}
        \left\|
            f ( x^{\star} , s) - f ( \hat{x}^{*}(s) , s)
        \right\|_2
        \left\|
            \theta - \theta^\prime
        \right\|_2   
        \notag 
        \\    
        & \leq
        d^H \left(
            f ( \mathcal{X} (\phi^\prime , s) , s) , \{ f ( \hat{x}^{*}(s) , s) \}
        \right)
        \left\|
            \theta - \theta^\prime
        \right\|_2   
        \notag  
        \\    
        & \leq
        L_f
        d^H \left(
            \mathcal{X} (\phi^\prime , s) , \{ \hat{x}^{*}(s) \}
        \right)
        \left\|
            \theta - \theta^\prime
        \right\|_2   
        .
        \notag  
    \end{align}
    By \Cref{prop:land-morphism_monotone} and $\phi^\prime \geq \phi^{\mathrm{true}}$, we have
    \begin{align}
        L_f
        d^H \left(
            \mathcal{X} (\phi^\prime , s) , \{ \hat{x}^{*}(s) \}
        \right)
        \left\|
            \theta - \theta^\prime
        \right\|_2   
        \leq 
        L_f
        d^H \left(
            \mathcal{X} (\phi^{\mathrm{true}} , s) , \{ \hat{x}^{*}(s) \}
        \right)
        \left\|
            \theta - \theta^\prime
        \right\|_2  
        .
        \label{eq:suboptimality_zero_Lipschitz_030}
    \end{align}
    
    The proposition follows from \cref{eq:suboptimality_zero_Lipschitz_010,eq:suboptimality_zero_Lipschitz_020,eq:suboptimality_zero_Lipschitz_030}.
\end{proof}

\begin{proposition}
    \label{prop:suboptimality_Lipschitz}
    Assume that $\sup_{\theta \in \Theta} \| \theta \|_2 \leq 1$.
    Let $\hat{x}^{*}(s) = x^*(\theta^{\mathrm{true}}, \phi^{\mathrm{true}}, s)$.
    For any $s \in \mathcal{S}$ and any $x,x^\prime \in \mathcal{X}$, assume that
    \[
        \| f (x, s) - f (x^\prime, s) \|
        \leq L_f \| x - x^\prime \|
    \]
    holds.
    Moreover, assume that for any $x \in \mathcal{X}$ and $s \in \mathcal{S}$, the map $g(x,\bullet,s)\colon \Phi \to \mathbb{R}^J$ is a lattice homomorphism.
    Then, for any $s \in \mathcal{S}$, $\theta, \theta' \in \Theta$, and any $\phi, \phi' \geq \phi^{\mathrm{true}}$ satisfying $g (\hat{x}^{*} (s) , \phi , s) \leq 0$ and $g (\hat{x}^{*} (s) , \phi^\prime , s) \leq 0$, we have
    \begin{align*}
        & \left|
            \ell^{\mathrm{sub},\lambda} ( \theta, \phi, s )
            -
            \ell^{\mathrm{sub},\lambda} ( \theta', \phi', s )
        \right|
        \\
        & \leq
        L_f d^H (\mathcal{X} ( \phi, s ), \mathcal{X} ( \phi', s ) )
        +
        L_f d^H ( \mathcal{X} ( \phi^{\mathrm{true}}, s ), \{ \hat{x}^{*}(s) \} )
        \| \theta - \theta' \|
        .
    \end{align*}
\end{proposition}

\begin{proof}
    By assumption, for $j = 1 , \ldots , J$ we have $g_j ( \hat{x}^{*} (s) , \phi , s) \leq 0$%
    , and hence
    \begin{equation}
        \mathrm{ReLU} ( g_j ( \hat{x}^{*}(s), \phi, s ) )
        = 0
        .
        \notag
    \end{equation}
    The same holds for $\phi^\prime$.
    For such $\phi$, the suboptimality loss satisfies
    \begin{align}
        \ell^{\mathrm{sub},\lambda} ( \theta, \phi, s )
        &=
        \ell^{\mathrm{sub},0} ( \theta, \phi, s )
        + \lambda \sum_{j=1}^J \mathrm{ReLU} ( g_j( \hat{x}^{*}(s), \phi, s ) ) = \ell^{\mathrm{sub},0} ( \theta, \phi, s )
        .
        \notag
    \end{align}
    The proposition then follows from \Cref{prop:suboptimality_zero_Lipschitz}.
\end{proof}

\begin{remark}
    \label{rem:penalty_vanish}
    When $\phi^{\mathrm{true}} \leq \phi \leq \phi^{\sup} (\{ s \})$, the condition $g (\hat{x}^{*} (s) , \phi , s) \leq 0$ is satisfied by the definition of $\phi^{\sup} (\{ s \})$.
    Moreover, when $\phi \sim \phi^{\mathrm{true}}$, for $\mathbb{P}_{\mathcal{S}}$-almost every $s$ we have $\mathcal{X} (\phi , s) = \mathcal{X} (\phi^{\mathrm{true}} , s) \ni \hat{x}^{*} (s)$, so the condition is satisfied for $\mathbb{P}_{\mathcal{S}}$-almost every $s$.
    In particular, if $\Phi (\delta) / \sim = \{ [ \phi^{\mathrm{true}} ] \}$ and $\Phi (\delta)$ is at most countable, then for any $\phi \in \Phi (\delta)$ and $\mathbb{P}_{\mathcal{S}}$-almost every $s$, we have $g (\hat{x}^{*} (s) , \phi , s) \leq 0$.
\end{remark}%

\begin{proposition}
    \label{prop:suboptimalityloss_subgaussian}
    For any $s \in \mathcal{S}$, $\theta, \theta' \in \Theta$, and $\phi, \phi' \in \Phi$, define
    \[
        d_s \big( (\theta, \phi), (\theta', \phi') \big)
        =
        L_f d^H \left( \mathcal{X} ( \phi^{\mathrm{true}}, s ), \{ \hat{x}^{*} (s) \} \right)
        \| \theta - \theta' \|
        +
        L_f d^H ( \mathcal{X} ( \phi, s ), \mathcal{X} ( \phi', s ) )
        .
    \]
    Then,
    \begin{align*}
        & d_{\mathcal{S}} \big( (\theta, \phi), (\theta', \phi') \big)
        \\
        & \leq 
        L_f \left\| d^H \left( \mathcal{X} ( \phi^{\mathrm{true}}, S ), \{ \hat{x}^{*} (S) \} \right) \right\|_{\psi_2}
        \| \theta - \theta' \|
        +
        L_f \left\| d^H ( \mathcal{X} ( \phi, S ), \mathcal{X} ( \phi', S ) ) \right\|_{\psi_2}
        .
    \end{align*}
\end{proposition}

\begin{proof}
    Since the sub-Gaussian norm satisfies the triangle inequality, the proposition follows.
\end{proof}

\begin{proposition}
    \label{prop:convering_integral_suboplimality}
    Let $d_s$ be the pseudometric in \Cref{prop:suboptimalityloss_subgaussian}.
    Then, (1)
        \begin{align*}
            & \int_0^{\infty}
            \sqrt{ \log N \left( \Theta \times \Phi, d_{\mathcal{S}}, \varepsilon \right) } d \varepsilon
            \notag \\
            & \leq 
            2L_f
            \left\| d^H \left( \mathcal{X} ( \phi^{\mathrm{true}}, S ), \{ \hat{x}^{*}(S) \} \right) \right\|_{\psi_2}
            \int_0^{\infty}
                \sqrt{ \log N \left( \Theta, \| \bullet \|_2, \varepsilon \right) } d \varepsilon
            \\
            & \qquad
            + 2L_f
            \int_0^{\infty}
                \sqrt{ \log N \left( \Phi, d_{\Phi}, \varepsilon \right) } d \varepsilon
            .
        \end{align*}
    (2)
        If $\Phi / \sim = \{ [ \phi^{\mathrm{true}} ] \}$, then
        \begin{align*}
            & \int_0^{\infty}
            \sqrt{ \log N \left( \Theta \times \Phi, d_{\mathcal{S}}, \varepsilon \right) } d \varepsilon
            \notag\\
            & =
            L_f
            \left\| d^H \left( \mathcal{X} ( \phi^{\mathrm{true}}, S ), \{ \hat{x}^{*}(S) \} \right) \right\|_{\psi_2}
            \int_0^{\infty}
                \sqrt{ \log N \left( \Theta, \| \bullet \|_2, \varepsilon \right) } d \varepsilon
            .
        \end{align*}
\end{proposition}

\begin{proof}
    (1) Statement (1) follows from \Cref{prop:suboptimalityloss_subgaussian,prop:covering_integral_on_product_sp_2,prop:covering_number_scaling}.

    (2) By the assumption, 
    \[
        N \left( \Theta \times \Phi, d_{\mathcal{S}}, \varepsilon \right)
        = N \left ( \Theta \times \{ \phi^{\mathrm{true}} \}, d_{\mathcal{S}}, \varepsilon \right).
    \]
    Therefore, statement (2) follows from \Cref{prop:covering_number_scaling}.
\end{proof}

\begin{proposition}
    \label{prop:separability_check}
    Suppose \Cref{assu:WIRL} and $\sup_{\theta \in \Theta} \| \theta \|_2 \leq 1$.
    Let $\Phi^\prime \subset \Phi$ be an at most countable set, and assume that for any $\phi \in \Phi^\prime$ and $s \in \mathcal{S}$, we have
    $M (\phi , s) := \sup_{x \in \mathcal{X} (\phi , s)} \left\| f (x , s) - f (\hat{x}^{*} (s) , s) \right\| < \infty$, and that $\mathbb{E} M (\phi , S) < \infty$ for any $\phi \in \Phi^\prime$.
    Let $d_s$ be the pseudometric of \Cref{prop:suboptimalityloss_subgaussian}.
    Then the process
    \[
        X_{(\theta , \phi)} := \frac{1}{N} \sum_{n=1}^N \ell^{\mathrm{sub},\lambda} ( \theta , \phi , S^{(n)}) - \mathbb{E} \ell^{\mathrm{sub},\lambda} ( \theta , \phi , S)
    \]
    and the process $(- X_{(\theta , \phi)})_{(\theta , \phi)}$ are separable on the pseudometric space $(\Theta \times \Phi^\prime , d_{\mathcal{S}})$ (\Cref{def:separable_process}).
\end{proposition}

\begin{proof}
    For each $\phi \in \Phi^\prime$ and $s \in \mathcal{S}$, in the same way as the bound on the second term in the proof of \Cref{prop:suboptimality_zero_Lipschitz}, for any $\theta , \theta^\prime \in \Theta$ we have
    \[
        \left| \ell^{\mathrm{sub},\lambda} (\theta , \phi , s) - \ell^{\mathrm{sub},\lambda} (\theta^\prime , \phi , s) \right|
        \leq M (\phi , s) \, \| \theta - \theta^\prime \|_2
    \]
    (the penalty term does not depend on $\theta$).
    Taking the expectation, we also have $| \mathbb{E} \ell^{\mathrm{sub},\lambda} (\theta , \phi , S) - \mathbb{E} \ell^{\mathrm{sub},\lambda} (\theta^\prime , \phi , S) | \leq \mathbb{E} M (\phi , S) \, \| \theta - \theta^\prime \|_2$.
    Therefore, for any $\omega$ and $\phi \in \Phi^\prime$, the map $\theta \mapsto X_{(\theta , \phi)} (\omega)$ is continuous.
    Since $\Theta$ is separable as a subset of a Euclidean space, there exists an at most countable $\Theta_{0} \subset \Theta$ that is dense with respect to the Euclidean norm.
    The set $\Theta^\prime := \Theta_{0} \times \Phi^\prime$ is at most countable, and by the continuity in $\theta$, for any $(\theta^{\prime\prime} , \phi^{\prime\prime}) \in \Theta \times \Phi^\prime$ and any $\omega$ we have
    \[
        \sup_{(\theta , \phi) \in \Theta \times \Phi^\prime} \left( X_{(\theta , \phi)} - X_{(\theta^{\prime\prime} , \phi^{\prime\prime})} \right)
        = \sup_{(\theta , \phi) \in \Theta^\prime} \left( X_{(\theta , \phi)} - X_{(\theta^{\prime\prime} , \phi^{\prime\prime})} \right)
        .
    \]
    Moreover, since $d_{\mathcal{S}} ((\theta , \phi) , (\theta^\prime , \phi)) \leq L_f \| d^H ( \mathcal{X} (\phi^{\mathrm{true}} , S), \{ \hat{x}^{*} (S) \} ) \|_{\psi_2} \, \| \theta - \theta^\prime \|$, the set $\Theta^\prime$ is dense in $(\Theta \times \Phi^\prime , d_{\mathcal{S}})$.
    Since $\theta \mapsto - X_{(\theta , \phi)} (\omega)$ is also continuous for any $\omega$, the above equality also holds with $X_{(\theta , \phi)}$ replaced by $- X_{(\theta , \phi)}$. Hence both processes are separable.
\end{proof}%

\begin{proposition}
    \label{prop:learning_theory_prob_suboptimality}

    We assume that $\sup_{\theta \in \Theta} \| \theta \|_2 \leq 1$. Let $\hat{x}^{*} \colon \mathcal{S} \to \mathcal{X}$ be an optimal solution map. For any $s \in \mathcal{S}$ and any $x, x' \in \mathcal{X}$, assume that
    \[
        \| f(x, s) - f(x', s) \|
        \leq L_f \| x - x' \|
        .
    \]
    Furthermore, for any $x \in \mathcal{X}$ and $s \in \mathcal{S}$, let $g(x, \bullet, s) \colon \Phi \to \mathbb{R}^J$ be a lattice homomorphism.
    Assume that $\Phi$ is at most countable, that $\phi \geq \phi^{\mathrm{true}}$ for any $\phi \in \Phi$, that for any $\phi \in \Phi$ and $s \in \mathcal{S}$, $\sup_{x \in \mathcal{X} (\phi , s)} \left\| f (x , s) - f (\hat{x}^{*} (s) , s) \right\| < \infty$, and that for $\mathbb{P}_{\mathcal{S}}$-almost every $s \in \mathcal{S}$, $g (\hat{x}^{*} (s) , \phi , s) \leq 0$.
    Let $(\theta^{*(N)} , \phi^{*(N)})$ be a minimizer of the empirical loss $\frac{1}{N} \sum_{n=1}^N \ell^{\mathrm{sub},\lambda} ( \theta , \phi , S^{(n)} )$ over $\Theta \times \Phi$.%
    
    Let $C > 6 \sqrt{3}$.

    Then,
    (1)
        For any $u \geq 0$,
        \begin{align*}
            \mathbb{P}
            \left(
                 \mathbb{E}_S \ell^{\mathrm{sub},\lambda} ( \theta^{*(N)}, \phi^{*(N)}, S )
                >
                \varepsilon(u, C, N, \Phi )
            \right)
            & \leq 2 \zeta \left( \frac{C^2}{36} - 2 \right) \exp \left( - \frac{C^2}{36} u^2 \right)
            .
            \notag
        \end{align*}
        where
        \begin{align}
            \varepsilon(u, C, N, \Phi) & := \frac{ 2 \sqrt{6 (1 + \log 2)}\, L_f C }{ \sqrt{N} } \notag \\
            & \qquad \times \biggl( 2 \| d^H ( \mathcal{X} (\phi^{\mathrm{true}} , S), \{ \hat{x}^{*} (S) \} ) \|_{\psi_2} \int_0^{\infty} \sqrt{\log N (\Theta , \| \bullet \|_2 ,\varepsilon )} d \varepsilon \notag \\
            & \qquad \quad + 2 \int_0^{\infty} \sqrt{\log N (\Phi , d_{\Phi} ,\varepsilon )} d \varepsilon \notag \\
            & \qquad \quad + u \bigl( \| d^H ( \mathcal{X} (\phi^{\mathrm{true}} , S), \{ \hat{x}^{*} (S) \} ) \|_{\psi_2} \mathrm{diam} (\Theta ) + \mathrm{diam} (\Phi , d_{\Phi}) \bigr) \biggr) .
        \end{align}

        (2)
        If $\Phi / \sim = \{ [ \phi^{\mathrm{true}} ] \}$, then for any $u \geq 0$,
        \begin{align}
            \mathbb{P}
            \left(
                 \mathbb{E}_S \ell^{\mathrm{sub},\lambda} \left( \theta^{*(N)}, \phi^{*(N)}, S \right)
                >
                \varepsilon(u, C, N)
            \right)
            & \leq 2 \zeta \left( \frac{C^2}{36} - 2 \right) \exp \left( - \frac{C^2}{36} u^2 \right)
            .
            \notag
        \end{align}
        where
        \begin{align*}
            \varepsilon(u, C, N) & := \frac{ 2 \sqrt{6 (1 + \log 2)}\, L_f C }{ \sqrt{N} } \\
            & \qquad \times \biggl( \left\| d^H \left( \mathcal{X} ( \phi^{\mathrm{true}}, S ), \{ \hat{x}^{*}(S) \} \right) \right\|_{\psi_2} \int_0^{\infty} \sqrt{ \log N \left( \Theta, \| \bullet \|_2, \varepsilon \right) } d \varepsilon \\
            & \qquad \quad + u \left\| d^H \left( \mathcal{X} ( \phi^{\mathrm{true}}, S ), \{ \hat{x}^{*}(S) \} \right) \right\|_{\psi_2} \mathrm{diam} (\Theta) \biggr) .
        \end{align*}
\end{proposition}

\begin{proof}
    By \Cref{prop:suboptimality_imitation},
    \[
        \mathbb{E}_S
        \ell^{\mathrm{sub},\lambda} ( \theta^{\mathrm{true}} , \phi^{\mathrm{true}} , S ) = 0
    \]
    holds.
    By assumption, for any $\phi \in \Phi$ and $\mathbb{P}_{\mathcal{S}}$-almost every $s$ we have $g (\hat{x}^{*} (s) , \phi , s) \leq 0$; hence, by \Cref{prop:suboptimality_Lipschitz}, for any $(\theta , \phi) , (\theta^\prime , \phi^\prime) \in \Theta \times \Phi$, \Cref{eq:assu:loss_estimated_pseudometric} holds for $\mathbb{P}_{\mathcal{S}}$-almost every $s$ (with the pseudometric $d_s$ of \Cref{prop:suboptimalityloss_subgaussian}). If any of the sub-Gaussian norms on the right-hand side of the conclusion is infinite, the claim is trivial; hence we may assume that $d_{\mathcal{S}}$ is finite. Moreover, the separability of the sign-reversed process of the corresponding empirical process on $(\Theta \times \Phi , d_{\mathcal{S}})$ follows from \Cref{prop:separability_check} with $\Phi^\prime = \Phi$; here $\mathbb{E} M (\phi , S) < \infty$ holds for any $\phi \in \Phi$, since $\phi \geq \phi^{\mathrm{true}}$ gives $M (\phi , s) \leq L_f d^H ( \mathcal{X} (\phi^{\mathrm{true}} , s), \{ \hat{x}^{*} (s) \} )$ as in the proof of \Cref{prop:suboptimality_zero_Lipschitz}, and $\mathbb{E} d^H ( \mathcal{X} (\phi^{\mathrm{true}} , S), \{ \hat{x}^{*} (S) \} ) \leq \sqrt{\log 2} \, \| d^H ( \mathcal{X} (\phi^{\mathrm{true}} , S), \{ \hat{x}^{*} (S) \} ) \|_{\psi_2} < \infty$. Then, since $0 \leq \ell^{\mathrm{sub},\lambda}$, $\mathbb{E}_S \ell^{\mathrm{sub},\lambda} ( \theta^{\mathrm{true}} , \phi^{\mathrm{true}} , S ) = 0$, and \cref{eq:assu:loss_estimated_pseudometric}, we have $\mathbb{E}_S | \ell^{\mathrm{sub},\lambda} (\theta , \phi , S) | < \infty$ for any $(\theta , \phi) \in \Theta \times \Phi$, and $(\theta^{\mathrm{true}} , \phi^{\mathrm{true}})$ minimizes the expected loss.%
    
    By \Cref{theo:learning_theory_prob} and \Cref{prop:suboptimality_Lipschitz,prop:suboptimalityloss_subgaussian,prop:convering_integral_suboplimality,prop:diam_of_scaling,prop:diam_on_product_sp},
    (1) and (2) follow.
\end{proof}

\subsection{Statistical Learning Theory of Inverse Optimization}

\begin{theorem}
    \label{theo:suboptimality_learning_prob}
    Let $\delta \in \mathbb{R}_{\geq 0}^J$.
    Assume that $\sup_{\theta \in \Theta} \| \theta \|_2 \leq 1$.
    Let $\hat{x}^{*} (s) = x^*(\theta^{\mathrm{true}}, \phi^{\mathrm{true}}, s)$.
    Assume that for any $s \in \mathcal{S}$ and any $x , x^\prime \in \mathcal{X}$,
    \[
        \| f (x, s) - f (x^\prime, s) \|
        \leq L_f \| x - x^\prime \|
    \]
    holds.
    Assume further that for any $x \in \mathcal{X}$ and $s\in \mathcal{S}$, the mapping $g(x , \bullet , s ) \colon \Phi \to \mathbb{R}^J$ is a lattice homomorphism.
    Assume that $\Phi (\delta)$ is at most countable, and that for any $\phi \in \Phi (\delta)$ and $s \in \mathcal{S}$, $\sup_{x \in \mathcal{X} (\phi , s)} \left\| f (x , s) - f (\hat{x}^{*} (s) , s) \right\| < \infty$.
    Assume that for any $\phi \in \Phi (\delta)$ and $\mathbb{P}_{\mathcal{S}}$-almost every $s \in \mathcal{S}$, we have $g (\hat{x}^{*} (s) , \phi , s) \leq 0$.
    For a sample $(S^{(1)} , \ldots , S^{(N)} )$, let $\phi^{*(N)} := \phi^{\sup} (\mathcal{S}_N)$, and let $\theta^{*(N)} \in \Theta$ satisfy $\frac{1}{N} \sum_{n=1}^N \ell^{\mathrm{sub},\lambda} ( \theta^{*(N)} , \phi^{*(N)} , S^{(n)} ) = 0$ (it may depend on the sample).%
    
    Let $C > 6 \sqrt{3}$.

    Then, for any $u \geq 0$, (1)
    \begin{align*}
        & \mathbb{P} \left(
            \phi^{*(N)} \leq \phi^{\mathrm{true}} + \delta
            \text{ and }
            \mathbb{E} \ell^{\mathrm{sub},\lambda} ( \theta^{*(N)} , \phi^{*(N)} , S )
            \leq \varepsilon (u , C, N , \Phi (\delta) )
        \right)
        \\
        & \geq 1 - \sum_{j=1}^J \mathbb{P} \left(
            \phi_j^{\sup} (\{ S \} ) >  \phi_j^{\mathrm{true}} + \delta_j
        \right)^N
        - 2 \zeta \left( \frac{C^2}{36} -2 \right) \exp \left( -\frac{C^2}{36} u^2  \right)
        .
    \end{align*}
    (2) If $\Phi (\delta) / \sim  = \{ [ \phi^{\mathrm{true} } ] \}$, then
        \begin{align*}
        & \mathbb{P} \left(
            \phi^{*(N)} \sim \phi^{\mathrm{true}}
            \text{ and }
            \mathbb{E} \ell^{\mathrm{sub},\lambda} ( \theta^{*(N)} , \phi^{*(N)} , S )
            \leq \varepsilon (u , C, N )
        \right)
        \\
        & \geq 1 - \sum_{j=1}^J \mathbb{P} \left(
            \phi_j^{\sup} (\{ S \} ) >  \phi_j^{\mathrm{true}} + \delta_j
        \right)^N
         - 2 \zeta \left( \frac{C^2}{36} -2 \right) \exp \left( -\frac{C^2}{36} u^2  \right)
    \end{align*}
\end{theorem}

\begin{proof}
    First,
    \begin{align}
        & 1 - \mathbb{P} \left(
                \phi^{*(N)} \leq \phi^{\mathrm{true}} + \delta
                \text{ and }
                \mathbb{E} \ell^{\mathrm{sub},\lambda} ( \theta^{*(N)} , \phi^{*(N)} , S )
                \leq \varepsilon (u , C, N , \Phi (\delta) )
        \right)
        \notag \\
        & =
        \mathbb{P} \left(
                \exists j , \phi_j^{*(N)} > \phi_j^{\mathrm{true}} + \delta_j
                \text{ or }
                \mathbb{E} \ell^{\mathrm{sub},\lambda} ( \theta^{*(N)} , \phi^{*(N)} , S )
                > \varepsilon (u , C, N , \Phi (\delta) )
        \right)
        \notag \\
        & \leq
        \mathbb{P} \left(
            \exists j , \phi_j^{*(N)} > \phi_j^{\mathrm{true}} + \delta_j
        \right)
        \label{eq:subopt_stat_learning_theo_010-01}
        \\
        & \quad +
        \mathbb{P} \left(
                \phi^{*(N)} \leq \phi^{\mathrm{true}} + \delta
                \text{ and }
                \mathbb{E} \ell^{\mathrm{sub},\lambda} ( \theta^{*(N)} , \phi^{*(N)} , S )
                > \varepsilon (u , C, N , \Phi (\delta) )
        \right)
        .
        \label{eq:subopt_stat_learning_theo_010-02}
    \end{align}
    For \Cref{eq:subopt_stat_learning_theo_010-01}, we have
    \begin{align}
        & \mathbb{P} \left(
            \exists j , \phi_j^{*(N)} > \phi_j^{\mathrm{true}} + \delta_j
        \right)
        \notag\\
        & =
        \mathbb{P} \left(
            \bigcup_{j=1}^J \left\{  \phi_j^{*(N)} > \phi_j^{\mathrm{true}} + \delta_j \right\}
        \right)
        \leq
        \sum_{j=1}^J
        \mathbb{P} \left(
            \phi_j^{*(N)} > \phi_j^{\mathrm{true}} + \delta_j
        \right)
        \notag
        \\
        & \leq
        \sum_{j=1}^J
        \mathbb{P} \left(
            \forall n = 1, \ldots , N , \,
            \phi_j^{\sup} ( \{ S^{(n)} \} ) >  \phi_j^{\mathrm{true}} + \delta_j
        \right)
        .
        \notag
    \end{align}
    Here we used the fact that, by \Cref{prop:constraint_learning_lattice_map}, 
    \[
        \phi^{*(N)} = \bigwedge_{n=1}^{N} \phi^{\sup} (\{ S^{(n)} \}) \leq \phi^{\sup} (\{ S^{(n)} \}).
    \]
    Since the random variables $S^{(n)}$ are independent for $n = 1 , \ldots , N$, it follows that
    \begin{align}
        & \sum_{j=1}^J
        \mathbb{P} \left(
            \forall n = 1, \ldots , N , \,
            \phi_j^{\sup} ( \{ S^{(n)} \} ) >  \phi_j^{\mathrm{true}} + \delta_j
        \right)
        \notag \\
        & \leq
        \sum_{j=1}^J
        \prod_{n=1}^N
        \mathbb{P} \left(
            \phi_j^{\sup} ( \{ S^{(n)} \} ) >  \phi_j^{\mathrm{true}} + \delta_j
        \right)
        \notag
        \\
        & \leq
        \sum_{j=1}^J
        \mathbb{P} \left(
            \phi_j^{\sup} ( \{ S \} ) >  \phi_j^{\mathrm{true}} + \delta_j
        \right)^N
         .
        \label{eq:subopt_stat_learning_theo_020-01}
    \end{align}

    On the other hand, we bound \Cref{eq:subopt_stat_learning_theo_010-02}.
    Define the pair $(\check{\theta}^{(N)} , \check{\phi}^{(N)})$ to be $(\theta^{*(N)} , \phi^{*(N)})$ on the event $\{ \phi^{*(N)} \leq \phi^{\mathrm{true}} + \delta \}$, and $(\theta^{\mathrm{true}} , \phi^{\mathrm{true}})$ outside it.
    By \Cref{prop:true_constraints_para_dominated}, we have $\phi^{*(N)} \geq \phi^{\mathrm{true}}$, so on this event $\phi^{*(N)} \in \Phi (\delta)$, and by assumption the empirical loss of $(\theta^{*(N)} , \phi^{*(N)})$ is $0$.
    Moreover, by \Cref{prop:suboptimality_imitation}, $\mathbb{E}_S \ell^{\mathrm{sub},\lambda} ( \theta^{\mathrm{true}} , \phi^{\mathrm{true}} , S ) = 0$ and $\ell^{\mathrm{sub},\lambda} \geq 0$, so almost surely $\frac{1}{N} \sum_{n=1}^N \ell^{\mathrm{sub},\lambda} ( \theta^{\mathrm{true}} , \phi^{\mathrm{true}} , S^{(n)} ) = 0$.
    Therefore, $(\check{\theta}^{(N)} , \check{\phi}^{(N)})$ is almost surely a minimizer of the empirical loss over $\Theta \times \Phi (\delta)$, and on the event we have 
    \[
        \mathbb{E}_S \ell^{\mathrm{sub},\lambda} ( \theta^{*(N)} , \phi^{*(N)} , S ) = \mathbb{E}_S \ell^{\mathrm{sub},\lambda} ( \check{\theta}^{(N)} , \check{\phi}^{(N)} , S ).
    \]
    Hence,
    \begin{align}
        &
        \mathbb{P} \left(
                \phi^{*(N)} \leq \phi^{\mathrm{true}} + \delta
                \text{ and }
                \mathbb{E} \ell^{\mathrm{sub},\lambda} ( \theta^{*(N)} , \phi^{*(N)} , S )
                > \varepsilon (u , C, N , \Phi (\delta) )
        \right)
        \notag \\
        &
        \leq
        \mathbb{P} \left(
            \mathbb{E} \ell^{\mathrm{sub},\lambda} ( \check{\theta}^{(N)} , \check{\phi}^{(N)} , S )
            > \varepsilon (u , C, N , \Phi (\delta) )
        \right)
        \notag
    \end{align}
    holds.
    Applying \Cref{prop:learning_theory_prob_suboptimality} with $\Phi(\delta)$ in place of $\Phi$, we obtain
    \begin{align}
        & \mathbb{P}
        \left(
            \mathbb{E}_S \ell^{\mathrm{sub},\lambda} ( \check{\theta}^{(N)} , \check{\phi}^{(N)} , S )
            >
            \varepsilon (u , C, N , \Phi (\delta) )
        \right)
        \notag\\
        & \leq 2 \zeta \left( \frac{C^2}{36} -2 \right) \exp \left( -\frac{C^2}{36} u^2  \right)
        .
        \label{eq:subopt_stat_learning_theo_030-02}
    \end{align}
    Combining \Cref{eq:subopt_stat_learning_theo_010-01,eq:subopt_stat_learning_theo_010-02,eq:subopt_stat_learning_theo_020-01,eq:subopt_stat_learning_theo_030-02} yields (1).

    (2) can be proved in the same manner.
\end{proof}

\begin{theorem}
    \label{theo:suboptimality_learning_E}
    Let $\delta \in \mathbb{R}_{\geq 0}^J$.
    We assume that $\sup_{\theta \in \Theta} \| \theta \|_2 \leq 1$. Let $\hat{x}^{*} \colon \mathcal{S} \to \mathcal{X}$ be the optimal solution map. For any $s \in \mathcal{S}$ and any $x, x' \in \mathcal{X}$, assume that
    \[
        \| f(x, s) - f(x', s) \|
        \leq L_f \| x - x' \|
        .
    \]
    Furthermore, for any $x \in \mathcal{X}$ and $s \in \mathcal{S}$, let $g(x, \bullet, s) \colon \Phi \to \mathbb{R}^J$ be a lattice homomorphism.
    Assume that $\Phi (\delta)$ is at most countable, and that for any $\phi \in \Phi (\delta)$ and $s \in \mathcal{S}$, $\sup_{x \in \mathcal{X} (\phi , s)} \left\| f (x , s) - f (\hat{x}^{*} (s) , s) \right\| < \infty$.
    Assume that for any $\phi \in \Phi (\delta)$ and $\mathbb{P}_{\mathcal{S}}$-almost every $s \in \mathcal{S}$, we have $g (\hat{x}^{*} (s) , \phi , s) \leq 0$.
    For a sample $(S^{(1)} , \ldots , S^{(N)} )$, let $\phi^{*(N)} := \phi^{\sup} (\mathcal{S}_N)$, and let $\theta^{*(N)} \in \Theta$ satisfy $\frac{1}{N} \sum_{n=1}^N \ell^{\mathrm{sub},\lambda} ( \theta^{*(N)} , \phi^{*(N)} , S^{(n)} ) = 0$ (it may depend on the sample).%
    
    Assume that $\sum_{j=1}^J \mathbb{P} \left( \phi_j^{\sup} (\{ S \}) > \phi_j^{\mathrm{true}} + \delta_j \right)^N < 1$.
    Then, (1)
        \begin{align*}
            & \mathbb{E}_{S^{(1)} , \ldots , S^{(N)}} \left[
                \mathbb{E}_{S} \ell^{\mathrm{sub},\lambda} ( \theta^{*(N)}, \phi^{*(N)}, S )
            \,\middle|\,
                \phi^{*(N)} \in \Phi(\delta)
            \right]
            \\
            & \leq
            \left( 1 - \sum_{j=1}^J \mathbb{P} \left( \phi_j^{\sup} (\{ S \} ) >  \phi_j^{\mathrm{true}} + \delta_j \right)^N \right)^{-1}
            \frac{16 \sqrt{2 (1 + \log 2)}\, L_f }{\sqrt{N}}
            \\
            & \quad \times
            \left(
            \begin{array}{l}
                \displaystyle
                \left\| d^H \left( \mathcal{X} ( \phi^{\mathrm{true}}, S ), \{ \hat{x}^{*} (S) \} \right) \right\|_{\psi_2}
                \int_0^{\infty}
                    \sqrt{ \log N ( \Theta, \| \bullet \|_2, \varepsilon ) } d \varepsilon
                \\
                \displaystyle
                +
                \int_0^{\infty}
                    \sqrt{ \log N ( \Phi(\delta), d_{\Phi}, \varepsilon ) } d \varepsilon
            \end{array}
            \right)
            .
        \end{align*}

    (2)
        If $\Phi(\delta) / \sim = \{ [ \phi^{\mathrm{true}} ] \}$, then
        \begin{align*}
            & \mathbb{E}_{S^{(1)} , \ldots , S^{(N)}} \left[
                \mathbb{E}_{S} \ell^{\mathrm{sub},\lambda} ( \theta^{*(N)}, \phi^{*(N)}, S )
            \,\middle|\,
                \phi^{*(N)} \in \Phi (\delta)
            \right]
            \\
            & \leq
            \left( 1 - \sum_{j=1}^J \mathbb{P} \left( \phi_j^{\sup} (\{ S \} ) >  \phi_j^{\mathrm{true}} + \delta_j \right)^N \right)^{-1}
            \frac{8 \sqrt{2 (1 + \log 2)}\, L_f}{\sqrt{N}}
            \\
            & \quad \times
            \left\| d^H \left( \mathcal{X} ( \phi^{\mathrm{true}}, S ), \{ \hat{x}^{*} (S) \} \right) \right\|_{\psi_2}
            \int_0^{\infty}
                \sqrt{ \log N ( \Theta, \| \bullet \|_2, \varepsilon ) } d \varepsilon
            .
        \end{align*}
\end{theorem}

\begin{proof}
    By \Cref{prop:suboptimality_imitation},
    \[
        \mathbb{E}_S
        \ell^{\mathrm{sub},\lambda} ( \theta^{\mathrm{true}} , \phi^{\mathrm{true}} , S ) = 0
    \]
holds.
    Since $\ell^{\mathrm{sub},\lambda} \geq 0$, almost surely $\frac{1}{N} \sum_{n=1}^N \ell^{\mathrm{sub},\lambda} ( \theta^{\mathrm{true}} , \phi^{\mathrm{true}} , S^{(n)} ) = 0$.
    Define the pair $(\check{\theta}^{(N)} , \check{\phi}^{(N)})$ to be $(\theta^{*(N)} , \phi^{*(N)})$ on the event $A := \{ \phi^{*(N)} \in \Phi (\delta) \}$, and $(\theta^{\mathrm{true}} , \phi^{\mathrm{true}})$ outside it.
    By assumption, on $A$ the empirical loss of $(\theta^{*(N)} , \phi^{*(N)})$ is $0$, so $(\check{\theta}^{(N)} , \check{\phi}^{(N)})$ is almost surely a minimizer of the empirical loss over $\Theta \times \Phi (\delta)$.
    Since $\ell^{\mathrm{sub},\lambda} \geq 0$, we have
    \[
        \mathbb{E} \left[
            \mathbb{E}_S \ell^{\mathrm{sub},\lambda} ( \theta^{*(N)} , \phi^{*(N)} , S )
            \, \mathbf{1}_{A}
        \right]
        \leq
        \mathbb{E} \, \mathbb{E}_S \ell^{\mathrm{sub},\lambda} ( \check{\theta}^{(N)} , \check{\phi}^{(N)} , S )
        .
    \]
    Moreover, by the definition of the conditional expectation,
    \begin{align*}
        & \mathbb{E}_{S^{(1)} , \ldots , S^{(N)}} \left[
            \mathbb{E}_S \ell^{\mathrm{sub},\lambda} ( \theta^{*(N)} , \phi^{*(N)} , S )
        \,\middle|\, A
        \right]
        \\
        & = \frac{1}{\mathbb{P} (A)}
        \mathbb{E} \left[
            \mathbb{E}_S \ell^{\mathrm{sub},\lambda} ( \theta^{*(N)} , \phi^{*(N)} , S )
            \, \mathbf{1}_{A}
        \right]
        ,
    \end{align*}
    and since $\phi^{*(N)} \geq \phi^{\mathrm{true}}$ (\Cref{prop:true_constraints_para_dominated}) implies $A^{c} = \{ \exists j , \, \phi_j^{*(N)} > \phi_j^{\mathrm{true}} + \delta_j \}$, in the same way as \Cref{eq:subopt_stat_learning_theo_020-01} (in the proof of \Cref{theo:suboptimality_learning_prob}), we have
    \[
        \mathbb{P} (A)
        \geq 1 - \sum_{j=1}^J \mathbb{P} \left( \phi_j^{\sup} (\{ S \}) > \phi_j^{\mathrm{true}} + \delta_j \right)^N
        > 0
    \]
    (the positivity is by assumption).
    Applying \Cref{theo:learning_theory_expectation} with $\Phi = \Phi( \delta )$ to the right-hand side (the assumptions of \Cref{theo:learning_theory_expectation} are verified as in the proof of \Cref{prop:learning_theory_prob_suboptimality}), statements (1) and (2) follow from \Cref{prop:suboptimality_Lipschitz,prop:suboptimalityloss_subgaussian,prop:convering_integral_suboplimality}, \Cref{prop:separability_check}, and $\mathbb{E}_S \ell^{\mathrm{sub},\lambda} ( \theta^{\mathrm{true}} , \phi^{\mathrm{true}} , S ) = 0$.
\end{proof}

\subsection{Piecewise-Linear Functions}

\begin{theorem}
    \label{theo:subpoptimality_learning_prob_pl}
    Let $\delta \in \mathbb{R}_{\geq 0}^J$.
    Moreover, assume that either (i) $\mathcal{S}$ is finite, or (ii) $\Theta = \Delta^{D-1}$, $\mathcal{X} \subset \mathbb{Z}^{D_{\mathcal{X}}}$, and $f (\mathcal{X} \times \mathcal{S}) \subset \mathbb{Z}^{D}$ (i.e., the forward problem is an ILP).
    Suppose \Cref{assu:WIRL,assu:WIRL-uniqueness,assu:IOP_std}.
    For a sample $(S^{(1)} , \ldots , S^{(N)} )$, let $\theta^{*(N)}\in \Theta$ and $\phi^{*(N)}\in \Phi$ denote, respectively, the weight and the constraint parameter output by \Cref{alg:Solve_IOP} with $\varepsilon = 0$, where line~3 is implemented by \Cref{alg:psgd_suboptimality_loss} with a sufficiently large sample-independent $K$ and with the learning rate chosen as either SRSS or SRSL.
    Assume that
    $\Phi (\delta) / \sim  = \{ [ \phi^{\mathrm{true} } ] \}$.
    Assume that $\Phi (\delta)$ is at most countable, and that for any $\phi \in \Phi (\delta)$ and $s \in \mathcal{S}$, $\sup_{x \in \mathcal{X} (\phi , s)} \left\| f (x , s) - f (\hat{x}^{*} (s) , s) \right\| < \infty$.
    Assume that $\sup_{\theta \in \Theta} \| \theta \|_2 \leq 1$ and that $\Theta$ is contained in some $(D-1)$-dimensional affine subspace (e.g., $\Theta = \Delta^{D-1}$).
    Let $C > 6 \sqrt{3}$.
    Then, for any $u \geq 0$,
        \begin{align*}
        & \mathbb{P} \biggl( \phi^{*(N)} \sim \phi^{\mathrm{true}} \\
        & \qquad \text{and } \mathbb{E} \ell^{\mathrm{sub},\lambda} ( \theta^{*(N)} , \phi^{*(N)} , S ) \\
        & \qquad \leq \frac{ 2 \sqrt{6 (1 + \log 2)}\, L_f C}{\sqrt{N}} \Bigl( 2 \sqrt{2 \log 2} \| d^H ( \mathcal{X} (\phi^{\mathrm{true}} , S), \{ \hat{x}^{*} (S) \} ) \|_{\psi_2} \sqrt{D-1} \\
        & \qquad \quad + u \| d^H ( \mathcal{X} (\phi^{\mathrm{true}} , S), \{ \hat{x}^{*} (S) \} ) \|_{\psi_2} \mathrm{diam} (\Theta ) \Bigr) \biggr) \\
        & \geq 1 - \sum_{j=1}^J \mathbb{P} \left( \phi_j^{\sup} ( \{ S \} ) >  \phi_j^{\mathrm{true}} + \delta_j \right)^N - 2 \zeta \left( \frac{C^2}{36} -2 \right) \exp \left( -\frac{C^2}{36} u^2  \right) .
    \end{align*}
\end{theorem}

\begin{proof}
    The output of \Cref{alg:Solve_IOP} satisfies $\phi^{*(N)} = \phi^{\sup} (\mathcal{S}_N)$, and applying \Cref{theo:alg:Solve_IOP_is_solve_IOP} to the uniform distribution on $\mathcal{S}_N$, the empirical loss of the output with sufficiently large $K$ is $0$.
    Moreover, by $\Phi (\delta) / \sim = \{ [ \phi^{\mathrm{true}} ] \}$ and \Cref{rem:penalty_vanish}, for any $\phi \in \Phi (\delta)$ and $\mathbb{P}_{\mathcal{S}}$-almost every $s$, we have $g (\hat{x}^{*} (s) , \phi , s) \leq 0$.
    Hence the assumptions of \Cref{theo:suboptimality_learning_prob} are satisfied.%
    
    Since $\Theta$ is contained in the intersection of a $(D-1)$-dimensional affine subspace and the closed unit ball $B^D$, which is a $(D-1)$-dimensional closed ball of radius at most $1$, $(\Theta , \| \bullet \|_2)$ is isometric to a subset of $B^{D-1}$. Hence, applying \Cref{prop:covering_integral_unit_ball} in dimension $D-1$, we obtain $\int_0^{\infty} \sqrt{\log N (\Theta , \| \bullet \|_2 ,\varepsilon )} d \varepsilon \leq 2 \sqrt{2 \log 2} \sqrt{D-1}$.
    The claim follows from \Cref{theo:suboptimality_learning_prob} and this bound.
\end{proof}

\begin{theorem}
    \label{theo:subpoptimality_learning_E_pl}
    Suppose the setting of \Cref{theo:subpoptimality_learning_prob_pl}.
    Assume that 
    \[
        \sum_{j=1}^J \mathbb{P} \left( \phi_j^{\sup} (\{ S \}) > \phi_j^{\mathrm{true}} + \delta_j \right)^N < 1 .
    \]
    Then
    \begin{align*}
        & \mathbb{E}_{S^{(1)} , \ldots , S^{(N)}} \left[
            \mathbb{E}_{S} \ell^{\mathrm{sub},\lambda} ( \theta^{*(N)} , \phi^{*(N)} , S )
        \, \middle| \, \phi^{*(N)} \in \Phi (\delta)
        \right]
        \notag \\
        & \leq
        \left( 1 - \sum_{j=1}^J \mathbb{P} \left( \phi_j^{\sup} (\{ S \} ) >  \phi_j^{\mathrm{true}} + \delta_j \right)^N \right)^{-1}
        \frac{32 \sqrt{(1 + \log 2) \log 2}\, L_f}{\sqrt{N}}
        \\
        & \quad \times
            \| d^H \left( \mathcal{X} (\phi^{\mathrm{true}} , S), \{ \hat{x}^{*} (S) \}\right) \|_{\psi_2}
            \sqrt{D-1}
            \notag
        .
    \end{align*}
\end{theorem}

\begin{proof}
    As in the proof of \Cref{theo:subpoptimality_learning_prob_pl}, the assumptions of \Cref{theo:suboptimality_learning_E} are satisfied.
    The claim follows from \Cref{theo:suboptimality_learning_E} and the bound
    \[
        \int_0^{\infty} \sqrt{\log N (\Theta , \| \bullet \|_2 ,\varepsilon )} d \varepsilon \leq 2 \sqrt{2 \log 2} \sqrt{D-1}
    \]
    obtained as in the proof of \Cref{theo:subpoptimality_learning_prob_pl}.
\end{proof}

\begin{remark}
\Cref{theo:subpoptimality_learning_E_pl_short,theo:subpoptimality_learning_prob_pl_short}
    apply to the case where $\Phi$ is a discrete space.
    When $\Phi$ contains continuous variables, $\Phi (\delta)$ is in general uncountable. In this case, by \Cref{prop:Phi_delta_maximum}, for statement (2) of \Cref{theo:suboptimality_learning_prob,theo:suboptimality_learning_E} and for \Cref{theo:subpoptimality_learning_prob_pl,theo:subpoptimality_learning_E_pl}, the countability of $\Phi (\delta)$ can be replaced by the existence of a maximum element of $\Phi (\delta)$. Since $\Phi (\delta)$ is bounded and closed under $\lor$, a maximum element exists whenever $\Phi (\delta)$ is closed: for each $j$, take $\phi^{[j]} \in \Phi (\delta)$ attaining $\max_{\phi \in \Phi (\delta)} \phi_j$; then $\bigvee_{j=1}^J \phi^{[j]} \in \Phi (\delta)$ is the maximum element.
\end{remark}

\begin{proposition}
    \label{prop:Phi_delta_maximum}
    In \Cref{theo:suboptimality_learning_prob,theo:suboptimality_learning_E}, replace the assumption that $\Phi (\delta)$ is at most countable and the assumption that $g (\hat{x}^{*} (s) , \phi , s) \leq 0$ for any $\phi \in \Phi (\delta)$ and $\mathbb{P}_{\mathcal{S}}$-almost every $s$ by the assumption that $\Phi (\delta) / \sim = \{ [ \phi^{\mathrm{true}} ] \}$ and $\Phi (\delta)$ has a maximum element $\bar{\phi}$.
    Then statement (2) of \Cref{theo:suboptimality_learning_prob} and statement (2) of \Cref{theo:suboptimality_learning_E} hold.
    Similarly, in \Cref{theo:subpoptimality_learning_prob_pl,theo:subpoptimality_learning_E_pl}, replace the assumption that $\Phi (\delta)$ is at most countable by the assumption that $\Phi (\delta)$ has a maximum element. Then the conclusions of \Cref{theo:subpoptimality_learning_prob_pl,theo:subpoptimality_learning_E_pl} hold.
\end{proposition}

\begin{proof}
    Let $\phi \in \Phi (\delta)$.
    Since $\phi^{\mathrm{true}} \leq \phi \leq \bar{\phi}$ and $g (x , \bullet , s)$ is order-preserving for any $x \in \mathcal{X}$ and $s \in \mathcal{S}$ (\Cref{prop:land-morphism_monotone}), we have $g (x , \phi^{\mathrm{true}} , s) \leq g (x , \phi , s) \leq g (x , \bar{\phi} , s)$, and hence
    \[
        \mathcal{X} (\bar{\phi} , s) \subset \mathcal{X} (\phi , s) \subset \mathcal{X} (\phi^{\mathrm{true}} , s)
        .
    \]
    Since $\bar{\phi} \in \Phi (\delta)$ and $\Phi (\delta) / \sim = \{ [ \phi^{\mathrm{true}} ] \}$, there exists a $\mathbb{P}_{\mathcal{S}}$-null set $\mathcal{S}_0 \subset \mathcal{S}$ such that $\mathcal{X} (\bar{\phi} , s) = \mathcal{X} (\phi^{\mathrm{true}} , s)$ for any $s \in \mathcal{S} \setminus \mathcal{S}_0$.
    Hence, for any $s \in \mathcal{S} \setminus \mathcal{S}_0$ and any $\phi \in \Phi (\delta)$, we have $\mathcal{X} (\phi , s) = \mathcal{X} (\phi^{\mathrm{true}} , s) \ni \hat{x}^{*} (s)$, so that $g (\hat{x}^{*} (s) , \phi , s) \leq 0$ and, by the definition of the suboptimality loss,
    \begin{equation}
        \ell^{\mathrm{sub},\lambda} ( \theta , \phi , s )
        = \ell^{\mathrm{sub},\lambda} ( \theta , \phi^{\mathrm{true}} , s )
        \quad \text{for any } \theta \in \Theta
        .
        \label{eq:Phi_delta_maximum_010}
    \end{equation}
    Taking the expectation, we obtain $\mathbb{E}_S \ell^{\mathrm{sub},\lambda} ( \theta , \phi , S ) = \mathbb{E}_S \ell^{\mathrm{sub},\lambda} ( \theta , \phi^{\mathrm{true}} , S )$ for any $(\theta , \phi) \in \Theta \times \Phi (\delta)$.

    Let $A := \{ \phi^{*(N)} \in \Phi (\delta) \}$ and $B := \{ S^{(n)} \notin \mathcal{S}_0 \text{ for } n = 1 , \ldots , N \}$; then $\mathbb{P} (B) = 1$.
    Define $\check{\theta}^{(N)}$ to be $\theta^{*(N)}$ on $A$ and $\theta^{\mathrm{true}}$ outside $A$.
    On $A \cap B$, by \Cref{eq:Phi_delta_maximum_010}, the empirical loss of $(\check{\theta}^{(N)} , \phi^{\mathrm{true}})$ equals that of $(\theta^{*(N)} , \phi^{*(N)})$, which is $0$ by assumption.
    Outside $A$, the empirical loss of $(\theta^{\mathrm{true}} , \phi^{\mathrm{true}})$ is almost surely $0$, as in the proof of \Cref{theo:suboptimality_learning_E}.
    Hence $(\check{\theta}^{(N)} , \phi^{\mathrm{true}})$ is almost surely a minimizer of the empirical loss over $\Theta \times \{ \phi^{\mathrm{true}} \}$.
    Moreover, on $A$ we have
    \begin{equation}
        \mathbb{E}_S \ell^{\mathrm{sub},\lambda} ( \theta^{*(N)} , \phi^{*(N)} , S )
        = \mathbb{E}_S \ell^{\mathrm{sub},\lambda} ( \check{\theta}^{(N)} , \phi^{\mathrm{true}} , S )
        .
        \label{eq:Phi_delta_maximum_020}
    \end{equation}
    The singleton $\{ \phi^{\mathrm{true}} \}$ satisfies the assumptions of \Cref{prop:learning_theory_prob_suboptimality} and \Cref{theo:learning_theory_expectation} with $\Phi = \{ \phi^{\mathrm{true}} \}$ as in the proofs of \Cref{theo:suboptimality_learning_prob,theo:suboptimality_learning_E}: it is at most countable, $\{ \phi^{\mathrm{true}} \} / \sim = \{ [ \phi^{\mathrm{true}} ] \}$, and $g (\hat{x}^{*} (s) , \phi^{\mathrm{true}} , s) \leq 0$ for any $s \in \mathcal{S}$.
    Moreover, $\phi^{*(N)} \geq \phi^{\mathrm{true}}$ (\Cref{prop:true_constraints_para_dominated}) implies $A^{c} = \{ \exists j , \, \phi_j^{*(N)} > \phi_j^{\mathrm{true}} + \delta_j \}$, and hence, as in \Cref{eq:subopt_stat_learning_theo_020-01},
    \[
        \mathbb{P} (A^{c})
        \leq \sum_{j=1}^J \mathbb{P} \left( \phi_j^{\sup} (\{ S \}) > \phi_j^{\mathrm{true}} + \delta_j \right)^N
        .
    \]

    Statement (2) of \Cref{theo:suboptimality_learning_prob}:
    on $A$, we have $\phi^{*(N)} \sim \phi^{\mathrm{true}}$ since $\Phi (\delta) / \sim = \{ [ \phi^{\mathrm{true}} ] \}$.
    Hence, by \Cref{eq:Phi_delta_maximum_020} and statement (2) of \Cref{prop:learning_theory_prob_suboptimality} with $\Phi = \{ \phi^{\mathrm{true}} \}$, for any $u \geq 0$,
    \begin{align*}
        & \mathbb{P} \left(
            \phi^{*(N)} \sim \phi^{\mathrm{true}}
            \text{ and }
            \mathbb{E} \ell^{\mathrm{sub},\lambda} ( \theta^{*(N)} , \phi^{*(N)} , S )
            \leq \varepsilon (u , C, N )
        \right)
        \\
        & \geq
        \mathbb{P} \left(
            A
            \text{ and }
            \mathbb{E} \ell^{\mathrm{sub},\lambda} ( \check{\theta}^{(N)} , \phi^{\mathrm{true}} , S )
            \leq \varepsilon (u , C, N )
        \right)
        \\
        & \geq
        1 - \mathbb{P} (A^{c})
        - \mathbb{P} \left(
            \mathbb{E} \ell^{\mathrm{sub},\lambda} ( \check{\theta}^{(N)} , \phi^{\mathrm{true}} , S )
            > \varepsilon (u , C, N )
        \right)
        \\
        & \geq 1 - \sum_{j=1}^J \mathbb{P} \left(
            \phi_j^{\sup} (\{ S \} ) >  \phi_j^{\mathrm{true}} + \delta_j
        \right)^N
        - 2 \zeta \left( \frac{C^2}{36} -2 \right) \exp \left( -\frac{C^2}{36} u^2  \right)
        .
    \end{align*}

    Statement (2) of \Cref{theo:suboptimality_learning_E}:
    by \Cref{eq:Phi_delta_maximum_020} and $\ell^{\mathrm{sub},\lambda} \geq 0$,
    \[
        \mathbb{E} \left[
            \mathbb{E}_S \ell^{\mathrm{sub},\lambda} ( \theta^{*(N)} , \phi^{*(N)} , S )
            \, \mathbf{1}_{A}
        \right]
        \leq
        \mathbb{E} \, \mathbb{E}_S \ell^{\mathrm{sub},\lambda} ( \check{\theta}^{(N)} , \phi^{\mathrm{true}} , S )
        .
    \]
    Applying \Cref{theo:learning_theory_expectation} with $\Phi = \{ \phi^{\mathrm{true}} \}$ to the right-hand side and dividing by $\mathbb{P} (A) \geq 1 - \sum_{j=1}^J \mathbb{P} ( \phi_j^{\sup} (\{ S \}) > \phi_j^{\mathrm{true}} + \delta_j )^N > 0$, we obtain statement (2) of \Cref{theo:suboptimality_learning_E} in the same manner as in its proof.

    \Cref{theo:subpoptimality_learning_prob_pl,theo:subpoptimality_learning_E_pl}:
    as in their proofs, the empirical loss of the output of \Cref{alg:Solve_IOP} is $0$ and $\phi^{*(N)} = \phi^{\sup} (\mathcal{S}_N)$.
    Hence the conclusions follow from statement (2) of \Cref{theo:suboptimality_learning_prob,theo:suboptimality_learning_E} shown above and the bound $\int_0^{\infty} \sqrt{\log N (\Theta , \| \bullet \|_2 ,\varepsilon )} d \varepsilon \leq 2 \sqrt{2 \log 2} \sqrt{D-1}$.
\end{proof}

\begin{remark}
    \label{rem:product_space_bound}
    If one applies \Cref{theo:learning_theory_expectation} to the product space $\Theta \times \Phi$, the generalization error can only be bounded as $O\left( \sqrt{(D + \dim \Phi -1) / N} \right)$.
    This bound is looser than that in \Cref{theo:subpoptimality_learning_E_pl_short}.
    Similarly, applying \Cref{theo:learning_theory_prob} to $\Theta \times \Phi$ yields a bound that is looser than \Cref{theo:subpoptimality_learning_prob_pl_short}.
\end{remark}

{
\renewcommand{\proofname}{Proof of \Cref{theo:subpoptimality_learning_prob_pl_short}}
\begin{proof}
    By taking $C = 8 \sqrt{2}$ in \Cref{theo:subpoptimality_learning_prob_pl} (which satisfies $C>6\sqrt{3}$), the conclusion follows (in this case, $2 \zeta \left( \frac{C^2}{36} - 2 \right) = 2 \zeta \left( \frac{14}{9} \right) \leq 5$ and $\frac{C^2}{36} = \frac{32}{9} \geq 3$).
\end{proof}
}

\section{Numerical experiment}

\subsection{Details of Implementation and Devices}
\label{sec:devices}

The fundamental libraries used in the experiment are OR-Tools v9.8 \citep{ortools}, NumPy 1.26.3 \citep{harris2020array}, and Python 3.9.0 \citep{Python3}.
They are distributed under the Apache License 2.0, the BSD 3-Clause License, and the PSF License, respectively.
Our computing environment is a machine with 192 Intel CPUs and 1.0TB CPU memory.

We describe the common settings of the optimization implementation.
The initial point of \Cref{alg:psgd_suboptimality_loss} is $\theta^1 = (1/D , \ldots , 1/D) + 10^{-3} (1 , \ldots , 1)$ (a shifted barycenter of the probability simplex); unless otherwise stated, the learning-rate constant is $\alpha = 1$ and the penalty coefficient is $\lambda = 0$ (after constraint learning, the observed solutions are feasible on the sample, so the penalty term is identically $0$).
The stopping criterion is that the $\ell_2$ norm of the subgradient $F (\theta^k)$ (\Cref{prop:SL_is_Lipschitz_convex}) is less than $10^{-6}$.
Since the features are integer-valued, this is equivalent to $F (\theta^k) = 0$, and by the convexity of the suboptimality loss, to the empirical loss being $0$, that is, to \Cref{eq:IOP_linear} being solved exactly on the sample (\Cref{prop:suboptimality_imitation}).
The MSEF (mean squared error with respect to the features) on the vertical axes of the figures is defined by
$\mathrm{MSEF} (\theta) := \frac{1}{N} \sum_{s \in \mathcal{S}_N} \| f (x^{*} (\theta , \phi^{\sup} , s) , s) - f (\hat{x}^{*} (s) , s) \|_2^2$,
and SL denotes the suboptimality loss $\mathbb{E} \ell^{\mathrm{sub}, 0}$.
The value reported at each iteration is the value at the iteration corresponding to the output (best-so-far) of \Cref{alg:psgd_suboptimality_loss}.
To display $0$ on logarithmic axes, the figures plot MSEF with an offset of $0.1$ and SL with an offset of $0.001$ (the vertical-axis labels read MSEF + 0.1 and SL + 0.001).
The curve for each setting is obtained by taking, at each iteration, the worst value over all episodes and shaping it to be monotonically nonincreasing.
The measured (wall) time is the total of constraint learning and objective learning.
In the experiments below, the weight set is $\Theta = \Delta^{D-1} + 10^{-3} (1, \ldots , 1)$, which does not satisfy $\sup_{\theta \in \Theta} \| \theta \|_2 \leq 1$ assumed in \Cref{sec:statistical_learning_theory}. However, since the optimal solutions of the forward problem are invariant under positive scaling of $\theta$ and the subgradient $F$ depends on $\theta$ only through them, running \Cref{alg:psgd_suboptimality_loss} on $\Theta$ with the learning-rate constant $\alpha$ is equivalent to running it on $c \Theta \subset \Delta^{D-1}$ with $c \alpha$, where $c := (1 + D \cdot 10^{-3})^{-1}$; the latter weight set satisfies the assumption, and the term of the suboptimality loss depending on $\theta$ is multiplied by $c$.%

\subsection{Single-Machine Scheduling with Release Dates and the Weighted Sum of Completion Times}
\label{sec:1-machine_appx}

\paragraph{Setting}
In the single-machine scheduling problem $1 \mid r_i \mid \sum \theta_i C_i$, we consider processing $D$ jobs on a single machine.
We assume that the machine can process at most one job at a time, and that once the processing of a job starts, it cannot be preempted.
Let $i=1 , \ldots , D$ index the jobs.
For job $i$, let $p_i$ denote the processing time, $\theta_i$ the importance weight, $r_i$ the release date (i.e., the earliest time at which the job can start processing), and let $\phi_{ik} \in \{0,1\}$.
The objective is to find an execution order (schedule) on the machine so as to minimize the weighted sum of completion times $C_i$.

Let the integer variable $b_i$ denote the start time of job $i$, and let $x_{ik}$ be a binary integer variable that equals $1$ if job $i$ precedes job $k$, and $0$ otherwise.
Define $M := \max_i r_i + \sum_i p_i$.
Then, the problem can be formulated as follows:
\begin{align*}
    \text{minimize}_{b,x} \,
    &
    \sum_{i=1}^D \theta_i ( b_i + p_i )
    \\
    \text{subject to} \,
    &
    b_i + p_i - M (1 - x_{ik}) \leq b_k , & \forall i \not = k, \\
    & x_{ik} + x_{ki} =1, \, x_{ik} \in \{ 0,1\} , & \forall i \not = k, \\
    & b_i \geq r_i, \quad b_i \leq 10^9 + 7, \quad b_i \in \mathbb{Z} & \forall i, \\
    & x_{ki} \leq \phi_{ik} & \forall i \not =k ,
\end{align*}
where $r_i$ are i.i.d.\ samples from the uniform distribution on $\{ 0 , 1 , \ldots , 9 \}$ (integers), and $p_i$ are i.i.d.\ samples from the uniform distribution on $\{ 1 , \ldots , 4 \}$ (integers).
Let $\mathcal{S}$ be the set of pairs $s=(p,r)$, from which we sampled $N$ states independently; we denote the resulting set by $\mathcal{S}_N$.
Let $\mathcal{X}(\phi,s)$ be the set of all $(b,x)$ satisfying the constraints above.
This problem is an instance of \Cref{exa:lattice_homomorphism_affine_check}.
We define
$
    \Theta = \Delta^{D-1} + 10^{-3} (1, \ldots , 1)
$.
We describe the generation of the observed data.
The true weight $\theta^{\mathrm{true}}$ is obtained by taking a standard Gaussian vector, normalizing it in $\ell_2$, taking the absolute value of each component, and adding $10^{-3} (1 , \ldots , 1)$ (since the optimal solutions of the suboptimality loss are invariant under positive scalar multiplication of the weights, the choice of normalization does not affect the results).
The true precedence constraint $\phi^{\mathrm{true}}$ allows all pairs ($\phi_{ik}^{\mathrm{true}} = 1$) except that precedence is enforced ($\phi_{i^{\dagger} k^{\dagger}}^{\mathrm{true}} = 0$) for a single pair $(i^{\dagger} , k^{\dagger})$ chosen uniformly at random.
For this choice, we used the random seeds $144 + t$ and $145 + t$ for episode number $t$.
The forward problem is solved with pywraplp of OR-Tools (the BOP solver with default settings), and the iteration limit of \Cref{alg:psgd_suboptimality_loss} is $K = 2000$. The other settings are as in \Cref{sec:devices}.%

Under this setting, we run \Cref{alg:Solve_IOP}.
Specifically, we first compute $\phi$ from the expert features $b^{(s)} = f(\hat{x}^*(s),s)$ according to
\begin{equation}
    \phi_{ik} = \begin{cases}
        0, & \text{ if } \forall s \in \mathcal{S}_N \, b_i^{(s)} \leq b_k^{(s)} ,\\
        1 , & \text{otherwise}.
    \end{cases}
\end{equation}
We then run \Cref{alg:psgd_suboptimality_loss}.

\paragraph{Experiment 1}

\begin{figure}
    \centering
    \begin{minipage}{0.48\linewidth}
            \includegraphics[width=\linewidth]{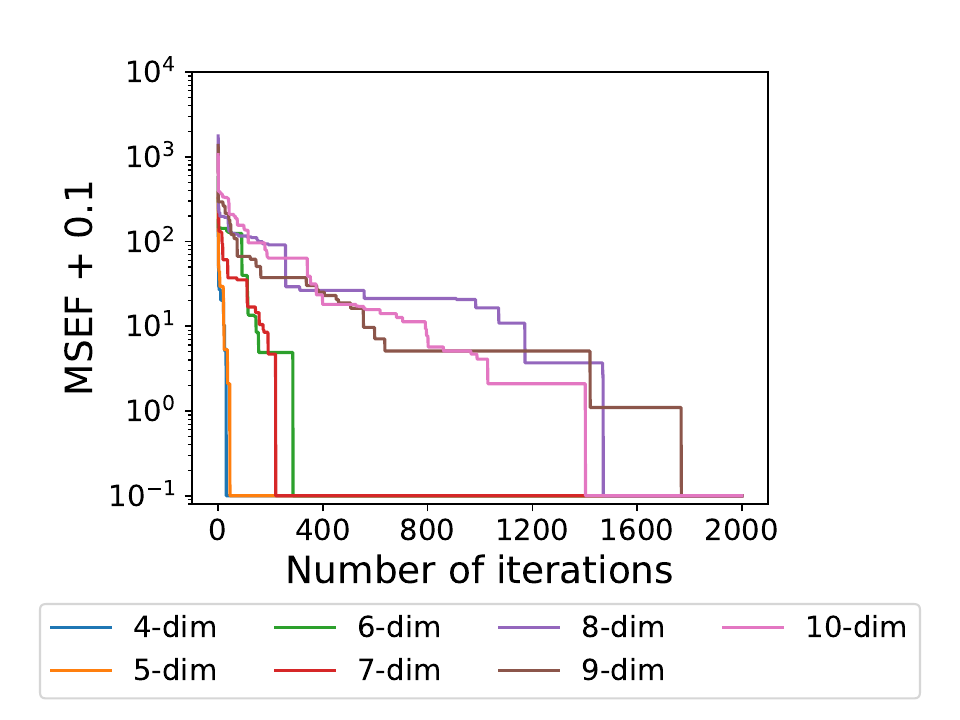}
    \caption{Worst-case MSEF over all episodes at each iteration of \Cref{alg:psgd_suboptimality_loss} (Experiment~1: $N = 5$, $D = 4, \ldots, 10$, $10$ episodes, SRSL). The iteration at which a curve reaches $0$ is the worst-case number of iterations required to solve \Cref{eq:IOP_linear}. The vertical axis reports the mean squared error with respect to the feature vector (MSEF).}
    \label{fig:1-machine_N_detail_MSEF}
    \end{minipage}
    \hfill
    \begin{minipage}{0.48\linewidth}
        \includegraphics[width=\linewidth]{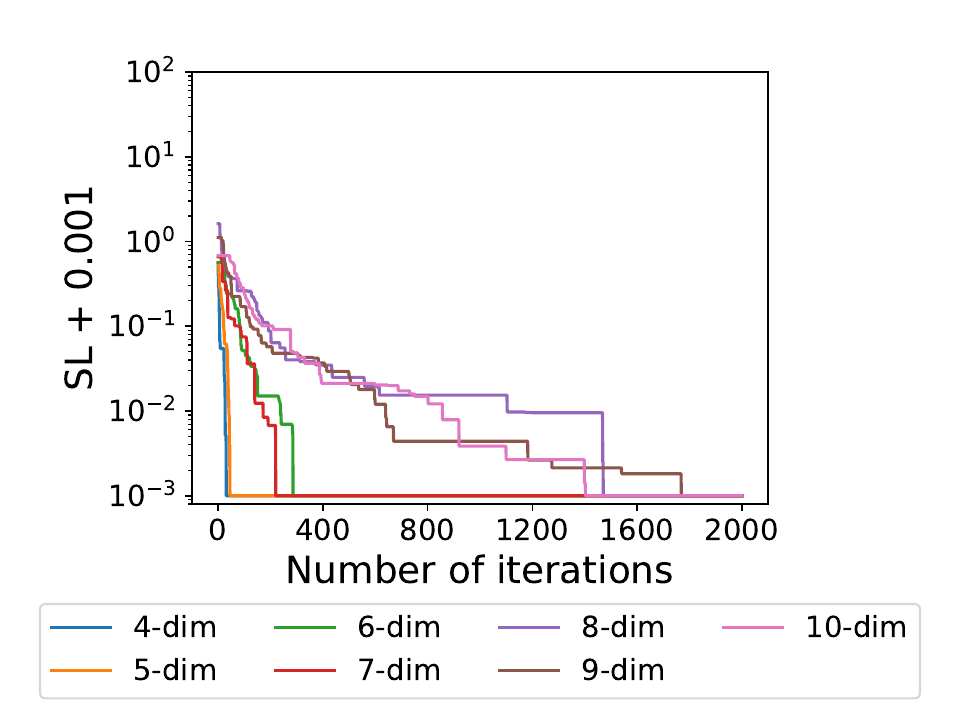}
    \caption{Worst-case suboptimality loss over all episodes at each iteration of \Cref{alg:psgd_suboptimality_loss} (Experiment~1: $N = 5$, $D = 4, \ldots, 10$, $10$ episodes, SRSL). The iteration at which a curve reaches $0$ is the worst-case number of iterations required to solve \Cref{eq:IOP_linear}. The vertical axis reports the suboptimality loss.}
    \label{fig:1-machine_N_detail_SL}
    \end{minipage}
\end{figure}

We implemented \Cref{alg:Solve_IOP} using \Cref{alg:psgd_suboptimality_loss} with the SRSL learning-rate schedule.
Using this implementation, we conducted $10$ episodes with $N=5$ for $D=4,5,6,7,8,9,10$, and report in \Cref{fig:1-machine_N_detail_MSEF} (and the suboptimality loss in \Cref{fig:1-machine_N_detail_SL}) the number of iterations in \Cref{alg:psgd_suboptimality_loss} required until \Cref{eq:IOP_linear} is solved.
As shown in \Cref{fig:1-machine_N_detail_MSEF}, we empirically verified that \Cref{eq:IOP_linear} can be solved for all values of $D$ considered.
We also summarize in \Cref{tab:1-machine_N5_detail} the wall-clock time required to solve \Cref{eq:IOP_linear}.

\begin{table}[ht]
    \centering
    \caption{The time required to solve \Cref{eq:IOP_linear} in Experiment~1.}
    \label{tab:1-machine_N5_detail}
    \begin{tabular}{llllll}
    \toprule
    $D$ & Decision variables & Constraints & Mean (s) & Max (s) & Median (s) \\
    \midrule
    4  & 16  & 40  & 0.125   & 0.391    & 0.089   \\
    5  & 25  & 65  & 0.313   & 0.737    & 0.322   \\
    6  & 36  & 96  & 2.245   & 6.492    & 1.627   \\
    7  & 49  & 133 & 3.314   & 11.830   & 2.296   \\
    8  & 64  & 176 & 62.761  & 237.787  & 51.446  \\
    9  & 81  & 225 & 194.870 & 1040.056 & 55.118  \\
    10 & 100 & 280 & 325.188 & 2244.227 & 99.000  \\
    \bottomrule
    \end{tabular}
\end{table}

\paragraph{Experiment 2}

\begin{figure}
    \centering
    \begin{minipage}{0.48\linewidth}
            \includegraphics[width=\linewidth]{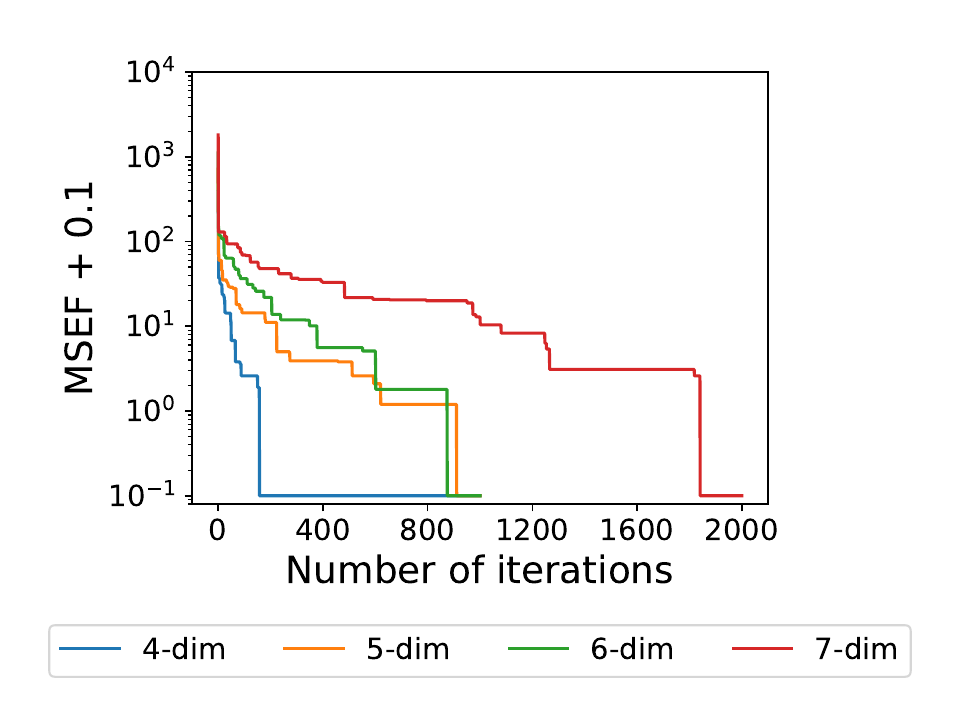}
    \caption{Worst-case MSEF over all episodes at each iteration of \Cref{alg:psgd_suboptimality_loss} (Experiment~2: $N = 10$, $D = 4, \ldots, 7$, $25$ episodes, SRSL). The iteration at which a curve reaches $0$ is the worst-case number of iterations required to solve \Cref{eq:IOP_linear}. The vertical axis reports the mean squared error with respect to the feature vector (MSEF).}
    \label{fig:1-machine_N10_detail_MSEF}
    \end{minipage}
    \hfill
    \begin{minipage}{0.48\linewidth}
        \includegraphics[width=\linewidth]{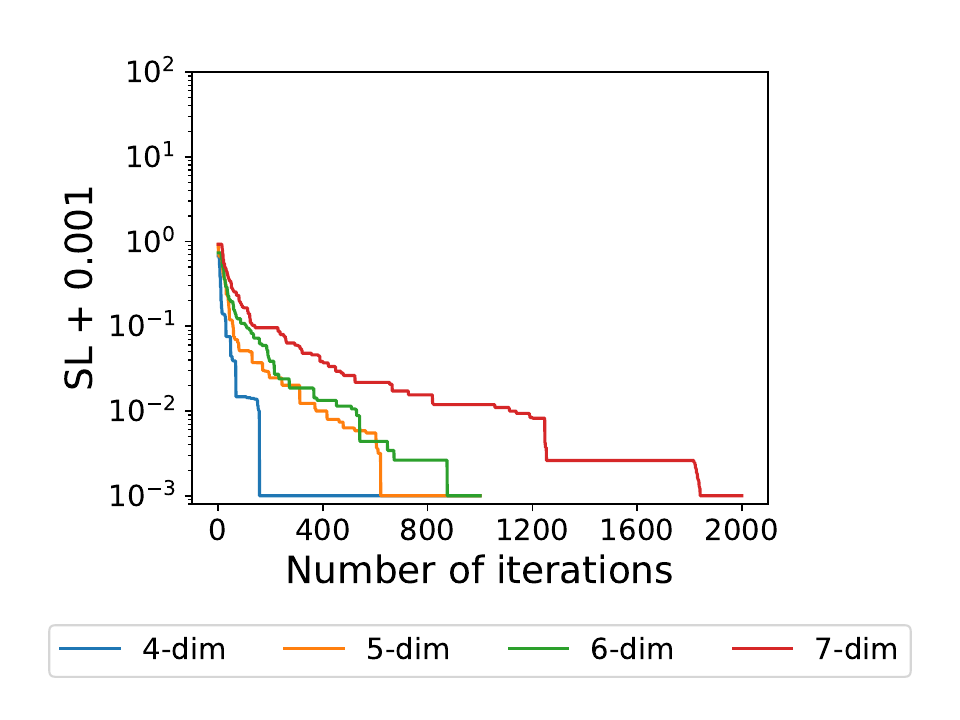}
    \caption{Worst-case suboptimality loss over all episodes at each iteration of \Cref{alg:psgd_suboptimality_loss} (Experiment~2: $N = 10$, $D = 4, \ldots, 7$, $25$ episodes, SRSL). The iteration at which a curve reaches $0$ is the worst-case number of iterations required to solve \Cref{eq:IOP_linear}. The vertical axis reports the suboptimality loss.}
    \label{fig:1-machine_N10_detail_SL}
    \end{minipage}
\end{figure}

We implemented \Cref{alg:Solve_IOP} using \Cref{alg:psgd_suboptimality_loss} with the SRSL learning-rate schedule.
Using this implementation, we conducted $25$ episodes with $N=10$ for $D=4,5,6,7$, and report in \Cref{fig:1-machine_N10_detail_MSEF} (and the suboptimality loss in \Cref{fig:1-machine_N10_detail_SL}) the number of iterations in \Cref{alg:psgd_suboptimality_loss} required until \Cref{eq:IOP_linear} is solved.
As shown in \Cref{fig:1-machine_N10_detail_MSEF}, we empirically verified that \Cref{eq:IOP_linear} can be solved for all values of $D$ considered.
We also summarize in \Cref{tab:1-machine_N10_detail} the wall-clock time required to solve \Cref{eq:IOP_linear}.

\begin{table}[ht]
    \centering
    \caption{The time required to solve \Cref{eq:IOP_linear} in Experiment~2.}
    \label{tab:1-machine_N10_detail}
    \begin{tabular}{l|rrrrr}
    $D$ & 4 & 5 & 6 & 7  \\
    \hline
    Decision variables & 16 & 25 & 36 & 49  \\
    Constraints & 40 & 65 & 96 & 133  \\
    Mean (s) & 1.04 & 6.16 & 9.24 & 63.00  \\
    Max (s) & 3.79 & 28.44 & 43.05 & 202.19  \\
    Median (s) & 0.33 & 3.62 & 6.42 & 44.81 \\
    Completion rate & 100\% & 100\% & 100\% & 100\% \\
    \end{tabular}
\end{table}

\subsection{Single-Machine Scheduling with Release Dates and the Weighted Sum of Tardiness}
\label{sec:tardiness_appx}

\begin{figure}[t]
    \vspace{-\intextsep}
    \centering
    \begin{minipage}[c]{0.48\textwidth}
        \centering
        \includegraphics[width=\linewidth]{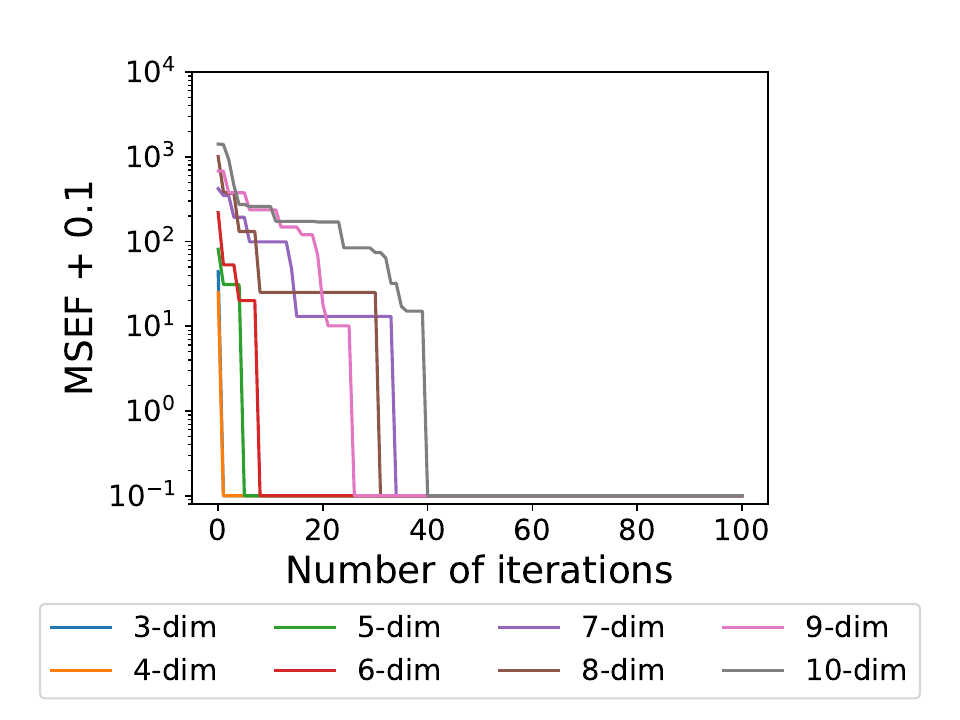}
        \vspace{-2\intextsep}
        \caption{Worst-case MSEF over all episodes at each iteration of \Cref{alg:psgd_suboptimality_loss} (\Cref{sec:tardiness}: $N = 1$, $D = 3, \ldots, 10$, $25$ episodes, SRSS). The iteration at which a curve reaches $0$ is the worst-case number of iterations required to solve \Cref{eq:IOP_linear}. The vertical axis reports the MSEF, plotted with an offset of $0.1$ to display $0$ on the logarithmic axis.}
        \label{fig:tardiness_SRSS_N1}
    \end{minipage}
    \hfill
    \begin{minipage}[c]{0.48\textwidth}
        \centering
        \includegraphics[width=\linewidth]{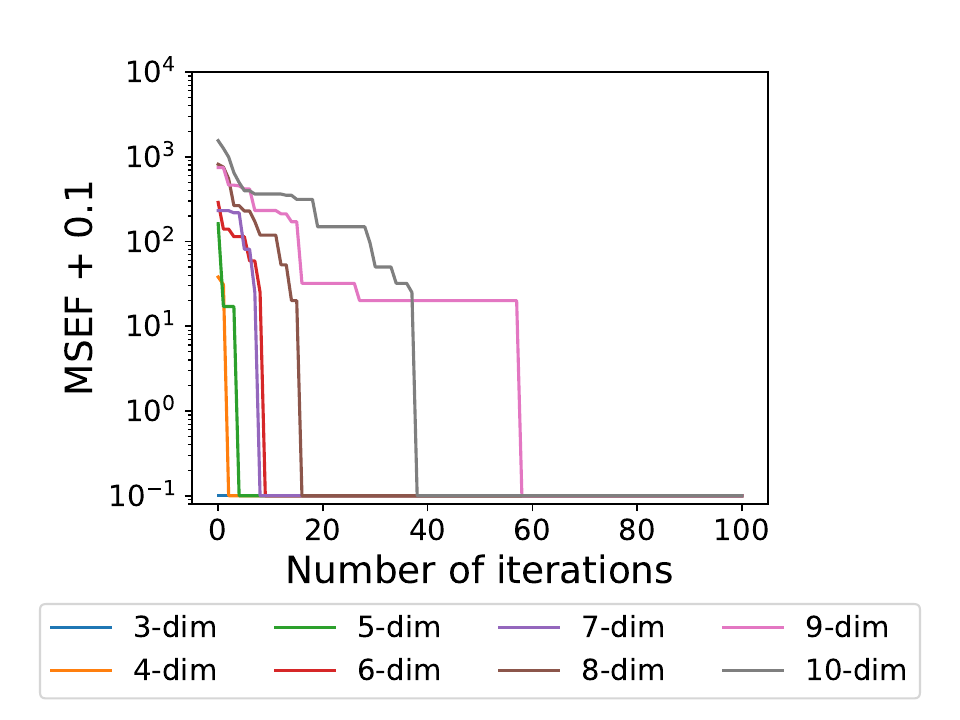}
        \vspace{-2\intextsep}
        \caption{Worst-case MSEF over all episodes at each iteration of \Cref{alg:psgd_suboptimality_loss} (\Cref{sec:tardiness}: $N = 1$, $D = 3, \ldots, 10$, $25$ episodes, SRSL). The iteration at which a curve reaches $0$ is the worst-case number of iterations required to solve \Cref{eq:IOP_linear}. The vertical axis reports the MSEF, plotted with an offset of $0.1$ to display $0$ on the logarithmic axis.}
        \label{fig:tardiness_SRSL_N1}
    \end{minipage}
    \vspace{-\intextsep}
\end{figure}

We consider the single-machine scheduling problem with release dates of minimizing the weighted sum of tardiness \citep[\S~3]{PostekZoccaAMPL2024}.

\paragraph{Setting}
Let $D$ denote the number of jobs, indexed by $i,j =1 , \ldots , D$.
We define the constants and decision variables as follows (\Cref{tab:tardiness_constants,tab:tardiness_variables}).
\begin{table}[ht]
    \begin{minipage}[t]{0.475\textwidth}
    \centering
    \caption{List of constants}\label{tab:tardiness_constants}
    \begin{tabular}{p{42pt}p{103pt}}
        \textbf{Symbol} & \textbf{Description} \\
        $\text{release}_i$ & time at which job $i$ becomes available \\
        $\text{duration}_i$ & processing time of job $i$ \\
        $\text{due}_i$ & due date of job $i$ \\
        $\theta_i$ & cost incurred per unit tardiness for job $i$ \\
        $M$ & big-$M$
    \end{tabular}
    \end{minipage}
    \hfill
    \begin{minipage}[t]{0.475\textwidth}
    \centering
    \caption{List of decision variables}\label{tab:tardiness_variables}
    \begin{tabular}{p{42pt}p{103pt}}
        \textbf{Symbol} & \textbf{Description} \\
        $\text{start}_i$ & start time of job $i$ \\
        $\text{finish}_i$ & completion time of job $i$ \\
        $y_i$ & tardiness of job $i$ \\
        $z_{ij}$ & $1$ if job $i$ precedes job $j$, and $0$ otherwise
    \end{tabular}
    \end{minipage}
\end{table}
The scheduling problem is formulated as follows:
{\allowdisplaybreaks
\begin{align}
    \text{minimize } & \sum_{i=1}^D \theta_i y_i
    \\
    \text{subject to } &
    \text{start}_i \geq \text{release}_i , & \forall i = 1, \ldots, D,
    \label{eq:min_pass_sum_st_1}
    \\
    & \text{finish}_i - \text{start}_i = \text{duration}_i , & \forall i = 1, \ldots, D,
    \label{eq:min_pass_sum_st_2}
    \\
    & y_i \geq \text{finish}_i - \text{due}_i & \forall i = 1, \ldots, D,
    \label{eq:min_pass_sum_st_3}
    \\
    & y_i \geq 0 , & \forall i = 1, \ldots, D,
    \label{eq:min_pass_sum_st_4}
    \\
    & z_{ij} \in \{0,1\}, \quad  z_{ij} + z_{ji} = 1 , & \forall i,j = 1, \ldots, D, \text{ s.t. } i \neq j ,
    \label{eq:min_pass_sum_st_5}
    \\
    & \text{finish}_i \leq \text{start}_j + M (1-z_{ij}) , & \forall i,j = 1, \ldots, D \text{ s.t. } i \neq j .
    \label{eq:min_pass_sum_st_6}
\end{align}
}
The variables $\text{start}_i$, $\text{finish}_i$, and $y_i$ take values in $\{ 0 , 1 , \ldots , H \}$, where $H := \sum_{i} \text{duration}_i + \max_i \text{release}_i + \max_i \text{due}_i$.
Here, \Cref{eq:min_pass_sum_st_1} enforces that each start time must be no earlier than the corresponding release time.
Constraint \Cref{eq:min_pass_sum_st_2} defines the completion time as the sum of the start time and the processing time.
Constraints \Cref{eq:min_pass_sum_st_3,eq:min_pass_sum_st_4} define the tardiness $y_i$ as the amount by which the completion time exceeds the due date, and enforce $y_i=0$ when the job is not tardy.
Constraint \Cref{eq:min_pass_sum_st_5} specifies the binary precedence variables.
Constraint \Cref{eq:min_pass_sum_st_6} imposes the precedence relationship: when $z_{ij}=1$, job $i$ precedes job $j$, and when $z_{ij}=0$, it does not.

We assume that $\text{due}_i$ can be written as $\text{due}_i = \text{release}_i + \text{duration}_i + \phi_i^1$ for some constant $\phi_i^1 \in \mathbb{Z}_{\geq 0}$.
Let $\mathcal{S}$ denote a subset of the set of release-time vectors $\{ \text{release}_i\}_i$.
Substituting \Cref{eq:min_pass_sum_st_2} and $\text{due}_i = \text{release}_i + \text{duration}_i + \phi_i^1$ into \Cref{eq:min_pass_sum_st_3}, the latter is equivalent to $y_i \geq \text{start}_i - \text{release}_i - \phi_i^1$, which depends only on $\phi_i^1$. We also write the equality constraint \Cref{eq:min_pass_sum_st_2} as the two inequalities $\text{finish}_i - \text{start}_i \leq \text{duration}_i^{+}$ and $\text{finish}_i - \text{start}_i \geq \text{duration}_i^{-}$ with two parameters $\text{duration}_i^{\pm}$. Let $\Phi = \{ 0 \}^{J^0} \times \mathbb{Z}_{\geq 0}^D \times \mathbb{Z}_{\geq 0}^D \times \mathbb{Z}_{\geq 0}^D$ be the parameter space of the constraint parameters $\phi = ( 0 , \phi^1 , (\text{duration}_i^{+})_i , (\text{duration}_i^{-})_i )$. Here, $J^0$ denotes the number of constraints that contain no parameter, and the corresponding components are set to $\{ 0 \}$ as in \Cref{exa:lattice_homomorphism_affine}. Then each constraint depends monotonically on a single coordinate of $\phi$; after reversing the signs of $\phi^1$ and $\text{duration}^{+}$ as in \Cref{exa:lattice_homomorphism_affine_check}, the constraint map is a lattice homomorphism by \Cref{theo:characterization_lattice_homomorphism}. The true parameter satisfies $\text{duration}_i^{+\mathrm{true}} = \text{duration}_i^{-\mathrm{true}} = \text{duration}_i$.
For $\phi\in \Phi$ and $s \in \mathcal{S}$, let $\mathcal{X} (\phi , s)$ denote the set of decision variables satisfying constraints \Cref{eq:min_pass_sum_st_1,eq:min_pass_sum_st_2,eq:min_pass_sum_st_3,eq:min_pass_sum_st_4,eq:min_pass_sum_st_5,eq:min_pass_sum_st_6} (with \Cref{eq:min_pass_sum_st_2,eq:min_pass_sum_st_3} rewritten as above).
We set
$
    \Theta = \Delta^{D-1} + 10^{-3} (1, \ldots , 1)
$.

We describe the sampling procedure for the state set $\mathcal{S}$.
Each $\text{release}_i$ is sampled uniformly at random from the integers between $0$ and $5$.
Each $\text{duration}_i$ is sampled uniformly at random from the integers between $1$ and $4$.
Each $\phi_i^{1,\mathrm{true}}$ is sampled uniformly at random from the integers between $0$ and $8$.

We next describe the generation of the weight vector $\theta^{\mathrm{true}}$.
We first sample each component $\theta_i^{\mathrm{true}}$ uniformly at random from the integers between $1$ and $3$, and then rescale the resulting vector so that $\theta^{\mathrm{true}} \in \Delta^{D-1}$.
Since $\theta^{\mathrm{true}}$ is drawn from a finite set, the guarantees for almost every $\theta^{\mathrm{true}}$ (\Cref{theo:alg:Solve_IOP_with_min_suboptloss_is_solve_IOP_LP,prop:iteration_count_solve_IOP}) do not directly apply to this setting, and \Cref{assu:WIRL-uniqueness} is not guaranteed to hold.
The implementation, stopping criterion, and metrics of the optimization are as in \Cref{sec:devices}.
The forward problem is solved with pywraplp of OR-Tools (the SCIP backend with default settings), the big-M constant is $M = 10^5$, and the objective coefficients are multiplied by $10^4$ and rounded to integers when passed to the solver.
Since any schedule can be shifted to the left so that all jobs finish by $\max_i \text{release}_i + \sum_i \text{duration}_i$ without increasing any tardiness, an instance-dependent value such as $M = \max_i \text{release}_i + \sum_i \text{duration}_i$ ($\leq 5 + 4D$ under the sampling below) would be more appropriate in view of the LP relaxation and numerical stability.
The iteration limit of \Cref{alg:psgd_suboptimality_loss} is $K = 1000$, and the learning-rate constant is $\alpha = 1$ for SRSL and $\alpha = 1 / \| F (\theta^{1}) \|_2$ (the reciprocal of the norm of the subgradient at the first iteration) for SRSS.
The parameter $\text{duration}^{\sup}$ is computed from the first observation, and $\phi^{1,\sup}$ as the maximum over all observations; since the durations do not depend on the state, $\text{duration}_i^{+\sup} = \max_{s} ( \text{finish}_i^{(s)} - \text{start}_i^{(s)} )$ and $\text{duration}_i^{-\sup} = \min_{s} ( \text{finish}_i^{(s)} - \text{start}_i^{(s)} )$ both coincide with the value $\text{duration}_i^{\sup}$ of the first observation.

Under this setting, we run \Cref{alg:Solve_IOP}.
For each state $s \in \mathcal{S}_N$, let
\[
\hat{x}^{*} (s) = (y^{(s)},z^{(s)}, (\text{start}_i^{(s)})_i ,(\text{finish}_i^{(s)})_i )= x^* (\theta^{\mathrm{true}},\phi^{\mathrm{true}},s) \in \mathcal{X} (\phi^{\mathrm{true}},s).
\]
We then compute
\begin{align}
    \text{duration}_i^{\sup} = \text{finish}_i^{(s)} - \text{start}_i^{(s)},
    \quad
    \phi_i^{1,\sup} = \max_{s \in \mathcal{S}_N} \left( \text{start}_i^{(s)} - \text{release}_i^{(s)} - y_i^{(s)} , 0   \right)
\end{align}
We subsequently run \Cref{alg:psgd_suboptimality_loss}.

\begin{figure}[ht]
    \centering
    \begin{minipage}{0.48\linewidth}
        \centering
            \includegraphics[width=\linewidth]{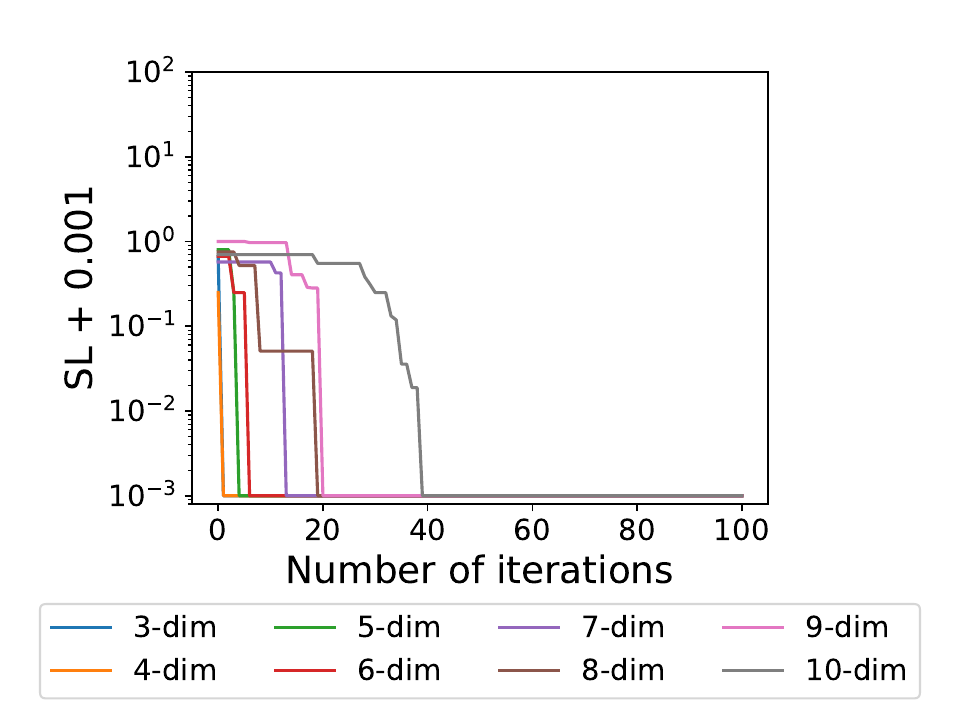}
    \caption{Worst-case suboptimality loss over all episodes at each iteration of \Cref{alg:psgd_suboptimality_loss} (\Cref{sec:tardiness}: $N = 1$, $D = 3, \ldots, 10$, $25$ episodes, SRSS). The iteration at which a curve reaches $0$ is the worst-case number of iterations required to solve \Cref{eq:IOP_linear}. The vertical axis reports the suboptimality loss.}
    \label{fig:tardiness_SRSS_N1_detail_SL}
    \end{minipage}
    \hfill
    \begin{minipage}{0.48\linewidth}
        \centering
        \includegraphics[width=\linewidth]{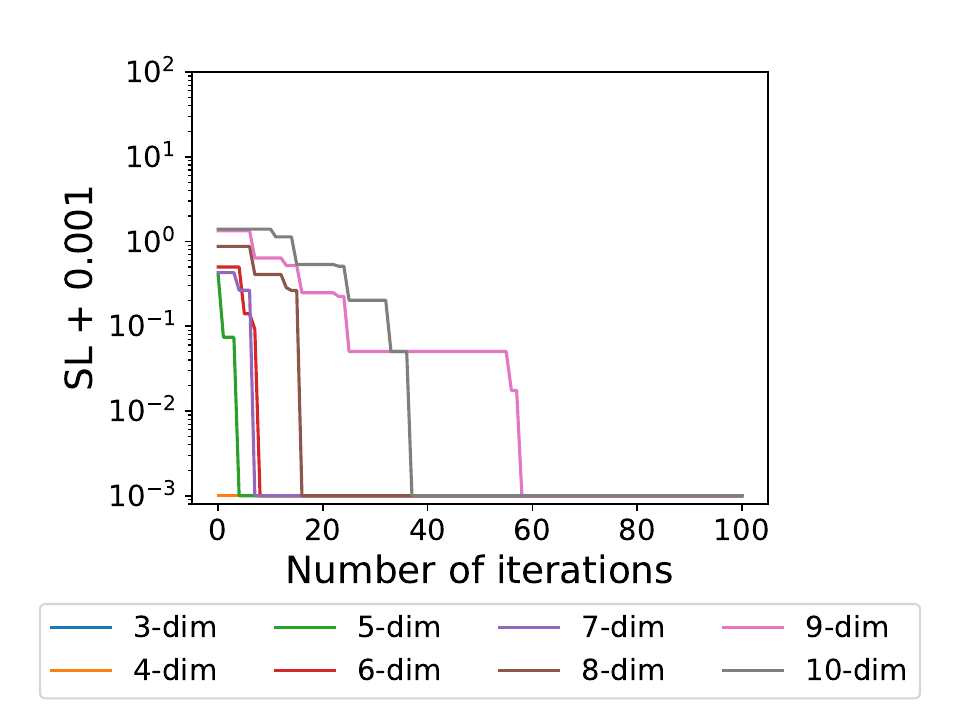}
    \caption{Worst-case suboptimality loss over all episodes at each iteration of \Cref{alg:psgd_suboptimality_loss} (\Cref{sec:tardiness}: $N = 1$, $D = 3, \ldots, 10$, $25$ episodes, SRSL). The iteration at which a curve reaches $0$ is the worst-case number of iterations required to solve \Cref{eq:IOP_linear}. The vertical axis reports the suboptimality loss.}
    \label{fig:tardiness_SRSL_N1_detail_SL}
    \end{minipage}
\end{figure}

\paragraph{Experiments}
Using this implementation, we set $N=1$ and varied $D$ from $3$ to $10$, conducting $25$ episodes for each value of $D$.
For \Cref{alg:Solve_IOP} equipped with \Cref{alg:psgd_suboptimality_loss} using the SRSS learning-rate schedule, we report the wall-clock time required to solve \Cref{eq:IOP_linear} in \Cref{tab:tardiness_SRSS_N1_detail}, and the evolution of the suboptimality loss until \Cref{eq:IOP_linear} is solved in \Cref{fig:tardiness_SRSS_N1_detail_SL}.
For \Cref{alg:Solve_IOP} equipped with \Cref{alg:psgd_suboptimality_loss} using the SRSL learning-rate schedule, we report the wall-clock time required to solve \Cref{eq:IOP_linear} in \Cref{tab:tardiness_SRSL_N1_detail}, and the evolution of the suboptimality loss until \Cref{eq:IOP_linear} is solved in \Cref{fig:tardiness_SRSL_N1_detail_SL}.

\begin{table}[ht]
\centering
    \begin{minipage}{0.48\linewidth}
        \centering
        \caption{Wall-clock time required to solve \Cref{eq:IOP_linear} with the SRSS learning-rate schedule}
        \label{tab:tardiness_SRSS_N1_detail}
        \begin{tabular}{llll}
        \toprule
        $D$ & Mean (s) & Max (s) & Median (s) \\
        \midrule
        3  & 0.011  & 0.014   & 0.011  \\
        4  & 0.018  & 0.069   & 0.014  \\
        5  & 0.025  & 0.094   & 0.021  \\
        6  & 0.098  & 0.430   & 0.047  \\
        7  & 0.294  & 2.293   & 0.108  \\
        8  & 1.738  & 10.860  & 0.996  \\
        9  & 10.830 & 88.349  & 3.365  \\
        10 & 60.700 & 369.823 & 31.816 \\
        \bottomrule
        \end{tabular}
    \end{minipage}
    \hfill
    \begin{minipage}{0.48\linewidth}
        \centering
        \caption{Wall-clock time required to solve \Cref{eq:IOP_linear} with the SRSL learning-rate schedule}
        \label{tab:tardiness_SRSL_N1_detail}
        \begin{tabular}{llll}
        \toprule
        $D$ & Mean (s) & Max (s) & Median (s) \\
        \midrule
        3  & 0.012  & 0.020   & 0.012  \\
        4  & 0.014  & 0.029   & 0.013  \\
        5  & 0.030  & 0.101   & 0.021  \\
        6  & 0.076  & 0.441   & 0.036  \\
        7  & 0.438  & 3.496   & 0.324  \\
        8  & 1.516  & 5.001   & 1.240  \\
        9  & 5.541  & 28.730  & 3.391  \\
        10 & 43.111 & 207.598 & 31.381 \\
        \bottomrule
        \end{tabular}
    \end{minipage}
\end{table}

\subsection{Details of the Verification of the Generalization Error (\Cref{sec:heldout})}
\label{sec:heldout_appx}

\paragraph{Details of the results in the main text}
First, since $D = 5$ is fixed, this experiment examines only the $N$-dependence of the bound: in terms of the log-log slope with respect to $N$, the loss on the held-out states decays at least as fast as the theoretical rate $O(1/\sqrt{N})$ (the log-log slope of the mean is about $-1.9$, and about $-1.5$ for the conditional mean corresponding to \Cref{theo:subpoptimality_learning_E_pl_short}).
This suggests that, in our noiseless and realizable setting, the uniform-convergence-type rate $O (1/\sqrt{N})$ is conservative and faster rates can hold (see (i) in the discussion below).
Second, we observe the decomposition structure of the generalization error (\Cref{theo:subpoptimality_learning_prob_pl_short}): when $N$ is small, the penalty term (the constraint side) is dominant (about $99\%$ of the total at $N = 1$), but it decays rapidly, roughly as $N^{-2.2}$, and at $N = 32$ the max term (the objective side) takes over as the dominant part.
That is, constraint learning is statistically easy, and for large $N$ the learning of the objective function governs the generalization error.
Third, the exact recovery rate of the constraint parameters increases monotonically from $0\%$ at $N = 1$ to $92\%$ at $N = 32$ (\Cref{fig:phi_recovery}).

\paragraph{Procedure}
The true parameters, the state distribution, and the learning settings (SRSL, $\alpha = 1$, $K = 1000$, and the stopping criterion) are the same as in \Cref{sec:tardiness_appx}.
For each episode, the random seed is fixed to $20260724 + 1000 N + (\text{episode number})$.
In each episode, we sample $N$ training states and $100$ held-out states independently, learn $(\theta^{*(N)} , \phi^{*(N)})$ using only the training states, and then, for each held-out state $s$, compute
\begin{align*}
    \text{max term}
    & := \mathrm{ReLU} \left(
        \max_{x \in \mathcal{X} (\phi^{*(N)} , s)} \theta^{*(N) \top} f (x , s)
        - \theta^{*(N) \top} f (\hat{x}^{*} (s) , s)
    \right) ,
    \\
    \text{penalty term}
    & := \sum_{i=1}^D \mathrm{ReLU} \left( \text{finish}_i - \text{due}_i^{\sup} - y_i \right)
    \\
    & \quad + \sum_{i=1}^D \left| ( \text{finish}_i - \text{start}_i ) - \text{duration}_i^{\sup} \right|
\end{align*}
(where $(\text{start} , \text{finish} , y)$ are the components of the expert solution at the held-out state $s$, and $\text{due}^{\sup}$ is the due date given by the learned thresholds; the penalty of the equality constraints uses the absolute value because violations are two-sided), and set $\ell^{\mathrm{sub}, 1} = \text{max term} + \text{penalty term}$.
The MSEF on the held-out states is computed analogously.

\paragraph{Results}
\Cref{tab:heldout_results} summarizes the results (unless otherwise noted, the values are averages over $25$ episodes).
\begin{table}[ht]
    \centering
    \caption{Results of the verification of the generalization error ($D = 5$, $25$ episodes for each $N$, $100$ held-out states, $\lambda = 1$). Standard errors in parentheses.}
    \label{tab:heldout_results}
    \begin{tabular}{p{0.5cm}p{1.8cm}p{1.5cm}p{1.5cm}p{1.2cm}p{2.3cm}p{1.4cm}}
        \toprule
        $N$ & Exact attainment (train) & $\phi$ recovery & Penalty term & Max term & $\ell^{\mathrm{sub},1}$ & MSEF (median) \\
        \midrule
        1  & $100\%$ & $0\%$  & $3.69$   & $0.047$  & $3.74\ (0.46)$     & $7.53$ \\
        2  & $100\%$ & $8\%$  & $1.87$   & $0.058$  & $1.92\ (0.30)$     & $6.95$ \\
        4  & $96\%$  & $52\%$ & $0.354$  & $0.054$  & $0.409\ (0.093)$   & $4.16$ \\
        8  & $100\%$ & $68\%$ & $0.085$  & $0.028$  & $0.113\ (0.029)$   & $0.99$ \\
        16 & $84\%$  & $84\%$ & $0.038$  & $0.012$  & $0.051\ (0.019)$   & $1.29$ \\
        32 & $68\%$  & $92\%$ & $0.0012$ & $0.0030$ & $0.0042\ (0.0010)$ & $0.36$ \\
        \bottomrule
    \end{tabular}
\end{table}

\begin{figure}[ht]
    \centering
    \begin{minipage}[c]{0.48\textwidth}
        \centering
        \includegraphics[width=\linewidth]{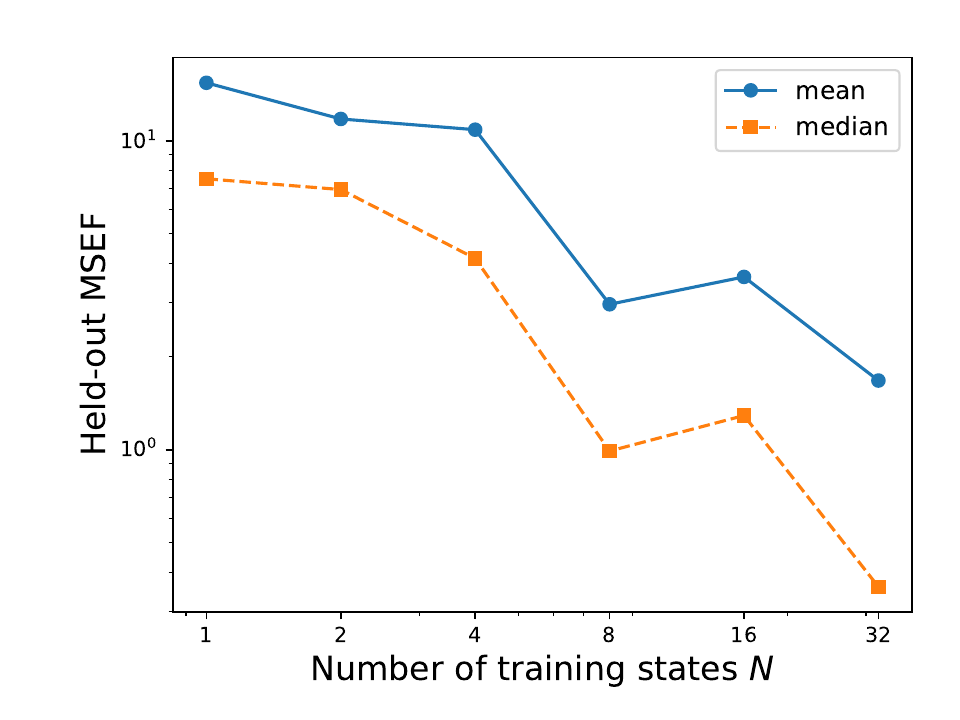}
        \caption{Mean and median of the MSEF on held-out states ($D = 5$, log-log scale).}
        \label{fig:heldout_MSEF}
    \end{minipage}
\end{figure}

\paragraph{Discussion}
(i) Convergence rates: the log-log slope of the loss $\ell^{\mathrm{sub},1}$ on the held-out states is about $-1.9$ for the mean (about $-2.2$ for the median).
The conditional mean $\hat{\mathbb{E}} \left[ \text{max term} \,\middle|\, \phi^{*(N)} = \phi^{\mathrm{true}} \right]$, which directly corresponds to \Cref{theo:subpoptimality_learning_E_pl_short}, also decays from $0.079$ ($N = 4$) to $0.0033$ ($N = 32$) with a slope of about $-1.5$.
In terms of the slope with respect to $N$ (with $D = 5$ fixed), the observed decay is no slower than the theoretical rate, which is thus conservative in this setting.
Whether faster rates hold in our noiseless and realizable setting is left for future work.
(ii) Decomposition structure: the penalty term decays roughly as $N^{-2.2}$ and the max term roughly as $N^{-0.8}$ (about $N^{-1.4}$ for $N \geq 4$), and the dominant term switches at $N = 32$ (\Cref{tab:heldout_results}).
(iii) Behavior of the non-recovery probability: the fraction $q_N := 1 - \hat{\mathbb{P}} ( \phi^{*(N)} = \phi^{\mathrm{true}} )$ of episodes in which the constraint parameter is not recovered satisfies $N q_N \approx 2.6$, which is approximately constant, for $N \geq 8$; this is consistent with $q_N \propto 1/N$. However, these estimates are based on only $8$, $4$, and $2$ unrecovered episodes out of $25$ for $N = 8$, $16$, and $32$, respectively, so this behavior should be regarded as a rough tendency.
The probability term $\sum_j \mathbb{P} (\cdot)^N$ in \Cref{theo:subpoptimality_learning_prob_pl_short} gives an exponential decay for a \emph{fixed} problem instance, whereas in this experiment $(\theta^{\mathrm{true}} , \phi^{\mathrm{true}})$ is resampled in each episode, so a mixture of instance-dependent exponential rates can yield a polynomial decay on average. The two are not in contradiction.
(iv) Gap with the agreement of the solutions: the decay of the MSEF on the held-out states is slow (slope about $-0.7$, \Cref{fig:heldout_MSEF}), and at $N = 32$ the median MSEF is still $0.36$, whereas the mean of $\ell^{\mathrm{sub},1}$ is $0.0042$ (\Cref{tab:heldout_results}).
This quantifies the fact that a small suboptimality loss does not immediately imply agreement of the solutions (actions).
(v) Training side: the exact attainment rate on the training data decreases as $N$ grows ($68\%$ at $N = 32$, hitting the iteration limit $K = 1000$, with an average computation time of about $334$ seconds).
This is consistent with the $N$-dependence of the iteration-count bound (\Cref{prop:iteration_count_solve_IOP}).
On the other hand, even in the episodes where exact attainment was not reached during training, the loss on the held-out states remains small (about $0.006$ on average over the $8$ unattained episodes at $N = 32$), so generalization remains good even under early termination by the iteration limit.
(vi) Caveat: the slopes above are point estimates from a log-log least-squares fit over $6$ points, and the loss distributions are heavy-tailed to the right (e.g., the non-monotonicity of the MSEF at $N = 16$). The conclusions are confirmed with both the mean and the median.%

\end{document}